\documentclass{elsarticle}
\usepackage{graphicx}
\usepackage{subfigure}
\usepackage{color}
\usepackage{tikz}
\usepackage{newlfont}
\usepackage{multirow}
\usepackage{array}
\usepackage{longtable} %Allows for long tables that stretch across many pages
\usepackage{ifthen}
\usepackage{alltt}
\usepackage{float}
\usepackage{amsmath}
\usepackage{xcolor}

\usepackage{enumerate}
 \usepackage[text={6.5in,8.5in},centering]{geometry} 
\newcolumntype{C}[1]{>{\centering\let\newline\\\arraybackslash\hspace{0pt}}m{#1}}
\newcolumntype{L}[1]{>{\raggedright\let\newline\\\arraybackslash\hspace{0pt}}m{#1}}
\newcolumntype{R}[1]{>{\raggedleft\let\newline\\\arraybackslash\hspace{0pt}}m{#1}}
\newcommand{\ba}{\begin{array} }
\newcommand{\ea}{\end{array} }
\newcommand{\bae}{\begin{eqnarray}}
\newcommand{\eae}{\end{eqnarray}}
\newcommand{\bea}{\begin{eqnarray*}}
\newcommand{\eea}{\end{eqnarray*}}
\newcommand{\be}{\begin{equation}}
\newcommand{\ee}{\end{equation}}

\newcommand{\lx}{\left(}
\newcommand{\rx}{\right)}
\newcommand{\lz}{\left[ }
\newcommand{\rz}{\right] }

\newcommand{\pr}{{\bf Proof}~~}

\usepackage{amssymb}
\usepackage{amsthm}
\usepackage{ mathrsfs }
\usepackage{amsmath}
\usepackage[mathcal]{eucal}

\usepackage{natbib}
\usepackage{graphicx}
\usepackage{color}
\usepackage[colorlinks]{hyperref}
\usepackage{lscape}
\usepackage{graphicx}
\usepackage{epstopdf}
\newtheorem{theorem}{\hskip\parindent\bf Theorem}[section]

\newtheorem{definition}{\hskip\parindent\bf Definition}[section]

\usepackage{multicol}
\usepackage{tabularx}
\usepackage{ctable}
\usepackage{booktabs}
\usepackage{setspace}
\usepackage{bm}
\begin{document}
 \markboth{A two-patch prey-predator interaction mode}{Collaborations}
\title{A two-patch prey-predator model with dispersal in predators driven by the strength of predation}
\author{Yun Kang\footnote{Sciences and Mathematics Faculty, College of Letters and Sciences, Arizona State University, Mesa, AZ 85212, USA ({\tt yun.kang@asu.edu})}, Sourav Kumar Sasmal\footnote{Agricultural and Ecological Research Unit, Indian Statistical Institute, 203, B. T. Road, Kolkata, 700108, INDIA.}, and Komi Messan\footnote{Simon A. Levin Mathematical and Computational Modeling Sciences Center, Arizona State University, Tempe, AZ. }}
%\institute{Yun Kang \at Sciences and Mathematics Faculty, School of Letters and Sciences, Arizona State University, Mesa, AZ 85212, USA. \\
      %    \email{yun.kang@asu.edu} \and Sourav Kumar Sasmal \at Agricultural and Ecological Research Unit, Indian Statistical Institute, 203, B. T. Road, Kolkata, 700108, INDIA. \\\email{sourav.sasmal@gmail.com} \and  Komi Messan \at Simon A. Levin Mathematical and Computational Modeling Sciences Center, Arizona State University, Tempe, AZ. \\ \email{komimessan@gmail.com}}
%\maketitle

\begin{abstract}
Foraging movements of predator play an important role in population dynamics of prey-predator interactions, which have been considered as mechanisms that contribute to spatial self-organization of prey and predator. In nature, there are many examples of prey-predator interactions where prey is immobile while predator disperses between patches non-randomly through different factors such as stimuli following the encounter of a prey. In this work, we formulate a Rosenzweig-MacArthur prey-predator two patch model with mobility only in predator and the assumption that predators move towards patches with more concentrated prey-predator interactions. We provide completed local and global analysis of our model. Our analytical results combined with bifurcation diagrams suggest that: (1) dispersal may stabilize or destabilize the coupled system; (2) dispersal may generate multiple interior equilibria that lead to rich bistable dynamics or may destroy interior equilibria that lead to the extinction of predator in one patch or both patches; (3) Under certain conditions, the large dispersal can promote the permanence of the system. In addition, we compare the dynamics of our model to the classic two patch model to obtain a better understanding how different dispersal strategies may have different impacts on the dynamics and spatial patterns.
%\keywords{Rosenzweig-MacArthur Prey-Predator Model \and Self-organization Effects \and Dispersal \and Persistence \and Non-random Foraging Movements}
%\subclass{34D20 \and 34D23 \and 34D45 \and 92D25 \and 92D40}
\end{abstract}
\bigskip
\maketitle
\section{Introduction}\label{sec_introduction}
Spatial heterogeneity, dispersal patterns, and biotic interactions influence the distribution of species within a landscape \citep{bolker1999spatial, levin1974dispersion, rees1996quantifying, seabloom2005spatial, shahak2008photoselective}. Spatial self-organization results from local interactions between organisms and the environment and emerges at patch-scales \citep{rietkerk2008regular, van2008experimental}. For example, limited dispersal ability and its related dispersal patterns \citep{miller2002old} is considered to be one of the key factors that promotes the development of self-organized spatial patterns \citep{aarssen1985biotic, kefi2007local, sole2006self, soro1999species}.\\

In nature, especially for ecological communities of insects, dispersal of a predator is usually driven by its non-random foraging behavior which can often response to prey-contact stimuli \citep{hassell1978foraging}, including spatial variation in prey density \citep{kareiva1987swarms} and different type of signals arising directly from prey \citep{Waage1977}. For instances, bloodsucking insects respond to the carbon dioxide output and the visual signals of a moving animal, which in tsetse flies (\emph{Glossina spp.}) lead to the formation of a ``following swarm" associated with herds of grazing ungulates \citep{Ford1971, Gatehouse1972}. Most mosquitoes were attracted over a larger distance by the odor of the host \citep{gillies1970range, gillies1972range, gillies1974range}. The wood-wasp, \emph{Sirex noctilio}, is attracted by  the concentration of the scent \citep{gillies1974range, madden1977physiological}. Social ants excite ``pheromone trails" to encourage other individuals to visit the same food source \citep{carroll1973ecology}. Plant-feeding insects commonly detect food items by gustatory signals \citep{Schoonhoven1974, Schoonhoven1976, Schoonhoven1977}. These non-random foraging behaviors driven by prey-mediated patch attractants, prey attractants themselves, and arrestant stimuli following the encounter of a prey, can lead to predation rates that are greater in regions where prey are more abundant (i.e., density-dependent predation), thus regulate population dynamics of both prey and predator.\\

Recent experimental work on population dynamics of immobile Aphids and Coccinellids by \citet{kummel2013aphids} show that the foraging movements of predator Coccinellids are combinations of passive diffusion, conspecific attraction, and retention on plants with high aphid numbers which is highly dependent on the strength of prey-predator interaction. Their study also demonstrates that predation by coccinellids was responsible for self-organization of aphid colonies. Many ecological systems exhibit similar foraging movements of predator. For example, Japanese beetles are attracted to feeding induced plant volatiles and congregate where feeding is taking place \citep{loughrin1996role}. Motivated by these field studies, we propose a two-patch prey-predator model incorporating foraging movements of predator driven by the strength of prey-predator interaction, to explore how this non-random dispersal behavior of predator affect population dynamics of prey and predator. \\

Dispersal of predator plays an important role in regulating, stabilizing, or destabilizing population dynamics of both prey and predator. There are fair amount literature on mathematical models of prey-predator interactions in patchy environments. For example, see work of \cite{levin1974dispersion, hassell1991spatial, bascompte1994spatially, hassell1995appropriate, Hastings1983, Hastings1993, ruxton1996density, tilman1997spatial, hanski1999metapopulation, hanski1997metapopulationbiology, Rohani1999dispersal, Kang2011a, Kang_C2012a} and also see \cite{kareiva1990population} for literature review.  Many studies examine how the interactions between patches affect the synchronicity of the oscillations in each patch, e.g. see the work of \cite{hauzy2010density, jansen1995regulation}, and how interactions may stabilize or destabilize the dynamics. For instances, \cite{jansen1994theoretical,jansen1995regulation} studied a model with two patches, each with the well known prey-predator Rosenzweig-McArthur dynamics, linked by density independent  dispersal (i.e., dispersal is driven by the difference of species' population densities in two patches). His study showed that this type of spatial predator-prey interactions might exhibit self-organization capable of producing stabilizing heterogeneities in prey distribution, and spatial populations can be regulated through the interplay of local dynamics and migration. \\

However, due to the intricacies that arise in density-dependent dispersal models, there are relatively limited work on models with non-random foraging behavior of predator or non-linear dispersal behavior \citep{kareiva1987swarms} but see the two patch model with predator attraction to prey, e.g. \cite{huang2001predator}), or predator attraction to conspecific, e.g. \cite{murdoch1992aggregation}, or only predators migrate who are attracted to regions with concentrated food resources, see the work of \cite{hassell1974aggregation, chesson1986aggregation}. \citet{kareiva1987swarms} proposed a non random foraging PDE model through a mechanistic approach to demonstrate that area-restricted search does yield predator aggregation, and explore the the consequences of area-restricted search for predator-prey dynamics. In addition, they provided many supporting ecological examples (e.g. Coccinellids, blackbirds, etc.) that abide by their theory.  \cite{huang2001predator} studied a two-patch predator-prey Rosenzweig-MacArthur  model with nonlinear density-dependent migration in the predator. The migration term of the predator is derived by extending the Holling time budget argument to migration. Their study showed that the extension of the Holling time budget argument to movement has essential effects on the dynamics. By extending the model of \cite{huang2001predator},  \citet{ghosh2011two} formulated a similar two patch prey-predator model with density-independent migration in prey and  density-dependent migration in the predator. Their study shows that  several foraging parameters such as handling time, dispersal rate can have important consequences in stability of prey-predator system.  \citet{cressman2013two} investigated the population-dispersal dynamics for predator-prey interactions in a two patch environment with assumptions that both predators and their prey are mobile and their dispersal between patches is directed to the higher fitness patch. They proved that such dispersal, irrespectively of its speed, cannot destabilize a locally stable predator- prey population equilibrium that corresponds to no movement at all. \\

In this paper, we formulates a new version of Rosenzweig-MacArthur two patch prey predator model with mobility only in predator. Our model is distinct from others as we assume that the non-random foraging movements of predator is driven by the strength of prey-predator interactions, i.e., predators move towards patches with more concentrated prey and predator.
Our model can apply to many insects systems such as Aphids and Coccinellids, Japanese beetles and their host plants, etc. The main focus of our study of such prey-predator interactions in heterogeneous environments is to explore the following ecological questions:
\begin{enumerate}
\item How does our proposed nonlinear density-dependent dispersal of predator stabilize or destabilize the system?
\item How does dispersal of predator affect the extinction and persistence of prey and predator in both patches?
\item How may dispersal promote the coexistence of prey and predator when predator goes extinct in the single patch?
\item What are potential spatial patterns of prey and predator?
\item How are the effects of our proposed nonlinear density-dependent dispersal of predator on population dynamics different from the effects of traditional density-independent dispersal?
\end{enumerate}

The rest of the paper is organized as follows:  In Section \eqref{sec_model_derivations}, we derive our two patch prey-predator model  and provide its basic dynamics properties. In Section \eqref{sec_Math_Analysis}, we perform completed local and global dynamics of our model, and derive sufficient conditions that lead to the persistence and extinction of predator as well as permanence of the model. In Section \eqref{sec_effects_of_dispersal}, we perform bifurcation simulations to explore the dynamical patterns and compare the dynamics of our model to the traditional model studied by \citet{jansen2001dynamics}. In Section \eqref{sec_discussion}, we conclude our study and discuss the potential future study. The detailed proofs of our analytical results are provided in the last section.

%%%%%%%%%%%%%%%%%%%%%%%%%%%%%%%%%%%%%%%%%%%%%%%%%%%%%%%%%%%%%%%%%%%%%%%%%%%%%%
%%%%%%%%%%%%%%%%%%%%%%%%%%%%%%%%%%%%%%%%%%%%%%%%%%%%%%%%%%%%%%%%%%%%%%%%%%%%%%
%%%%%%%%%%%%%%%%%%%%%%%%%%%%%%%%%%%%%%%%%%%%%%%%%%%%%%%%%%%%%%%%%%%%%%%%%%%%%%

\section{Model derivations}\label{sec_model_derivations}
Let $x_i, y_i$ be the population of prey, predator at Patch $i$, respectively. We consider the following two-patch prey-predator interaction model after rescaling (see similar rescaling approaches by \cite{Liuyuanyuan2010})\eqref{2patch}
\begin{equation}
\begin{aligned}
\frac{dx_1}{dt} & = x_1\left(1 - \frac{x_1}{K_1}\right) - \frac{a_1 x_1 y_1}{1 + x_1} \\
\frac{dy_1}{dt} & =  \frac{a_1 x_1 y_1}{1 + x_1} -d_1 y_1+ \rho_1\left(  \underbrace{ \frac{a_1 x_1 y_1}{1 + x_1}}_\text{attraction strength to Patch 1} y_2- \underbrace{\frac{a_2 x_2 y_2}{1 + x_2}}_\text{attraction strength to Patch 2} y_1\right)\\
\frac{dx_2}{dt} & = r x_2\left(1 - \frac{x_2}{K_2}\right) - \frac{a_2 x_2 y_2}{1 + x_2} \\
\frac{dy_2}{dt} & =  \frac{a_2 x_2 y_2}{1 + x_2} -d_2 y_2+\rho_2\left(   \underbrace{\frac{a_2 x_2 y_2}{1 + x_2}}_\text{attraction strength to Patch 2} y_1-  \underbrace{\frac{a_1 x_1 y_1}{1 + x_1}}_\text{attraction strength to Patch 1} y_2\right)
 \label{2patch}
\end{aligned}
\end{equation} where $K_i$ is the relative carrying capacity of prey in the absence of predation; $a_i$ is the relative predation rate at Patch 1; $d_i$ is the relative death rate of predator at Patch $i$; $\rho_i$ is the relative dispersal rate of predator at Patch $i$; and $r$ is the relative maximum growth rate of prey at Patch 2.  All parameters are nonnegative.
The ecological assumptions of Model \eqref{2patch} can be stated as follows:
\begin{enumerate}
\item In the absence of dispersal, Model \eqref{2patch} is reduced to the following uncoupled Rosenzweig-MacArthur prey-predator single patch models
\begin{equation}
\begin{aligned}
\frac{dx_i}{dt} & =r_i x_i\left(1 - \frac{x_i}{K_i}\right) - \frac{a_i x_i y_i}{1 + x_i} \\
\frac{dy_i}{dt} & =  \frac{a_i x_i y_i}{1 + x_i} -d_i y_i \label{Onepatch}
\end{aligned}
\end{equation} where $r_1=1$ and $r_2=r$ and its ecological assumptions \cite{rosenzweig1963graphical} can be stated as follows:
\begin{enumerate}
\item In the absence of predation, population of prey $x_i$ follows the logistic growth model.
\item Predator is specialist (i.e., predator $y_i$ goes extinct without prey $x_i$) and the functional response between prey and predator follows Hollying Type II functional response.
\end{enumerate}
\item There is no dispersal in prey species. This assumption fits in many prey-predator (or plant-insects) interactions in ecosystems such as Aphid and Ladybugs, Japanese beetles and its feeding plants, etc.
\item The dispersal of predator from Patch $i$ to Patch $j$ is driven by prey-predation interaction strength in Patch $j$ termed as \emph{attraction strength}. In Model \eqref{2patch}, we assume that predator in Patch $i$ disperse to Patch $j$ is determined by the predation term $\frac{a_j x_j y_j}{1 + x_j}$ in Patch $j$, thus the dispersal term of predator from Patch $i$ to Patch $j$ is described by $\rho_i\frac{a_j x_j y_j}{1 + x_j} y_i$ which gives the net dispersal of predator at Patch $i$ as
$$\rho_i\left(\frac{a_j x_j y_j}{1 + x_j} y_i-\frac{a_i x_i y_i}{1 + x_i} y_j\right)$$ by assuming that the dispersal constant $\rho_i$ of Patch $i$ is the same for predator arriving Patch $i$ from other patches as predator leaving Patch $i$ to other patches. This assumption is motivated by the fact that dispersal of a predator is usually driven by its non-random foraging behavior which can often response to prey-contact stimuli \citep{hassell1978foraging} which has been supported in many field studies including the recent work by \citet{kummel2013aphids}.

\end{enumerate}
The state space of Model \eqref{2patch} is $\mathbb R^4_+$. Let $\mu_i=\frac{d_i}{a_i-d_i},\nu_1=\frac{(K_1-\mu_1)(1+\mu_1)}{a_1K_1}, \nu_2=\frac{r(K_2-\mu_2)(1+\mu_2)}{a_2K_2}$. We have the following theorem regarding the dynamics properties of Model \eqref{2patch}:
%%%%%%%%%%%%%%%%%%%%%%%%%%%%%%%%%%%%%%%%%%%%%%%%%%%%%%%%
\begin{theorem}\label{th1:dp} Assume all parameters are nonnegative and $r, a_i, K_i, d_i, i=1,2$ are strictly positive. Model \eqref{2patch} is positively invariant and bounded in $\mathbb R^4_+$ with
$\limsup_{t\rightarrow\infty} x_i(t)\leq K_i$ for both $i=1,2.$
In addition, it has the following properties:
\begin{enumerate}
\item If there is no dispersal in predator, i.e., $\rho_1=\rho_2=0$, then Model \eqref{2patch} is reduced to Model \eqref{Onepatch} whose dynamics can be classified in the following three cases:
\begin{enumerate}
\item Model \eqref{Onepatch} always has the extinction equilibrium $(0,0)$ which is a saddle.
\item If $\mu_i>K_i$ or $\mu_i<0$, then the boundary equilibrium $(K_i,0)$ is globally asymptotically stable.
\item If $\frac{K_i-1}{2}<\mu_i<K_i$, then the boundary equilibrium $(K_i,0)$ is a saddle while the interior equilibrium $(\mu_i, \nu_i)$ is globally asymptotically stable.
\item If $0<\mu_i<\frac{K_i-1}{2}$, then the boundary equilibrium $(K_i,0)$ is a saddle; the interior equilibrium $(\mu_i, \nu_i)$ is a source, and the system has a unique limit cycle which is globally asymptotically stable. In addition, the Hopf bifurcation occurs at $\mu_i=\frac{K_i-1}{2}$.
\end{enumerate}
\item The sets $\{(x_1,y_1,x_2,y_2)\in\mathbb R^4_+:x_i=0 \}$ and $\{(x_1,y_1,x_2,y_2)\in\mathbb R^4_+: y_i=0 \}$ are invariant for both $i=1,2$. If $x_j=0$, Model \eqref{2patch} is reduced to the single patch model Model \eqref{Onepatch}.
If $y_j=0$, Model \eqref{2patch} is reduced to the following two uncoupled models:
\begin{equation}
\begin{aligned}
\frac{dx_i}{dt} & =r_i x_i\left(1 - \frac{x_i}{K_i}\right) - \frac{a_i x_i y_i}{1 + x_i} \\
\frac{dy_i}{dt} & =  \frac{a_i x_i y_i}{1 + x_i} -d_i y_i \\
\frac{dx_j}{dt} & =r_j x_j\left(1 - \frac{x_j}{K_j}\right) \label{yj0}
\end{aligned}
\end{equation} where $\lim_{t\rightarrow\infty}x_j(t)=K_j$ and the dynamics of $x_i, y_i$ is the same as Model \eqref{Onepatch}.
\end{enumerate}
\end{theorem}

\noindent\textbf{Notes:} Theorem \ref{th1:dp} provides a foundation on our further study of local stability and global dynamics of Model \eqref{2patch}. In addition, Item 2 of Theorem \ref{th1:dp} implies that Model \eqref{2patch} has the same the invariant sets $x_i=0$ and $y_i=0$ for both $i=1,2$ as the single patch models \eqref{Onepatch}. In addition, the results of the single patch models \eqref{Onepatch} indicate that prey is always persist while predator $i$ is persist if $0<\mu_i<K_i$ hold.

%%%%%%%%%%%%%%%%%%%%%%%%%%%%%%%%%%%%%%%%%%%%%%%%%%%%%%%%%%%%%%%
%%%%%%%%%%%%%%%%%%%%%%%%%%%%%%%%%%%%%%%%%%%%%%%%%%%%%%%%%%%%%%%
%%%%%%%%%%%%%%%%%%%%%%%%%%%%%%%%%%%%%%%%%%%%%%%%%%%%%%%%%%%%%%%

\section{Mathematical analysis}\label{sec_Math_Analysis}
Now we start with the boundary equilibria of Model \eqref{2patch}. Recall that 
$$\mu_i=\frac{d_i}{a_i-d_i},\nu_1=\frac{(K_1-\mu_1)(1+\mu_1)}{a_1K_1}, \nu_2=\frac{r(K_2-\mu_2)(1+\mu_2)}{a_2K_2}.$$ We define the following notations for all possible boundary equilibria of Model \eqref{2patch}:
$$\begin{array}{cccc}
E_{0000}=(0,0,0,0), &E_{K_1000}=(K_1,0,0,0),&E_{\mu_1\nu_100}=(\mu_1,\nu_1,0,0),&E_{K_10\mu_2\nu_2}=(K_1,0,\mu_2,\nu_2)\\
E_{K_10K_20}=(K_1,0,K_2,0), &E_{00K_20}=(0,0,K_2,0),&E_{00\mu_2\nu_2}=(0,0,\mu_2,\nu_2),&E_{\mu_1\nu_1K_20}=(\mu_1,\nu_1,K_2,0)
\end{array}.$$
The following theorem provides sufficient conditions on the existence and stability of these boundary equilibria:
%%%%%%%%%%%%%%%%%%%%%%%%%%%%%%%%%%%%%%%%%%%%%%%%%%%%%%%%%%%%%%%
\begin{theorem}\label{th2:be}[Boundary equilibria of Model \eqref{2patch}]Model \eqref{2patch} always has the following four boundary equilibria
$$E_{0000},\,E_{K_1000},\,E_{00K_20},\,E_{K_10K_20}$$ where the first three ones are saddles while $E_{K_10K_20}$ is locally asymptotically stable if $\mu_i>K_i$ and it is a saddle if $\left(\mu_1-K_1\right)\left(\mu_2-K_2\right)<0$ or $\mu_i<K_i, i=1,2$. Let $i,j=1,2,\,i\neq j$, and
$$E_{11}^b=E_{\mu_1\nu_100},\,E_{12}^b=E_{\mu_1\nu_1K_20},\,E_{21}^b=E_{00\mu_2\nu_2} \mbox{ and }E_{22}^b=E_{K_10\mu_2\nu_2}.$$
Then if $0<\mu_i<K_i$, then Model \eqref{2patch} has additional two boundary equilibria $E^b_{i1}$ and $E_{i2}^b$ where $E_{i1}^b$ is always a saddle. The boundary equilibrium $E_{i2}^b$
is locally asymptotically stable if $\frac{K_i-1}{2}<\mu_i<K_i$ and one of the following conditions holds:
 \begin{enumerate}
 \item[sa:]$a_j\leq d_i, K_j<\mu_j$.
 \item[sb:] $0<K_j<\min\Big\{\mu_j,\frac{d_i}{a_j-d_i}\Big\}.$
 \item[sc:] $0<\frac{d_i}{a_j-d_i}<K_j<\mu_j \mbox{ and }\rho_j<\frac{d_j-K_j(a_j-d_j)}{\nu_i\left[K_j(a_j-d_i)-d_i\right]}.$
 \item[sd:] $0<\mu_j<K_j<\frac{d_i}{a_j-d_i} \mbox{ and }\rho_j>\frac{K_j(a_j-d_j)-d_j}{\nu_i\left[d_i-K_j(a_j-d_i)\right]}.$
 \end{enumerate}
 And $E_{i2}^b$ is a saddle if $0<\mu_i<\frac{K_i-1}{2}$ or one of the following conditions holds:
 \begin{enumerate}
 \item[ua:] $K_j>\max\Big\{\mu_j,\frac{d_i}{a_j-d_i}\Big\}.$
  \item[ub:] $0<\frac{d_i}{a_j-d_i}<K_j<\mu_j\mbox{ and }\rho_j>\frac{d_j-K_j(a_j-d_j)}{\nu_i\left[K_j(a_j-d_i)-d_i\right]}.$
 \item[uc:] $\mu_j<K_j<\frac{d_i}{a_j-d_i} \mbox{ and }\rho_j<\frac{K_j(a_j-d_j)-d_j}{\nu_i\left[d_i-K_j(a_j-d_i)\right]}.$
 \end{enumerate}
In addition, if $0<\mu_i<K_i$ for both $i=1$ and $i=2$, the boundary equilibria $E_{12}^b$ and $E_{22}^b$ exist but they cannot be locally asymptotically stable at the same time while if $r_i=1, a_i=a,d_i=d, K_i=d$ for both $i=1,2$, the boundary equilibria $E_{12}^b$ and $E_{22}^b$ can not be locally asymptotically stable at all if they exist.\\\vspace{10pt}
\end{theorem}

\noindent\textbf{Notes:} Theorem \ref{th2:be} implies the following points regarding the effects of dispersal in predators:
\begin{enumerate}
\item Dispersal has no effects on the local stability of the boundary equilibrium $E_{K_10K_20}$.
\item Large dispersal of predator in its own patch may have stabilizing effects from the results of Item sd: In the absence of dispersal, the dynamics of Patch $j$ is unstable at $(K_j,0)$ since $0<\mu_j<K_j$. However, in the presence of dispersal, large values of $\rho_j$ can lead to the local stability of the boundary equilibrium $E^b_{i2}$ where $i,j=1,2$ and $i\neq j$, under conditions of $\mu_j<K_j<\frac{d_i}{a_j-d_i}$.
\item Large dispersal of predator in its own patch may have destabilizing effects from the results of Item ub: In the absence of dispersal, the dynamics of Patch $j$ is local stable at $(K_j,0)$ since $K_j<\mu_j$. However, in the presence of dispersal, large values of $\rho_j$ can drive the boundary equilibrium $E^b_{j2}$ being unstable, under conditions of $0<\frac{d_i}{a_j-d_1}<K_j<\mu_j$.
\item Under conditions of $\mu_i<K_i$, the boundary equilibria $E_{12}^b$ and $E_{22}^b$ can not be asymptotically stable at the same time.
\end{enumerate}

\subsection{Global dynamics}
In this subsection, we focus on the extinction and persistence dynamics of prey and predator of Model \eqref{2patch}. First we show the following theorem regarding the boundary equilibrium $E_{K_10K_20}$:
%%%%%%%%%%%%%%%%%%%%%%%%%%%%%%%%%%%%%%%%%%%%%%%%%%%%%%%%%%%%%%%

\begin{theorem}\label{th3:gb}
Model \eqref{2patch} has global stability at $E_{K_10K_20}$ if $\mu_i>K_i, i=1,2$.
\end{theorem}

\noindent\textbf{Notes:} Theorem \ref{th3:gb} implies that the dispersal of predators does not effect the global stability of the boundary equilibrium $E_{K_10K_20}$.\\
%%%%%%%%%%%%%%%%%%%%%%%%%%%%%%%%%%%%%%%%%%%%%%%%%%%%%%%%%%%%%%%

To proceed the statement and proof of our results on persistence, we provide the definition of \emph{persistence} and \emph{permanence} as follows:
\begin{definition}[Persistence of single species]We say species $z$ is persistent in $\mathbb R^4_+$ for Model \eqref{2patch} if there exists constants $0<b<B$, such that for any initial condition with $z(0)>0$, the following inequality holds
$$b\leq\liminf_{\tau\rightarrow\infty} z(\tau)\leq \limsup_{\tau\rightarrow\infty} x(\tau)\leq B.$$
\end{definition} where $z$ can be $x_i,y_i, i=1,2$ for Model \eqref{2patch}.
\begin{definition}[Permanence of a system]We say Model \eqref{2patch} is permanent in $\mathbb R^4_+$ if there exists constants $0<b<B$, such that for any initial condition taken in $\mathbb R^4_+$ with $x_1(0)y_1(0)x_2(0)y_2(0)>0$, the following inequality holds
$$b\leq\liminf_{\tau\rightarrow\infty}\min\{ x_1(\tau), y_1(\tau), x_2(\tau),y_2(\tau)\}\leq \limsup_{\tau\rightarrow\infty} \max\{x_1(\tau), y_1(\tau), x_2(\tau),y_2(\tau)\}\leq B.$$
\end{definition}The permanence of Model \eqref{2patch} indicates that all species in the system are persistence.
%%%%%%%%%%%%%%%%%%%%%%%%%%%%%%%%%%%%%%%%%%%%%%%%%%%%%%%%%%%%%%%

\begin{theorem}\label{th4:p}[Persistence of prey and predator]
Prey $x_i, i=1,2$ of Model \eqref{2patch} are always persistent for all $r>0$. Predator $y_j$ is persistent if one of the following inequalities hold
 \begin{enumerate}
\item $\mu_j<K_j, \mu_i>K_i$. Or
\item $\frac{K_i-1}{2}<\mu_i<K_i \mbox{ and } K_j>\max\Big\{\mu_j,\frac{d_i}{a_j-d_i}\Big\}$. Or
\item $\frac{K_i-1}{2}<\mu_i<K_i,\,\,\mu_j<K_j<\frac{d_i}{a_j-d_i} \mbox{ and }\rho_j<\frac{K_j(a_j-d_j)-d_j}{\nu_i\left[d_i-K_j(a_j-d_i)\right]}$. % where $i=1,j=2$ or $i=2,j=1$.
\end{enumerate}
\end{theorem}
\noindent\textbf{Notes:} Theorem \ref{th4:p} indicates that the dispersal of predators does not affect the persistence of preys, while small dispersal of predator $j$, under condition of $\frac{K_i-1}{2}<\mu_i<K_i,\,\,\mu_j<K_j<\frac{d_i}{a_j-d_i}$, can keep the persistence of predator $j$. This is consistent with the results of Item uc in Theorem \ref{th2:be}.\\
%%%%%%%%%%%%%%%%%%%%%%%%%%%%%%%%%%%%%%%%%%%%%%%%%%%%%%%%%%%%%%%

\begin{theorem}\label{th5:perm}[Permanence of the two patch dispersal model]
 Model \eqref{2patch} is permanent if one of the following inequalities hold
 \begin{enumerate}
 \item $\frac{K_j-1}{2}<\mu_j<K_j, \,0<\frac{d_j}{a_i-d_j}<K_i<\mu_i \mbox{ and } \rho_i>\frac{d_i-K_i(a_i-d_i)}{\nu_j\left[K_i(a_i-d_j)-d_j\right]}$ where $i=1,j=2$ or $i=2,j=1$. Or
 \item $\frac{K_j-1}{2}<\mu_j<K_j,\, \mu_i>\frac{K_i-1}{2} \mbox{ and } K_i>\max\Big\{\mu_i,\frac{d_j}{a_i-d_j}\Big\}$ for both $i=1,j=2$ and $i=2,j=1$. Or
\item $\frac{K_i-1}{2}<\mu_i<K_i,\,K_i>\max\Big\{\mu_i,\frac{d_j}{a_i-d_j}\Big\},\,\frac{K_j-1}{2}<\mu_j<K_j<\frac{d_i}{a_j-d_i} \mbox{ and }\rho_j<\frac{K_j(a_j-d_j)-d_j}{\nu_i\left[d_i-K_j(a_j-d_i)\right]}$ where $i=1,j=2$ or $i=2,j=1$.
\end{enumerate}
\end{theorem}
\noindent\textbf{Notes:}  According to Theorem \ref{th4:p}, we can conclude that Model \eqref{2patch} is permanent whenever both predators are persistent. Theorem \ref{th5:perm} provides such sufficient conditions that can guarantee the coexistence of bother predator for the two patch model \eqref{2patch}, thus provides sufficient conditions of its permanence. Item 1 of this theorem implies that if predator $j$ is persistent, and large dispersal of predator $i$ can promote its persistence, thus, promote the permanence since in the absence of dispersals in predator, predator $i$ goes extinct due to $\mu_i>K_i$. This is consistent with the results of Item ub of Theorem \ref{th2:be} that large dispersal of predator $i$ can have destabilize effects on the boundary equilibria $E^b_{j2}$. \\
%%%%%%%%%%%%%%%%%%%%%%%%%%%%%%%%%%%%%%%%%%%%%%%%%%%%%%%%%%%%%%%
\subsection{Interior equilibrium and stability}
Let $p_i(x)=\frac{a_i x}{1+x}$ and $q_i(x)=\frac{r_i(K_i-x)(1+x)}{a_iK_i}$, then we have
$$\begin{array}{lcl}
\frac{dx_i}{dt}&=&r_ix_i\left(1-\frac{x_i}{K_i}\right)-\frac{a_i x_iy_i}{(1+x_i)}=\frac{a_i x_i}{1+x_i}\left[\frac{r_i(K_i-x_i)(1+x_i)}{a_iK_i}-y_i\right]=p_i(x_i)\left[q_i(x_i)-y_i\right]\\\\
\frac{dy_i}{dt}&=&y_i\left[\frac{a_i x_i}{1+x_i}-d_i+\rho_i y_j\left(\frac{a_i x_i}{1+x_i}-\frac{a_j x_j}{1+x_j}\right)\right]=y_i\left[p_i(x_i)-d_i+\rho_i y_j\left(p_i(x_i)-p_j(x_j)\right)\right]\end{array}$$

If $(x_1^*,y_1^*,x_2^*, y_2^*)$ is an interior equilibrium of Model \eqref{2patch}, then it satisfies the following equations:
\bae\label{interior-eq1}
\begin{array}{lcl}
q_i(x_i)-y_i=0\Leftrightarrow q_i(x_i)=y_i,\,\\
p_i(x_i)-d_i+\rho_i y_j\left(p_i(x_i)-p_j(x_j)\right)=0 \Leftrightarrow p_i(x_i)=\frac{a_i x_i}{1+x_i}=\frac{\rho_i y_jp_j(x_j)+d_i}{1+\rho_i y_j}=\frac{\rho_i q_j(x_j)p_j(x_j)+d_i}{1+\rho_i q_j(x_j)}
\end{array}
\eae which gives:
\bae\label{interior-eq2}
x_i&=&\frac{\rho_i q_j(x_j)p_j(x_j)+d_i}{a_i(1+\rho_i q_j(x_j))-(\rho_i q_j(x_j)p_j(x_j)+d_i)}=\frac{\rho_i q_j(x_j)p_j(x_j)+d_i}{\rho_i q_j(x_j)\left[a_i-p_j(x_j)\right]+a_i-d_i}
\eae
Since $\limsup_{t\rightarrow\infty} x_i(t)\leq K_i$ for both $i=1,2$ and $y_i=q_i(x_i)$, therefore, positive solutions of $x_i\in (0,K_i)$ for \eqref{interior-eq2} determine interior equilibrium of Model \eqref{2patch}. By substituting the explicit forms of $p_i, q_i$ into \eqref{interior-eq2}, we obtain the following null clines:
\bae\label{interior-eq3}
\begin{array}{lcl}
x_1&=&\frac{a_2\left[r_2\rho_1x_2\left(K_2-x_2\right)+K_2d_1\right]}{r_2\rho_1x_2\left(K_2a_1-K_2a_2-a_1\right)-r_2\rho_1x_2^2(a_1-a_2)+K_2(a_1r_2\rho_1+a_1a_2-a_2d_1)}=\frac{f_t(x_2)}{f_b(x_2)}=F(x_2)\\\\
x_2&=&\frac{a_1\left[r_1\rho_2x_1\left(K_1-x_1\right)+K_1d_2\right]}{r_1\rho_2x_1\left(K_1a_2-K_1a_1-a_2\right)-r_1\rho_2x_1^2(a_2-a_1)+K_1(a_2r_1\rho_2+a_1a_2-a_1d_2)}=\frac{g_t(x_1)}{g_b(x_1)}=G(x_1)
\end{array}
\eae with $r_1=1, r_2=r$ and the following properties:
\begin{enumerate}
\item $F(0)=\frac{f_t(0)}{f_b(0)}=\frac{a_2K_2d_1}{K_2(a_1r_2\rho_1+a_1a_2-a_2d_1)}=\frac{a_2d_1}{a_1r_2\rho_1+a_1a_2-a_2d_1}$ and $F(K_2)=\frac{f_t(K_2)}{f_b(K_2)}=\frac{d_1}{a_1-d_1}=\mu_1$.
\item $f_t(x_2)=a_2\left[r_2\rho_1x_2\left(K_2-x_2\right)+K_2d_1\right]\geq a_2K_2 d_1>0 \mbox{ for } x_2\in [0, K_2]$ and
$$f_b(x_2)\big\vert_{a_1=a_2=a}=a\left[r_2\rho_1(K_2-x_2)+K_2(a-d_1)\right].$$
\item $G(0)=\frac{g_t(0)}{g_b(0)}=\frac{a_1K_1d_2}{K_1(a_2r_1\rho_2+a_1a_2-a_1d_2)}=\frac{a_1d_2}{a_2r_1\rho_2+a_1a_2-a_1d_2}$ and $G(K_1)=\frac{g_t(K_1)}{g_b(K_1)}=\frac{d_2}{a_2-d_2}=\mu_2$.
\item $g_t(x_1)=a_1\left[r_1\rho_2x_1\left(K_1-x_1\right)+K_1d_2\right]\geq a_1 K_1d_1>0\mbox{ for } x_1\in [0, K_1]$ and
$$g_b(x_1)\big\vert_{a_1=a_2=a}=a\left[r_1\rho_2(K_1-x_1)+K_1(a-d_2)\right].$$
\end{enumerate}
\begin{theorem}\label{th6:interior}[Interior equilibrium]
 If $\mu_i>K_i$ for both $i=1,2$, then Model \eqref{2patch} has no interior equilibrium. Moreover, we have the following two cases:
 \begin{enumerate}
 \item Assume that $a_i>a_j$ where $i=1,j=2$ or $i=2,j=1$. Define
 $$x^c_i=\frac{K_i\left(r_i\rho_j+a_i-d_j-\sqrt{(a_i-d_j)(r_i\rho_j+a_i-d_j)}\right)}{r_i\rho_j}.$$
 %\mbox{ and }$$x^c_2=\frac{K_2\left(r_2\rho_1+a_2-d_1-\sqrt{(a_2-d_1)(r_2\rho_1+a_2-d_1)}\right)}{r_2\rho_1}.$$
 Model \eqref{2patch} has no interior equilibrium if $ a_i>a_j, \,\rho_i<\frac{4K_ja_j(a_i-a_j)(d_i-a_i)}{r_j\left(K_ja_i-K_ja_j+a_i\right)^2}$ hold.
 %, \,\rho_2>\frac{4K_1a_1(a_2-a_1)(d_2-a_2)}{r_1\left(K_1a_2-K_1a_1+a_2\right)^2}$ hold.
 %either $$ a_1>a_2, \,\rho_1<\frac{4K_2a_2(a_1-a_2)(d_1-a_1a)}{r_2\left(K_2a_1-K_2a_2+a_1\right)^2}, \,\rho_2>\frac{4K_1a_1(a_2-a_1)(d_2-a_2)}{r_1\left(K_1a_2-K_1a_1+a_2\right)^2}$$ or
%$$a_1<a_2,\,\rho_1>\frac{4K_2a_2(a_2-a_1)(a_1a_2-a_1d_2)}{r_2\left(K_2a_1-K_2a_2+a_1\right)^2}, \,\rho_2<\frac{4K_1a_1(a_2-a_1)(d_2-a_2)}{r_1\left(K_1a_2-K_1a_1+a_2\right)^2}$$hold.

And it has at least one interior equilibrium $(x_1^*,y_1^*,x_2^*,y_2^*)$ if  the following conditions hold for both $i=1, j=2$ and $i=2,j=1$
$$a_i>\max\{a_j,d_1,d_2\}, a_j>\max\{d_1,d_2\}, \,\rho_j<\frac{4K_ia_i(a_j-a_i)(d_j-a_j)}{r_i\left(K_ia_j-K_ia_i+a_j\right)^2}, \,F(x^c_2)<K_1,\,\mbox{ and } G(x^c_1)<K_2$$ where sufficient conditions for the inequalities $F(x^c_2)<K_1,\,\mbox{ and } G(x^c_1)<K_2$ hold are
$$\rho_i\leq \frac{4(K_ia_i-K_id_i-d_i)}{K_jr_j} \mbox{ and }\rho_j<\frac{4K_ja_j(K_ia_i-K_id_i-d_i)}{a_jr_jK_j^2+r_jK_i(K_ja_j-K_ja_i-a_i)^2}.$$
In addition, we have $\frac{a_id_j}{a_jr_i\rho_j+a_ia_j-a_id_j}<x_j^*<K_j$ for both $i=1,j=2$ and $i=2,j=1$.
 \item Assume that $a_1=a_2=a$.  Model \eqref{2patch} has no interior equilibrium if $d_1>a+r_2\rho_1$ or $d_2>a+r_1\rho_2$ while it has at least one interior equilibrium $(x_1^*,y_1^*,x_2^*,y_2^*)$ if the following inequalities hold
 $$a_1=a_2=a>\max\{d_1,d_2\},\,F(x^c_2)<K_1,\,\mbox{ and } G(x^c_1)<K_2$$where sufficient conditions for the inequalities $F(x^c_2)<K_1,\,\mbox{ and } G(x^c_1)<K_2$ hold are
 $$\rho_i<\frac{4(K_ia-K_id_i-d_i)}{K_jr_j}$$ for both $i=1,j=2$ and $i=2,j=1$.
 In addition, we have $\frac{d_j}{r_i\rho_j+a-d_j}<x_j^*<K_j$ for both $i=1,j=2$ and $i=2,j=1$.
%\item Assume that $a_1=a_2=a>\max\{d_1,d_2\}$, then sufficient conditions for the inequalities $F(x^c_2)<K_1,\,\mbox{ and } G(x^c_1)<K_2$ hold are $\rho_i<\frac{4(K_ia-K_id_i-d_i)}{K_jr_j}$ for both $i=1,j=2$ and $i=2,j=1$.

 \end{enumerate}
  \end{theorem}
\noindent\textbf{Notes:} Theorem \ref{th6:interior} provides sufficient conditions on the existence of no interior equilibrium when  $\mu_i>K_i$ for either $i=1$ or $i=2$; and at  least one interior equilibrium of Model \eqref{2patch} when $\mu_i<K_i$ for both $i=1,2$. The results indicate follows:
\begin{enumerate}
\item If $\mu_i>K_i$, then Model \eqref{2patch} has no interior equilibrium if the dispersal of its predator is too small.
\item If $\mu_i<K_i$ for both $i=1,2$, then large values of the predation rate $a_i, a_j$ and small values of dispersal of both predators can lead to at least one interior equilibrium.
 \end{enumerate}
 The question is how we can solve the explicit form of an interior equilibrium of Model \eqref{2patch}.
 The following theorem provides us an example of such interior equilibrium of Model \eqref{2patch}.
%%%%%%%%%%%%%%%%%%%%%%%%%%%%%%%%%%%%%%%%%%%%%%%%%%%%%%%%%%%%%%%

\begin{theorem}\label{th7:interior}[Interior equilibrium and the stability]Suppose that $d_1=d_2=d$. Let
 $$\mu_1=\frac{d}{a_1-d}, \nu_1=\frac{(K_1-\mu_1)(1+\mu_1)}{a_1K_1}, \mu_2=\frac{d}{a_2-d}, \nu_2=\frac{r(K_2-\mu_2)(1+\mu_2)}{a_2K_2}.$$
 If $0<\mu_i<K_i$ for both $i=1$ and $i=2$, then $E^i=(\mu_1,\nu_1,\mu_2,\nu_2)$ is an interior equilibrium of Model \eqref{2patch} and its stability can be classified in the following cases:
 \begin{enumerate}
 \item $E^i$ is locally asymptotically stable if $\frac{K_i-1}{2}<\mu_i<K_i$ hold for both $i=1$ and $i=2$ while it is unstable if  the following inequality holds
 $$\frac{\mu_1(K_1a_1-a_1-K_1d-d)}{a_1K_1}+\frac{r\mu_2(K_2a_2-a_2-K_2d-d)}{a_2K_2}>0.$$
 \item Assume that
 $(K_1a_1-a_1-K_1d-d)(K_2a_2-a_2-K_2d-d)<0$\mbox{ and } $$\frac{\mu_1(K_1a_1-a_1-K_1d-d)}{a_1K_1}+\frac{r\mu_2(K_2a_2-a_2-K_2d-d)}{a_2K_2}<0.$$
 If $K_ia_i-a_i-K_id-d>0$ (i.e., $\mu_i<\frac{K_i-1}{2}$) for either $i=1$ or $i=2$, then the large values of $\rho_i$ can make $E^i$ being locally asymptotically stable, i.e., %$i=1,j=2$
{\small$$\rho_i>\max\Big\{\frac{-\nu_j-r_j\mu_i\mu_j(K_ia_i-K_id-a_i-d)(K_ja_j-K_jd-a_j-d)}{(K_iK_ja_j\nu_jd\nu_i(a_i-d)^2)},\frac{-\frac{\mu_i\nu_jK_j(\nu_i\rho_j+1)(a_j-d)^2(K_ia_i-K_id-a_i-d)}{r_j\mu_j\nu_iK_i(a_i-d)^2(K_ja_j-K_jd-a_j-d)}-1}{\nu_j}\Big\} .$$}
 \end{enumerate}
 If, in addition, $a_1=a_2=a, K_1=K_2=K, r=1$, then $\mu_1=\mu_2=\mu=\frac{d}{a-d},\,\nu_1=\nu_2=\nu=\frac{(K-\mu)(1+\mu)}{aK}$, and $E^i=(\mu,\nu,\mu,\nu)$ is the only interior equilibrium for Model \eqref{2patch} which has the same local stability as the interior equilibrium $(\mu,\nu)$ for the single patch model \eqref{Onepatch}, i.e., $E^i$ is locally asymptotically stable if
 $\frac{K-1}{2}<\mu<K$ while it is unstable if $\mu<\frac{K-1}{2}$.
\end{theorem}
\noindent\textbf{Notes:}  Theorem \ref{th7:interior} implies Model \eqref{2patch} has an interior equilibrium $E^i=(\mu_1,\nu_1,\mu_2,\nu_2)$ if $d_1=d_2=d$ and $0<\mu_i<K_i$ for both $i=1,2$. In addition, Theorem \ref{th7:interior} indicates that dispersal of predators has no effects on the local stability if $\frac{K_i-1}{2}<\mu_i<K_i$ for both $i=1,2$ or one of the single patch models \eqref{Onepatch} is unstable and $\frac{\mu_1(K_1a_1-a_1-K_1d-d)}{a_1K_1}+\frac{r\mu_2(K_2a_2-a_2-K_2d-d)}{a_2K_2}>0.$ However, large dispersal of predator at Patch $i$ can stabilize the interior equilibrium when its single patch model model is unstable at $(\mu_i,\nu_i)$ with $\frac{\mu_1(K_1a_1-a_1-K_1d-d)}{a_1K_1}+\frac{r\mu_2(K_2a_2-a_2-K_2d-d)}{a_2K_2}<0.$\\

%%%%%%%%%%%%%%%%%%%%%%%%%%%%%%%%%%%%%%%%%%%%%%%%%%%%%%%%%%%%%%%
%%%%%%%%%%%%%%%%%%%%%%%%%%%%%%%%%%%%%%%%%%%%%%%%%%%%%%%%%%%%%%%
%%%%%%%%%%%%%%%%%%%%%%%%%%%%%%%%%%%%%%%%%%%%%%%%%%%%%%%%%%%%%%%

\section{Effects of dispersal on dynamics}\label{sec_effects_of_dispersal}
From mathematical analysis in the previous sections, we can have the following summary regarding the effects of dispersal of predators for Model \eqref{2patch}:%\\\vspace{10pt}
\begin{enumerate}
\item Large dispersal of predator at Patch $i$ can stabilize or destabilize the boundary equilibrium of $x_i^*=K_i, y_i^*=0, x_j^*=\mu_j, y_j^*=\nu_j$ depending on additional conditions.
\item Small dispersal of predator at Patch $i$ may preserve its persistence under certain conditions. On the other hand, large dispersal of predator at Patch $i$ may promote its persistence when the other predator is already persist even if $\mu_i>K_i$.
\item Dispersal has no effects on the persistence of prey and the number of boundary equilibrium. It has also no effects on the local stability of the boundary equilibrium $E_{K_10K_20}$ and  the symmetric interior equilibrium $(\mu,\nu,\mu,\nu)$ when it exists.
\item If $d_i>a_i$, then small dispersal of predator at Patch $i$ prevents the interior equilibrium while if $0<\mu_i<K_i$, large predations rates $a_i, a_j$ and small dispersal of predators at both patches can lead to at least one interior equilibrium.
\item If $d_i=d_j$ and $0<\mu_i<\frac{K_i-1}{2}, \frac{K_j-1}{2}<\mu_j<K_j$, then large dispersal of predator at Patch $i$ can stabilize the interior equilibrium $(\mu_i,\nu_i,\mu_j,\nu_j)$.
\end{enumerate}
To continue our study, we will  perform bifurcations diagrams and simulations to explore the effects on the dynamical patterns and compare dynamics of our model \eqref{2patch} to the classical two patch model \eqref{2patch2}.
%%%%%%%%%%%%%%%%%%%%%%%%%%%%%%%%%%%%%%%%%%%%%%%%%%%%%%%%%%%%%%%
\subsection{Bifurcation diagrams and simulations}
In this subsection, we perform bifurcation diagrams and simulations to obtain additional insights on the effects of dispersal on the dynamics of our proposed two patch model \eqref{2patch}. We fix $r_1=1, r_2=1.5, K_1=5, K_2=3, d_1=0.2, d_2=0.1$. Then according to Theorem \ref{th1:dp}, we know that in the absence of dispersal, the dynamics of Patch 1 has global stability at $(5,0)$ if $0<a_1<0.24$; it has global stability at its unique interior equilibrium $\left(\frac{0.2}{a_1-0.2},\frac{\left(5-\frac{0.2}{a_1-0.2}\right)\left(1+\frac{0.2}{a_1-0.2}\right)}{5}\right)$ if $2<\frac{0.2}{a_1-0.2}<5\Leftrightarrow 0.24<a_1<0.3$; and it has a unique limit cycle if $a_1>0.3$; %$(4,1)$;
while the dynamics of Patch 2 has global stability at $(3,0)$ if $0<a_2<0.133$; it has global stability its unique interior equilibrium $\left(\frac{0.1}{a_2-0.1},\frac{1.5\left(3-\frac{0.1}{a_2-0.1}\right)\left(1+\frac{0.1}{a_2-0.1}\right)}{3}\right)$ if $1<\frac{0.1}{a_2-0.1}<3\Leftrightarrow 0.133<a_2<0.2$ while it has a unique limit cycle if $a_2>0.2$. Now we consider the following cases:
\begin{enumerate}
\item Choose $a_1=0.25$ and $a_2=0.15$. In the absence of dispersal, the dynamics at both Patch 1 and 2 have global stability at its unique interior equilibrium $(4,1)$, $\left(2,1.5\right)$, respectively. After turning on the dispersal, the coupled two patch model can have one interior equilibrium (see the blue regions in Figure \ref{fig_number_inteior_SS_case}) which can be locally stable (see the blue dots in Figure \ref{fig1:stability}), or be a saddle (see the green dots in Figure \ref{fig1:stability}) where the coupled system has fluctuated dynamics; or it can have two interior equilibria (see the red regions in Figure \ref{fig_number_inteior_SS_case}) which could be saddles and generate bistability between fluctuated interior dynamics and the boundary attractor $(4,1,3,0)$ (see the examples of two interior saddles of Figure \ref{fig4:stability-ss}); or it can have three interior equilibria (see the black regions in Figure \ref{fig_number_inteior_SS_case}) which generate multiple interior attractors (see the examples of two saddles, one sink and two sinks, one saddle of Figure \ref{fig1:stability}) or it could have no interior equilibrium (see white and yellow regions of Figure \ref{fig_number_inteior_SS_case}). Bifurcation diagrams Figure \ref{fig_number_inteior_SS_case},\ref{fig1:stability}, and \ref{fig4:stability-su} suggest that dispersal may destabilize system and generate fluctuated dynamics; may generate multiple interior attractors (the case of three interior equilibria), thus generate multiple attractors; or even may drive extinction of predator in one or both patches (he case of two interior equilibria, no interior equilibrium, respectively).
\begin{figure}
\begin{center}
\subfigure[$\rho_1$ V.S. $\rho_2$ for the number of interior equilibria]{\includegraphics[height = 65mm, width = 65mm]{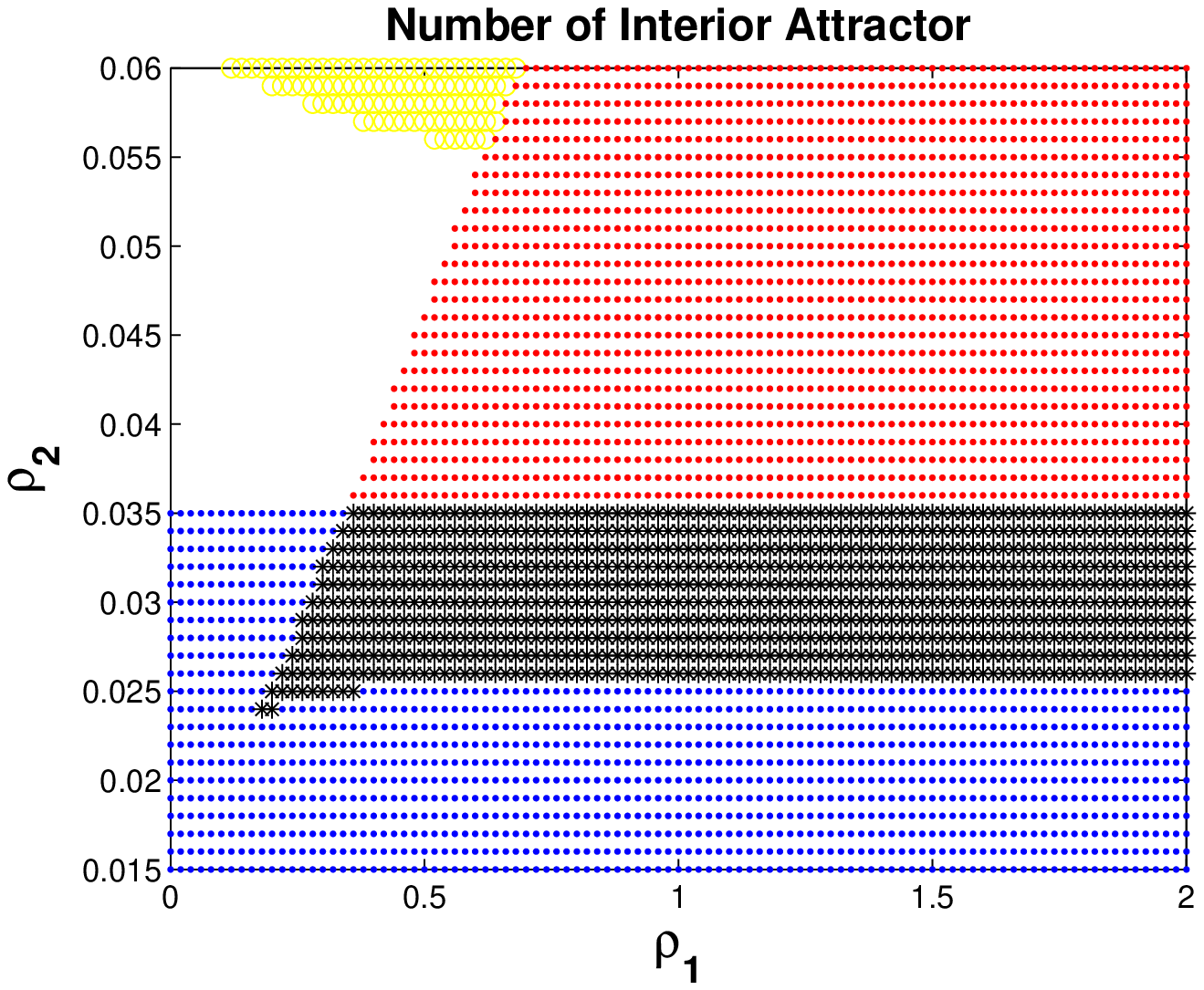}\label{fig_number_inteior_SS_case}}\hspace{2mm}
\subfigure[$\rho_1$ V.S. predator population $y_1$ when $\rho_2=0.025$]{\includegraphics[height = 65mm, width = 65mm]{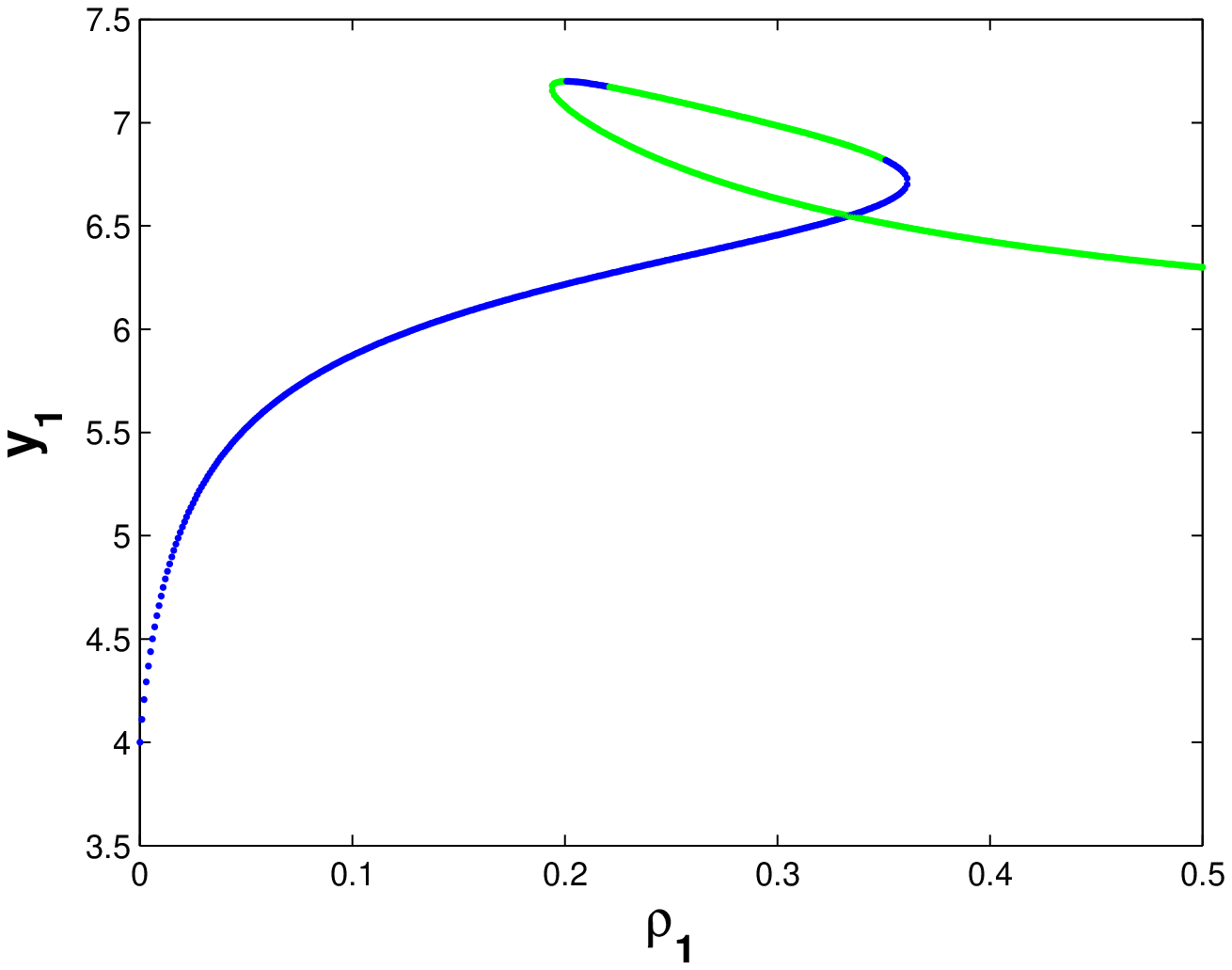}\label{fig1:stability}}
\end{center}
\caption{One and two bifurcation diagrams of Model \eqref{2patch} where $r=1.5$, $d_1=0.2$, $d_2=0.1$, $K_1=5$, $K_2=3$, $a_1=0.25$, and $a_2=0.15$. The left figure \ref{fig_number_inteior_SS_case} describes how number of interior equilibria changes for different dispersal values $\rho_i, i=1,2$: black regions have three interior equilibria; red regions have two interior equilibria; blue regions have unique interior equilibrium; yellow regions have no interior equilibrium and predator in Patch 2 dies out; white regions have no interior equilibrium and both predator die out. The right figure \ref{fig1:stability} describes the number of interior equilibria and their stability when $\rho_2=0.025$ and $\rho_1$ changes from 0 to 0.5 where $y$-axis is the population size of predator at Patch 1: Blue represents the sink; green represents the saddle; and red represents the source.}
\label{fig1:ss}
\end{figure}
\item Choose $a_1=0.25$ and $a_2=0.25$. In the absence of dispersal, the dynamics of Patch 1 has global stability at its unique interior equilibrium $(4,1)$ while the dynamics of Patch 2 has a unique stable limit cycle around $\left(2,1.5\right)$. After turning on the dispersal, the coupled two patch model can have one interior equilibrium (see the blue regions in Figure \ref{fig_number_inteior_SU_case}) which can be locally stable (see the blue dots in Figure \ref{fig2:stability}), or be a saddle (see the green dots in Figure \ref{fig2:stability}), or be a source (see the red dots in Figure \ref{fig2:stability}) where the coupled system has fluctuated dynamics for the later two cases; or it can have two interior equilibria (see the red regions in Figure \ref{fig_number_inteior_SS_case}) which could be two saddles or one sink, one saddle and generate bistability between the interior attractor and the boundary attractor (see Figure \ref{fig2:stability}); or it can have three interior equilibria (see the black regions in Figure \ref{fig_number_inteior_SS_case}) which generate multiple interior attractors (see the examples of two sinks and one saddle of Figure \ref{fig4:stability-su}) or it could have no interior equilibrium (see white regions of Figure \ref{fig_number_inteior_SU_case}). Bifurcation diagrams of Figure \ref{fig_number_inteior_SU_case}, \ref{fig2:stability}, and \ref{fig4:stability-su} suggest that dispersal may stabilize system and generate equilibrium dynamics; may generate multiple interior equilibria (the case of three interior equilibria), thus generate multiple attractors; or even may drive extinction of predator in one or two both patches (the case of two interior equilibria, no interior equilibrium, respectively).
\begin{figure}
\begin{center}
\subfigure[$\rho_1$ V.S. $\rho_2$ for the number of interior equilibria]{\includegraphics[height = 65mm, width = 65mm]{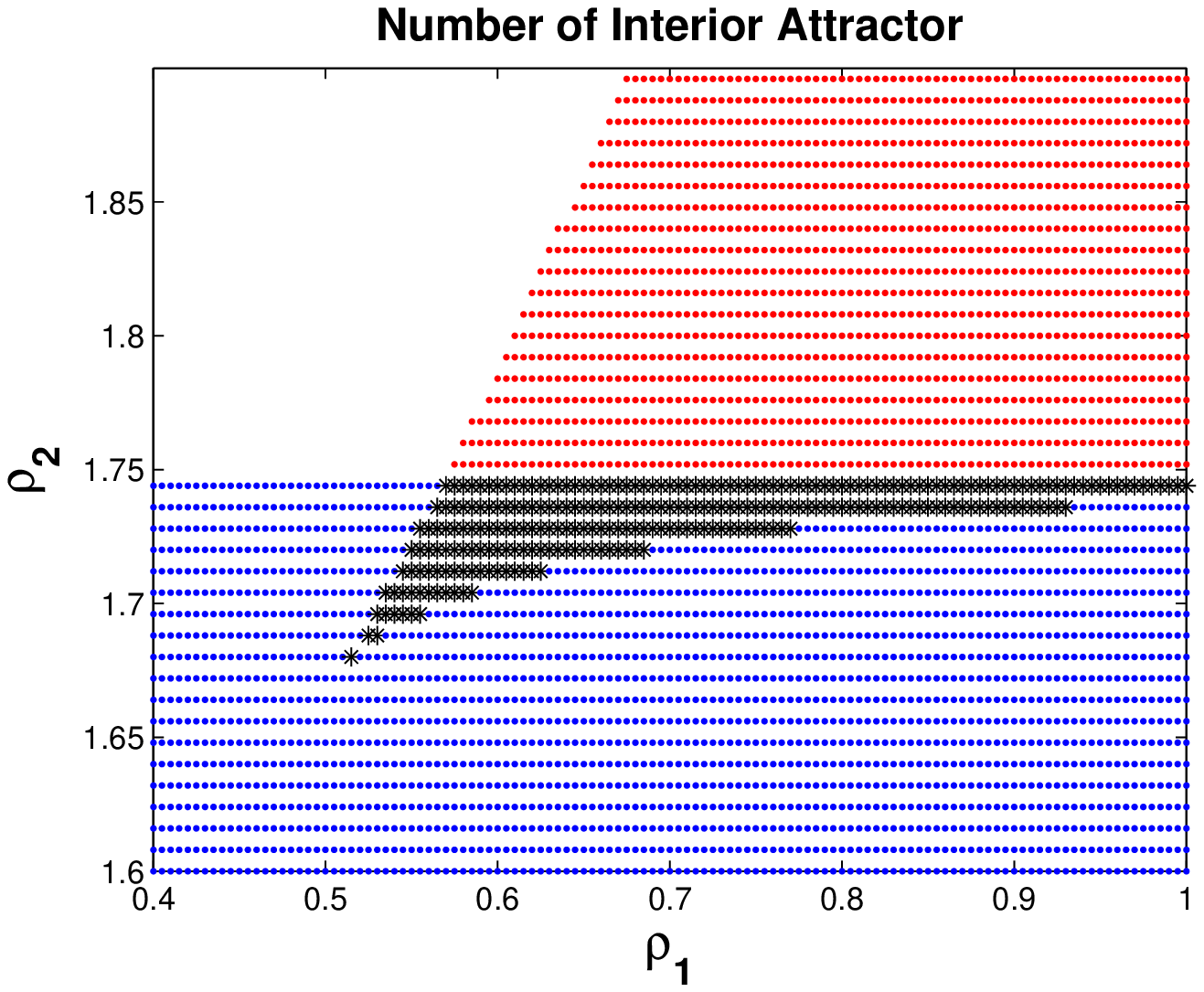}\label{fig_number_inteior_SU_case}}\hspace{2mm}
\subfigure[$\rho_2$ V.S. predator population $y_1$ when $\rho_1=1$]{\includegraphics[height = 65mm, width = 65mm]{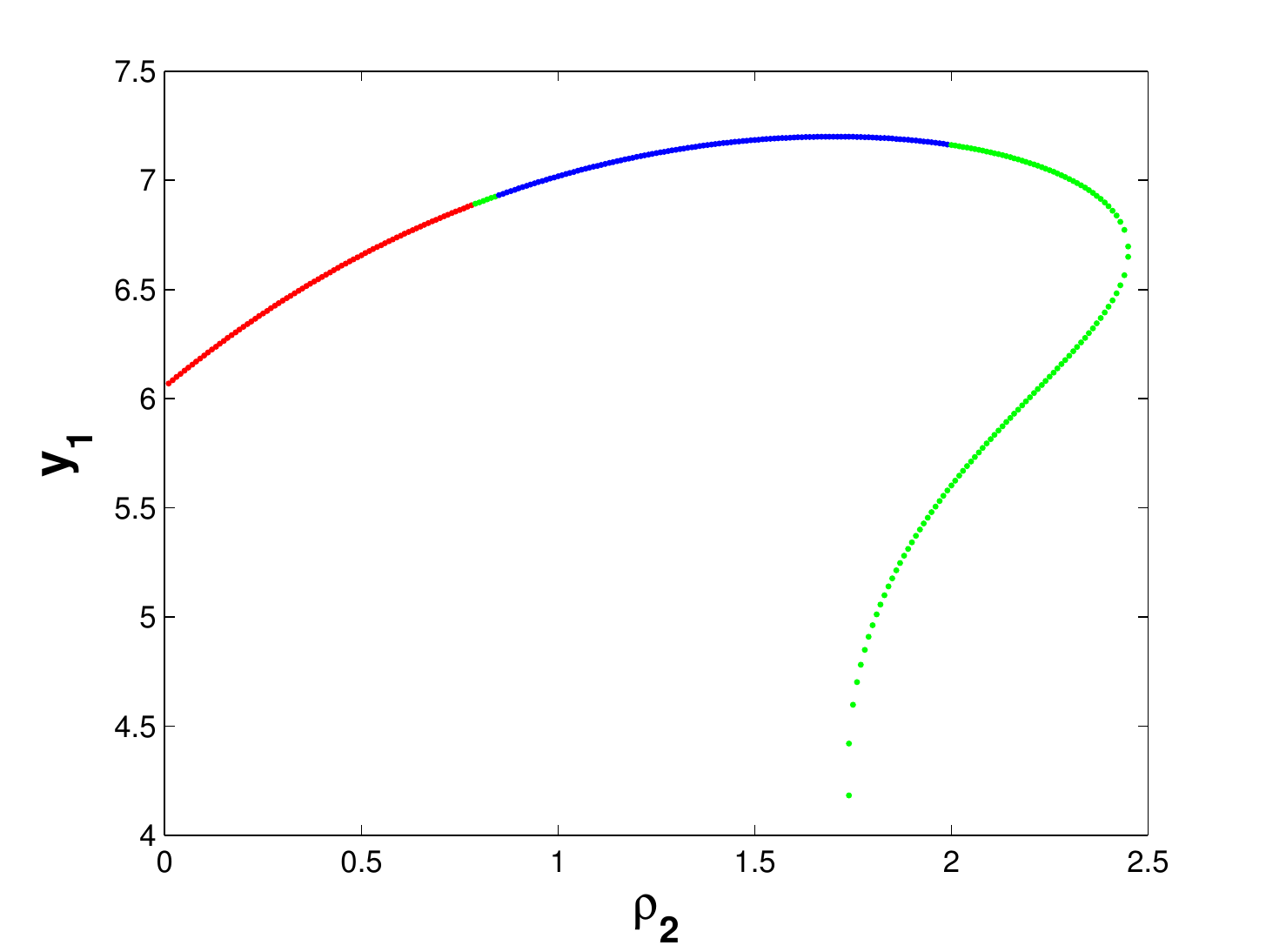}\label{fig2:stability}}
\end{center}
\caption{One and two bifurcation diagrams of Model \eqref{2patch} where $r=1.5$, $d_1=0.2$, $d_2=0.1$, $K_1=5$, $K_2=3$, $a_1=0.25$ and $a_2=0.25$. The left figure \ref{fig_number_inteior_SU_case} describes how number of interior equilibria changes for different dispersal values $\rho_i, i=1,2$: black regions have three interior equilibria; red regions have two interior equilibria; blue regions have unique interior equilibrium; yellow regions have no interior equilibrium and predator in Patch 2 dies out; white regions have no interior equilibrium and both predator die out.  The right figure \ref{fig2:stability} describes the number of interior equilibria and their stability when $\rho_1=1$ and $\rho_2$ changes from 0 to 2.5 where $y$-axis is the population size of predator at Patch 1: Blue represents the sink; green represents the saddle; and red represents the source.}
\label{fig2:su}
\end{figure}
\item Choose $a_1=0.35$ and $a_2=0.25$. In the absence of dispersal, the dynamics of both Patch 1 and 2 have a unique stable limit cycle. After turning on the dispersal, the coupled two patch model can have one interior equilibrium (see the blue regions in Figure \ref{fig_number_inteior_UU_case}) which can be a sink (see the red dots in Figure \ref{fig3:stability}) where the coupled system has fluctuated dynamics; or it can have two interior equilibria (see the red regions in Figure \ref{fig_number_inteior_SU_case}) which could be two saddles, one sink v.s. one saddle, one source v.s. one saddle and generate bistability between the interior attractors and the boundary attractor (see Figure \ref{fig3:stability}); it could have no interior equilibrium (see white and yellow regions of Figure \ref{fig_number_inteior_UU_case}). Bifurcation diagrams Figure \ref{fig_number_inteior_SU_case}-\ref{fig3:stability} suggest that dispersal may generate bistability between the interior attractor and the boundary attractor; or even may drive the extinction of predator in one or both patches (the case of two interior equilibria, no interior equilibrium, respectively).
\begin{figure}
\begin{center}
\subfigure[$\rho_1$ V.S. $\rho_2$ for the number of interior equilibria]{\includegraphics[height = 75mm, width =75mm]{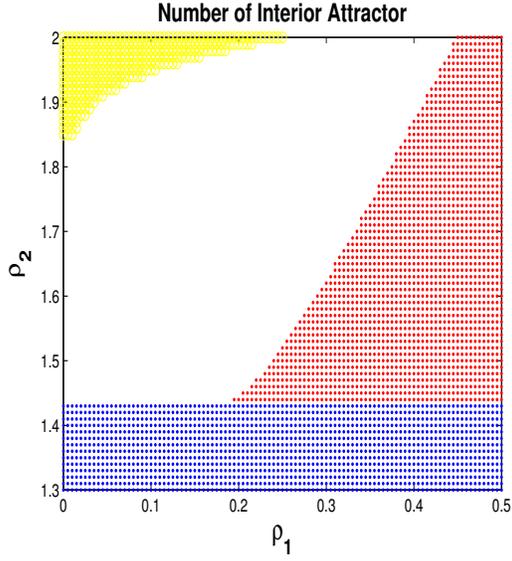}\label{fig_number_inteior_UU_case}}\hspace{2mm}
\subfigure[$\rho_2$ V.S. predator population $y_1$ when $\rho_1=1$]{\includegraphics[height = 75mm, width = 75mm]{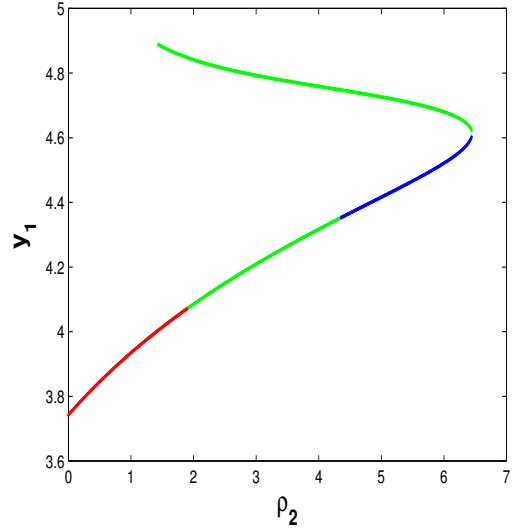}\label{fig3:stability}}
\end{center}
\caption{One and two bifurcation diagrams of Model \eqref{2patch} where $r=1.5$, $d_1=0.2$, $d_2=0.1$, $K_1=5$, $K_2=3$, $a_1=0.35$ and $a_2=0.25$. The left figure \ref{fig_number_inteior_UU_case} describes how number of interior equilibria changes for different dispersal values $\rho_i, i=1,2$: black regions have three interior equilibria; red regions have two interior equilibria; blue regions have unique interior equilibrium; yellow regions have no interior equilibrium and predator in Patch 2 dies out; white regions have no interior equilibrium and both predator die out.  The right figure \ref{fig3:stability} describes the number of interior equilibria and their stability when $\rho_1=1$ and $\rho_2$ changes from 0 to 7 where $y$-axis is the population size of predator at Patch 1: Blue represents the sink; green represents the saddle; and red represents the source.}
\label{fig3:uu}
\end{figure}
\end{enumerate}
In summary, Figure \ref{fig1:ss}, \ref{fig2:su}, \ref{fig3:uu}, and \ref{fig4:ss-su} suggest that dispersal of predator may stabilize or destabilize interior dynamics; it may drive the extinction of predator in one or both patches; and it may generate the following patterns of multiple attractors via two or three interior equilibria:

\begin{enumerate}
\item \textbf{Multiple interior attractors through three interior equilibria:} In the presence of dispersal, Model \eqref{2patch} can have the following types of interior equilibria and the corresponding dynamics:
\begin{itemize}
\item Two interior sinks and one interior saddle: Depending on the initial conditions with $x_1(0)y_1(0)x_2(0)y_2(0)>0$, Model \eqref{2patch} converges to one of two sinks for almost all initial conditions (see examples in  Figure \ref{fig1:stability}-\ref{fig4:stability-su}).
\item One interior sink and two interior saddles: Depending on the initial conditions with $x_1(0)y_1(0)x_2(0)y_2(0)>0$, Model \eqref{2patch} either converges to the sink or has fluctuated dynamics for almost all initial conditions
(see examples in Figure \ref{fig1:stability}, \ref{fig4:stability-ss}).
%\item three saddles: see Figure \ref{fig4:stability-su}
\end{itemize}
We should also expect the case of one sink v.s. one saddle v.s. one source and the case of two source v.s. one saddle when  the interior sink(s) become unstable and go through Hopf-bifurcation. In addition, Model \eqref{2patch} seems to be permanent whenever it processes three interior equilibria.
%[\textbf{Multiple interior attractors}:]
%Choose $a_2=0.5$ and $\rho_1=\rho_2=0.5$. In this case, Patch 2 is unstable at its unique interior equilibrium $(2,1.5)$ and has a unique stable limit cycle in the absence of the dispersal.
\begin{figure}
\begin{center}
\subfigure[$\rho_2$ V.S. predator population $y_2$ when $\rho_1=0.5, a_2=0.15$]{\includegraphics[height = 65mm, width = 65mm]{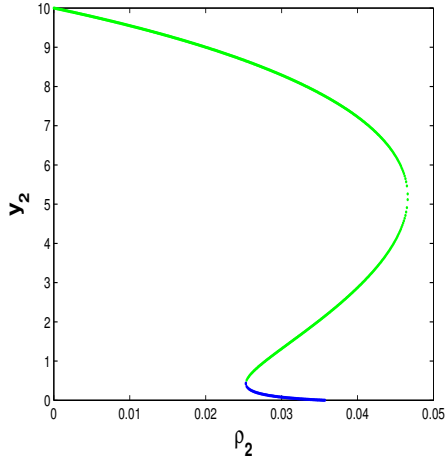}\label{fig4:stability-ss}}\hspace{2mm}
\subfigure[$\rho_2$ V.S. predator population $y_2$ when $\rho_1=0.6, a_2=0.25$]{\includegraphics[height = 65mm, width = 65mm]{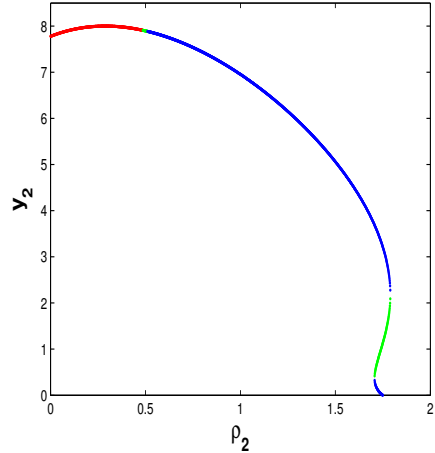}\label{fig4:stability-su}}
\end{center}
\caption{One dimensional bifurcation diagrams of Model \eqref{2patch} where $r=1.5$, $d_1=0.2$, $d_2=0.1$, $K_1=5$, $K_2=3$ and $a_1=0.25$. The left figure \ref{fig4:stability-ss} describes describes the number of interior equilibria and their stability when $\rho_1=0.5$ and $\rho_2$ changes from 0 to 0.05.  The right figure \ref{fig4:stability-su} describes the number of interior equilibria and their stability when $\rho_1=0.6$ and $\rho_2$ changes from 0 to 1.8. In both figures, blue represents the sink; green represents the saddle; and red represents the source.}
\label{fig4:ss-su}
\end{figure}
\item \textbf{Boundary attractors and interior attractors through three interior equilibria:}

\begin{itemize}
\item one interior sink and one interior saddle: Depending on the initial conditions with $x_1(0)y_1(0)x_2(0)y_2(0)>0$, Model \eqref{2patch} converges either  to the interior sinks or to the boundary attractors with one predator going extinct for almost all initial conditions (see examples in Figure \ref{fig2:stability},\ref{fig3:stability}, \ref{fig4:stability-su}).
\item two interior saddles: Depending on the initial conditions with $x_1(0)y_1(0)x_2(0)y_2(0)>0$, Model \eqref{2patch} converges either to the fluctuated interior attractors or to the boundary attractors with one predator going extinct for almost all initial conditions (see examples in Figure \ref{fig2:stability}, \ref{fig3:stability}, \ref{fig4:ss-su}).
\item one interior source and one interior saddle: Depending on the initial conditions with $x_1(0)y_1(0)x_2(0)y_2(0)>0$, Model \eqref{2patch} converges either to the fluctuated interior attractors or to the boundary attractors with one predator going extinct for almost all initial conditions (see examples in Figure \ref{fig3:stability}).
\end{itemize}
Model \eqref{2patch} has bistability between interior attractors and the boundary attractors whenever it processes two interior equilibria. This implies that depending on the initial conditions, predator at one patch can go extinct when the system has two interior equilibria.
\end{enumerate}
In general, simulations suggest that Model \eqref{2patch} is permanent when it processes one or three interior equilibria while it has bistability between interior attractors and the boundary attractors whenever it processes two interior equilibria.

%%%%%%%%%%%%%%%%%%%%%%%%%%%%%%%%%%%%%%%%%%%%%%%%%%%%%%%%%%%%%%%
\subsection{Comparisons to the classic model}
The dispersal of predator in our model is driven by the strength of prey-predator interactions. This is different from the classical dispersal model such as Model \eqref{2patch2} which has been introduced in \citep{jansen2001dynamics}:
\begin{equation}
\begin{aligned}
\frac{dx_i}{dt} & =r_i x_i\left(1 - \frac{x_i}{K_i}\right) - \frac{a_i x_i y_i}{1 + x_i} \\
\frac{dy_i}{dt} & =  \frac{a_i x_i y_i}{1 + x_i} -d_i y_i+ \rho_i( y_j-y_i)\\
\frac{dx_j}{dt} & =r_j x_j\left(1 - \frac{x_j}{K_j}\right) - \frac{a_jx_jy_j}{1 + x_j} \\
\frac{dy_j}{dt} & =  \frac{a_j x_j y_j}{1 + x_j} -d_j y_j- \rho_j( y_j-y_i)
\label{2patch2}
\end{aligned}
\end{equation}where $i=1, j=2$ or $i=2, j=1$ with $r_1=1, r_2=r$.  The symmetric case of Model \eqref{2patch2} (i.e., $r_i=r_j$, $a_i=a_j$, $K_i=K_j$, $d_i=d_j$, and $\rho_i=\rho_j$) has been discussed and studied by \citet{jansen2001dynamics} through simulations of different scenarios of local bifurcation analysis. Jansen's study shows that the classical two-patch model \eqref{2patch2} has a rich dynamical behavior where spatial predator-prey populations can be regulated through the interplay of local dynamics and migration: (i) for very small migration rates the oscillations always synchronize; (ii) For intermediate migration rates the synchronous oscillations are unstable and there are periodic, quasi-periodic, and intermittently chaotic attractors with asynchronous dynamics; and (ii) For large predator migration rates, attractors in the form of equilibria or limit cycles exist in which one of the patches contains no prey.\\

Recently, \cite{Liuyuanyuan2010} studied Model \eqref{2patch2} with both dispersal in prey and predator. Liu provide global stability of the interior equilibrium for the symmetric case and performed simulations for the asymmetric cases. Here we provide rigorous results on the persistence and permanence conditions that can be used for the comparisons to our Model \eqref{2patch} in the following theorem:

%%%%%%%%%%%%%%%%%%%%%%%%%%%%%%%%%%%%%%%%%%%%%%%%%%%%%%%%%%%%%%%
\begin{theorem}\label{th8:ds}[Summary of the dynamics of Model \eqref{2patch2}] Define
$$\hat{\mu_i}=\frac{\hat{d_i}}{a_i-\hat{d_i}},\, \hat{\nu_i}=q_i(\hat{\mu_i})=\frac{r_i(K_i-\hat{\mu_i})(1+\hat{\mu_i})}{a_iK_i}
\mbox{ and } \hat{\nu}_j^i=\frac{\rho_j\hat{\nu_i}}{d_j+\rho_j}$$ where $\hat{d_i}=d_i +\frac{\rho_id_j}{d_j+\rho_j}$.
Let $E^b_1=E_{x_1y_10y_2}=\left(\hat{\mu_1},\hat{\nu_1},0,\hat{\nu}_2^1\right)$ and $E^b_2=E_{0y_1x_2y_2}=\left(0,\hat{\nu}^2_1,\hat{\mu_2},\hat{\nu}_2\right)$.%$E_{x_1y_10y_2}=\left(\hat{\mu_1},\hat{\nu_1},0,\hat{\nu}_2^1\right)$ and $E_{0y_1x_2y_2}=\left(0,\hat{\nu}_1^2,\hat{\mu_1},\hat{\nu}_2^1\right)$.
Then we have the following summary on the dynamics of Model \eqref{2patch2}
\begin{enumerate}
\item Model \eqref{2patch2} is positively invariant and bounded in its state space $\mathbb R^4_+$ with $\limsup_{t\rightarrow\infty} x_i(t)\leq K_i$ for both $i=1$ and $i=2$.

\item \textbf{Boundary equilibria:} Model \eqref{2patch2} always has the following four boundary equilibria
$$E_{0000},\,E_{K_1000},\,E_{00K_20},\,E_{K_10K_20}$$ where the first three ones are saddles while $E_{K_10K_20}$ is locally asymptotically stable if
$$d_1+d_2+\rho_1+\rho_2>\frac{a_1K_1}{1+K_1}+\frac{a_2K_2}{1+K_2}$$ \mbox{ and }
$$\begin{array}{lcl}
&&\hat{d}_1-\frac{a_1K_1}{1+K_1} +\frac{K_2a_2\left(\frac{a_1K_1}{1+K_1}-d_1-\rho_1\right)}{(d_2+\rho_2)(1+K_2)}>0\\
&\Leftrightarrow&\\
&& \left[d_1-\frac{a_1K_1}{1+K_1}\right]\left[1-\frac{K_2a_2}{(d_2+\rho_2)(1+K_2)}\right]+\frac{\rho_1}{d_2+\rho_2}\left[d_2-\frac{K_2a_2}{(1+K_2)}\right]>0.\end{array}$$
and it is a saddle if one of the above inequalities does not hold. If $0<\hat{\mu}_i<K_i$, then the boundary equilibrium $E^b_i$
exists which is locally asymptotically stable if $\frac{K_i-1}{2}<\hat{\mu_i}<K_i, r_j<a_j \hat{\nu}_j^i$.
\item \textbf{Subsystem $i$:} If $x_j=0$, then Model \eqref{2patch2}  reduces to the following subsystem \eqref{2patch2s} with three species $x_i,y_i, y_j$:
\begin{equation}
\begin{aligned}
\frac{dx_i}{dt} & =r_i x_i\left(1 - \frac{x_i}{K_i}\right) - \frac{a_i x_i y_i}{1 + x_i} \\
\frac{dy_i}{dt} & =  \frac{a_i x_i y_i}{1 + x_i} -d_i y_i+ \rho_i( y_j-y_i)\\
%\frac{dx_j}{dt} & =r_j x_j\left(1 - \frac{x_j}{K_j}\right) - \frac{a_jx_jy_j}{1 + x_j} \\
\frac{dy_j}{dt} & = -d_j y_j- \rho_j( y_j-y_i)
\label{2patch2s}
\end{aligned}
\end{equation}whose global dynamics can be described as follows:
\begin{enumerate}
\item[3a] Prey $x_i$ is persistent for Model \eqref{2patch2s} with $\limsup_{t\rightarrow\infty} x_i(t)\leq K_i$.
\item[3b] Model \eqref{2patch2s} has global stability at $(K_i,0,0)$ if $\hat{\mu_i}>K_i$.
\item[3c] Model \eqref{2patch2s} has global stability at $(\hat{\mu_i},\hat{\nu_i},\hat{\nu}_j^i)$ if
$\frac{K_i-1}{2}<\hat{\mu_i}<K_i.$
\end{enumerate}

\item \textbf{Persistence of prey:} Prey $x_i$ persists if $\hat{\mu_j}<0$, or $\hat{\mu_j}>K_j$, or $\frac{K_j-1}{2}<\hat{\mu_j}<K_j, r_i>a_i \hat{\nu}_i^j$ hold.
%\item \textbf{Persistence of both prey $x_i$ and $x_j$:}
Both prey $x_i$ and $x_j$ persist if one of the following three conditions hold
\begin{enumerate}
\item[4(a)]\label{5a} $\hat{\mu_i}>K_i$ for both $i=1$ and $i=2$. Or $\hat{\mu_i}<0$ for both $i=1$ and $i=2$.
\item[4(b)]\label{5b} $\frac{K_i-1}{2}<\hat{\mu_i}<K_i, r_j>a_j \hat{\nu}_j^i$  for both $i=1,j=2$ and $i=2,j=1$.%and $\frac{K_j-1}{2}<\hat{\mu_j}<K_j, r_i>a_i \hat{\nu}_j^i$.
\item[4(c)]\label{5c} $\hat{\mu_i}>K_i$, or  $\hat{\mu_i}<0$, and $\frac{K_j-1}{2}<\hat{\mu_j}<K_j, r_i>a_i\hat{\nu}_i^j$ for either $i=1, j=2$ or $i=2,j=1$.
\end{enumerate}
\item \textbf{Extinction of prey $x_i$:} Prey $x_i$ goes to extinction if $\frac{K_j-1}{2}<\hat{\mu_j}<K_j \mbox{ and } \frac{r_i (K_i+1)^2}{4a_iK_i}<\hat{\nu}_j^i$. In addition, under these conditions, Model \eqref{2patch2} has global stability at $(0,\hat{\nu}_i^j,\hat{\mu_j},\hat{\nu_j})$.
\item \textbf{Persistence and extinction of predators:} Predator $y_i$ and $y_j$  have the same persistence and extinction conditions. Predators persist if $0<\mu_i<K_i$ for both $i=1$ and $i=2$ while both predators go extinct if  $\mu_i>K_i$ for both $i=1$ and $i=2$. In addition, Model \eqref{2patch2} has global stability at $(K_1,0,K_2,0)$ for the later case.
%whenever $E_{K_10K_20}$ is unstable, i.e.,
%\bae\label{persisy}\left[d_1+d_2+\rho_1+\rho_2-\left(\frac{a_1K_1}{1+K_1}+\frac{a_2K_2}{1+K_2}\right)\right]\left[ \hat{d}_1-\frac{a_1K_1}{1+K_1} +\frac{K_2a_2\left(\frac{a_1K_1}{1+K_1}-d_1-\rho_1\right)}{(d_2+\rho_2)(1+K_2)}\right]<0.\eae
%\item \textbf{Extinction of both predators:} If $\mu_i>K_i$ for both $i=1$ and $i=2$, then Model \eqref{2patch2s} has global stability at $(K_1,0,K_2,0)$.

\item \textbf{Permanence of Model \eqref{2patch2}:} Model \eqref{2patch2} is permanent if $0<\mu_i<K_i$ for both $i=1$ and $i=2$ and one of 4(a), 4(b), 4(c) hold.

\item \textbf{The symmetric case: }Let $r_1=r_2=1, a_1=a_2=a, b_1=b_2=b, d_1=d_2=d$ and $K_1=K_2=K$, then  $\mu_1=\mu_2=\mu$ and $\nu_1=\nu_2=\nu$. Therefore, we can conclude that Model \eqref{2patch2s} has global stability at $(\mu,\nu,\mu,\nu)$ if $\frac{K-1}{2}<\mu<K$. In addition, the local stability of $(\mu,\nu,\mu,\nu)$ for Model \eqref{2patch2} is the same as the local stability of $(\mu,\nu)$ for single patch models when $\rho_1=\rho_2=0$ of Model \eqref{2patch2}.
\end{enumerate}
\end{theorem}
\noindent\textbf{Notes:} Theorem \ref{th8:ds} indicates follows:
\begin{enumerate}
\item If $\mu_i>K_i$ and $0<\mu_j<K_j$, then the large dispersal of predator at Patch $i$ stabilizes $E_{K_10K_20}$.

\item Proper dispersal of predators can drive the extinction of prey in one patch.
\item Dispersal has no effects on the persistence of predator. This is different from our proposed model \eqref{2patch2}.
\end{enumerate}
To see how different types of strategies in dispersal of predators affect population dynamics of prey and predator, we start with the comparison of the boundary equilibria of our model \eqref{2patch} and the classic model \eqref{2patch2}.
Both Model \eqref{2patch} and \eqref{2patch2} always have four boundary equilibria $E_{0000} = \lx0,0,0,0\rx$, $E_{K_1000} = \lx K_1,0,0,0\rx$, $E_{00K_20} = \lx 0,0,K_2,0\rx$ and $E_{K_10K_20} = \lx K_1,0,K_2,0\rx$ among which first three are saddles in both the cases. % and the forth one is locally/globally asymptotically stable for some given conditions.
If $0<\mu_i<K_i$ for $i=1,2$, then Model \eqref{2patch} has another two boundary equilibria $E^b_{i1}$ and $E^b_{i2}$ where $E^b_{i1}$ is a saddle. If $0<\hat{\mu}_i<K_i$, then Model \eqref{2patch} has another boundary equilibrium $E^b_{i}$. We summarize and compare the dynamics of our model \eqref{2patch} with dispersal in predator driven by the strength of predation and the classical model \eqref{2patch2} with dispersal in predator driven by the difference of predator densities in Table \ref{table_comparison1}-Table \ref{table_comparison3}. We highlight effects of dynamical outcomes due to different dispersal strategies in predators between Model \eqref{2patch} and \eqref{2patch2} as follows:

\begin{enumerate}
\item The boundary equilibria: $E_{K_10K_20}$, $E_{i2}^b$ and $E_i^b$. The comparisons listed in Table \ref{table_comparison1} suggest that dispersal of predator has larger effects on the boundary equilibrium of the classic model than ours.
{\small\begin{table}[ht]
\centering
\begin{tabular}{|c|p{6.cm}|p{6.cm}|}
	\hline
Scenarios   &  Model \eqref{2patch} whose dispersal is driven by the strength of prey-predator interactions      &    Classical Model \eqref{2patch2}  whose dispersal is driven by the density of predators \\
\hline
$E_{K_10K_20}$ & LAS and GAS if $\mu_i>K_i$ for both $i=1,2$.  Dispersal has no effects on its stability. & GAS if $\mu_i>K_i$ for both $i=1,2$; While LAS if $d_1+d_2+\rho_1+\rho_2>\frac{a_1K_1}{1+K_1}+\frac{a_2K_2}{1+K_2}$ and $\lz d_1-\frac{a_1K_1}{1+K_1}\rz\lz1-\frac{a_2K_2}{(d_2+\rho_2)(1+K_2)}\rz+\frac{\rho_1}{d_2+\rho_2}\lz d_2-\frac{a_2K_2}{1+K_2}\rz>0$. {Large dispersal} may be able to stabilize the equilibrium.\\
\hline
$E_{i2}^b$ ($y_i=0$)& LAS if $\frac{K_i-1}{2}<\mu_i<K_i$ and one of the conditions \textbf{sa, sb, sc, sd} in Theorem \eqref{th2:be} holds.  {Large dispersal} has potential to either stabilize or stabilize the equilibrium. & Does not exists \\
\hline

$E_i^b$  ($x_i=0$)& Does not exists & LAS if $\frac{K_i-1}{2}<\widehat{\mu}_i<K_i$ and $r_j<a_j\hat{\nu}_j^i$. GAS if $\frac{K_i-1}{2}<\widehat{\mu}_i<K_i$ and $\frac{r_j(K_j+1)^2}{4a_jK_j}<\widehat{\nu}_i^j$. Large dispersal of predator in Patch $i$ will either destroy or destabilize the equilibrium while large dispersal of predator in Patch $j$ may stabilize the equilibrium.
\\
\hline

\end{tabular}
\caption{The comparison of boundary equilibria between Model \eqref{2patch} and Model \eqref{2patch2}. LAS refers to the local asymptotical stability, and GAS refers to the global stability.}
\label{table_comparison1}
\end{table}}

\item  Persistence and extinction of prey. According to the comparison of sufficient conditions leading either persistence or extinction of prey in a patch listed in Table \ref{table_comparison2}, we can conclude that the strength of dispersal ability of predator has huge impact on the prey for the classical model \eqref{2patch2} but not for our model \eqref{2patch}.
{\small\begin{table}[ht]
\centering
\begin{tabular}{|c|p{6.cm}|p{6.cm}|}
	\hline
Scenarios   &  Model \eqref{2patch} whose dispersal is driven by the strength of prey-predator interactions      &    Classical Model \eqref{2patch2}  whose dispersal is driven by the density of predators \\
\hline

Persistence of prey & Always persist, dispersal of predator has no effects & One or both prey persist if conditions 4. in Theorem \eqref{th8:ds} holds. Small dispersal of predator in Patch $i$ and large dispersal of predator in Patch $j$ can help the persistence of prey in Patch $i$.\\
\hline
Extinction of prey & Never extinct & $x_i$ extinct if $\frac{K_j-1}{2}<\widehat{\mu}_j<K_j$ and \\
& & $\frac{r_i(K_i+1)^2}{4a_iK_i}<\widehat{\nu}_i^j$. Large dispersal of predator in Patch $i$ can promote the extinction of prey in Patch $i$.\\
\hline

\end{tabular}
\caption{The comparison of prey persistence and extinction between Model \eqref{2patch} and Model \eqref{2patch2}}
\label{table_comparison2}
\end{table}}

\item  Persistence and extinction of predator. Simulations and the comparison of sufficient conditions leading either persistence or extinction of predators in a patch listed in Table \ref{table_comparison3}, suggest that the strength of dispersal ability of predator has profound impacts on the persistence of predator for our model \eqref{2patch} while it has no effects on  the persistence of predator for the classical model \eqref{2patch2}.
{\small\begin{table}[ht]
\centering
\begin{tabular}{|c|p{6.cm}|p{6.cm}|}
	\hline
Scenarios   &  Model \eqref{2patch} whose dispersal is driven by the strength of prey-predator interactions      &    Classical Model \eqref{2patch2}  whose dispersal is driven by the density of predators \\
\hline

Persistence of predator & Predator at Patch $j$ is persistent if Conditions in Theorem \eqref{th4:p} holds. Small dispersal of predator in Patch $j$ can help the persistence of predator in that patch. Dispersal is able to promote the persistence of predator when predator goes extinct in the single patch model. & Predators in both patches have the same persistence conditions. They persist if $0<{\mu}_i<K_i$ for $i=1,2$. Dispersal seems to have no effects in the persistence of predator. \\
\hline
Extinction of predator & Simulations suggestions (see the yellow regions of Figure \ref{fig_number_inteior_SS_case} and Figure \ref{fig_number_inteior_UU_case}) that the large dispersal of predator in Patch $i$ may lead to the its own extinction. &Predators in both patches have the same extinction conditions. They go extinct if ${\mu}_i>K_i$  or $\mu_i<0$ for $i=1,2$.\\
\hline
\end{tabular}
\caption{The comparison of predator persistence and extinction between Model \eqref{2patch} and Model \eqref{2patch2}.}
\label{table_comparison3}
\end{table}}
\item Permanence of a system depends on the persistence of each species involved in the system. Our comparisons of sufficient conditions leading to the persistence of prey and predator listed in Table \ref{table_comparison2}-\ref{table_comparison3}, indicate that dispersal of predator has important impacts in the persistence of predator in our model \eqref{2patch} while it has significant effects on the persistence of prey of the classical model \eqref{2patch2}. We can include that (i) the large dispersal of predator in a patch has potential lead to the extinction of prey (the classical model \eqref{2patch2}) or predator (our model \eqref{2patch2}) in that patch, thus destroy the permanence of the system; (ii) the small dispersal of predator in Patch $i$ with the large dispersal in Patch $j$ can promote the persistence of prey (the classical model \eqref{2patch2}) or predator (our model \eqref{2patch2}) in Patch $i$, thus promote the permanence of the system.

\item Interior equilibria: Both our model \eqref{2patch} and the classical model \eqref{2patch2} have the maximum number of three interior equilibria. However, for the symmetric case, our model \eqref{2patch}  can have the unique interior equilibrium (see Theorem \ref{th7:interior}) while the classical model can potentially process three interior equilibria \citep{jansen2001dynamics}.
\item Multiple attractors: Both our model \eqref{2patch} and the classical model \eqref{2patch2} have two types bi-stability:
(a) The boundary attractors where one of prey or predator can not sustain and the interior attractors where all four species can co-exist; and (b) Two distinct interior attractors. One big difference we observed is that for the symmetric case when each single patch model has global stability at its unique interior equilibrium, our model \eqref{2patch}  can have only one interior attractor while the classical model can potentially have two distinct interior attractors. This is due to the fact that Model \eqref{2patch} has unique interior equilibrium while Model \eqref{2patch2} can potentially process three interior equilibria as we mentioned earlier.
\end{enumerate}

%%%%%%%%%%%%%%%%%%%%%%%%%%%%%%%%%%%%%%%%%%%%%%%%%%%%%%%%%%%%%%%
%%%%%%%%%%%%%%%%%%%%%%%%%%%%%%%%%%%%%%%%%%%%%%%%%%%%%%%%%%%%%%%
%%%%%%%%%%%%%%%%%%%%%%%%%%%%%%%%%%%%%%%%%%%%%%%%%%%%%%%%%%%%%%%

\section{Discussion}\label{sec_discussion}
The idea of  ``metapopulation" originated from  \citet{levins1969some} where R. Levins used the concept to study the dynamics of pests in agricultural field in which insect pests move from site to site through migrations. Since Levin's work, many mathematical models have been applied to study prey-predator interactions between two or multiples patches that are connected through random dispersion, see examples in \cite{jansen1994theoretical, jansen2001dynamics, casal1994existence, klepac2007dispersal, jansen1995regulation, lengyel1991diffusion,levin1974dispersion, pascual1993diffusion, chewning1975migratory, Hastings1983,Adler1993, May1978,holt1985population}. The study of these metapopulation models help us get a better understanding of the dynamics of species interacting in a heterogeneous environment, and allow us to obtain a useful insight of random dispersal effects on the persistence and permanence of these species in the ecosystem. Recently, there has been increasing empirical and theoretical work on the non-random foraging movements of predators which often responses to prey-contact stimuli such as spatial variation in prey density \citep{Curio1976,kareiva1987swarms}, or different type of signals arising directly from prey \citep{Waage1977}. See more related examples of mathematical models in \citep{li2005impact, kuto2004multiple, Doebli1995,chesson1986aggregation, feng2011new, hassell1974aggregation, cressman2013two, kummel2013aphids, ghosh2011two, huang2001predator,hauzy2010density,mcmurtrie1978persistence}.  Kareiva \cite{kareiva1990population} provided a good review on varied mathematical models that deal with dispersal and spatially distributed populations and pointed out the needs of including non-random foraging movements in meta-population models. Motivated by this and the recent experimental work of immobile Aphids and Coccinellids by \cite{kummel2013aphids}, we formulate a two patch prey-predator model \eqref{2patch} with the following assumptions: (a) In the absence of dispersal the model reduced to the two uncoupled Rosenzweig-MacArthur prey-predator single patch models \eqref{Onepatch}; (b)  Prey is immobile;  and (c) Predator foraging movements are driven by the strength of prey-predator interaction. We provide basic dynamical properties such  positivity and boundedness of our model in  Theorem \ref{th1:dp}. \\

Based on our analytic results and bifurcation diagrams, we list our main findings regarding the following questions stated in the introduction how our proposed nonlinear density-dependent dispersal of predator stabilizes or destabilizes the system; how it affects the extinction and persistence of prey and predator in both patches; how it may promote the coexistence ; and how it can generate spatial population patterns of prey and predator:
 \begin{enumerate}
\item Theorem \eqref{th2:be} provides us the existence and local stability features of the eight boundary equilibria of our model \eqref{2patch}. This result indicates that large dispersal of predator in its own patch may have both stabilizing and destabilizing effects on the boundary equilibrium depending on certain conditions.  Theorem \eqref{th3:gb} gives sufficient conditions on the extinction of predator in both patches, which suggest that predator can not survive in the coupled system if predator is not able to survive at its single patch. In this case, dispersal of predator has no effect on promoting the persistence of predator but dispersal may drive predator extinct even if predator is able to persist at the single patch state (see white regions of Figure \ref{fig_number_inteior_SS_case}, \ref{fig_number_inteior_SU_case}, and \ref{fig_number_inteior_UU_case}).

%In addition to this we have summarize the dispersal effects on the boundary equilibria in the comparison Table \eqref{table_comparison_1}.

\item Theorem \eqref{th4:p} provides sufficient conditions of the persistence of prey and predator while Theorem \eqref{th5:perm}  provides sufficient conditions of the permanence of our two patch model. These results imply that under certain conditions, large dispersal of predator can promote its persistence, thus, promote the permanence of the coupled system while predator in that patch goes extinct in the absence of dispersal. Our numerical studies also suggests that large dispersal can also drive the extinction of predators in both patches (see white regions of Figure \ref{fig_number_inteior_SS_case}, \ref{fig_number_inteior_SU_case}, and \ref{fig_number_inteior_UU_case}).
\item Theorem \eqref{th6:interior} and Theorem \eqref{th7:interior} provide sufficient conditions on the existence and the local stability of the interior equilibria under certain conditions. Our analytic study shows that large dispersal of predator may be able to  stabilize the interior equilibrium when one of the single patch has global stable interior equilibrium while the other one has limit cycle dynamics. At the mean time, our bifurcation diagrams (see Figure \ref{fig1:stability}, \ref{fig2:stability}, and \ref{fig2:stability}) suggest that the stabilizing or destabilizing effects of predator's dispersal are not definite, i.e., dispersal can either stabilize or destabilize  the system depending on other life history parameters. Moreover,  our simulations also suggest that the dispersal of predator can either generate multiple interior equilibria or destroy the interior equilibrium which leads to the extinction of predator in one patch or predators in both patches.

\end{enumerate}

\noindent\textbf{Comparisons to the classic model \eqref{2patch2}:} We provide detailed comparison between the dynamics of our model \eqref{2patch} to the classic model \eqref{2patch2}. These comparisons suggest that the mode of forging movement of predator has profound impacts on the dynamics of the coupled two patch model.
 %persistence and permanence of the coupled system as well as the interior equilibrium for the symmetric case.
 Here we highlight two significant differences: (1) the strength of dispersal ability of predator has profound impacts on the persistence of predator for our model \eqref{2patch} while it has no effects on the persistence of predator for the classical model \eqref{2patch2}. However,  the dispersal of predator has huge impacts on the persistence of prey for the classical model \eqref{2patch2} while it has little or no effects on the persistence of prey for our model \eqref{2patch}. And (2) for the symmetric case, our model \eqref{2patch} has a unique interior equilibrium while the classical model \eqref{2patch2} can have up to three interior equilibria thus it is able to generate different spatial patterns.  \\

\noindent\textbf{Future work:} Our study combined with the literature study on the classical model \eqref{2patch2} by \cite{Liuyuanyuan2010, jansen1994theoretical, jansen2001dynamics,Hastings1983,levin1974dispersion}, provide us a better understanding on how different dispersal behavior of predator could have different effects on the dynamical outcomes and spatial pattens. In nature, predator may have different foraging behavior as pointed out by \citet{kummel2013aphids} that the foraging movements of predator Coccinellids are combinations of passive diffusion, conspecific attraction, and retention on plants with high aphid numbers. It will be interesting to the extended version of Model \eqref{2patch} and  Model \eqref{2patch2} by incorporating two different modes of foraging behavior. One potential example is showed as follows:

\begin{equation}
\begin{aligned}
\frac{dx_1}{dt} & = x_1\left(1 - \frac{x_1}{K_1}\right) - \frac{a_1 x_1 y_1}{1 + x_1} \\
\frac{dy_1}{dt} & =  \frac{a_1 x_1 y_1}{1 + x_1} -d_1 y_1+ \rho_1\textbf{s}\left(  \frac{a_1 x_1 y_1}{1 + x_1} y_2- \frac{a_2 x_2 y_2}{1 + x_2} y_1\right)+\rho_1\textbf{(1-s)}\left( y_2-y_1\right)\\
\frac{dx_2}{dt} & = r x_2\left(1 - \frac{x_2}{K_2}\right) - \frac{a_2 x_2 y_2}{1 + x_2} \\
\frac{dy_2}{dt} & =  \frac{a_2 x_2 y_2}{1 + x_2} -d_2 y_2+\rho_2\textbf{s}\left(  \frac{a_2 x_2 y_2}{1 + x_2} y_1- \frac{a_1 x_1 y_1}{1 + x_1} y_2\right)+\rho_2\textbf{(1-s)}\left( y_1-y_2\right)
 \label{2patch2d}
\end{aligned}
\end{equation} where $s$ is a real number in $[0,1]$ indicating the portion of predator using the dispersal strategy driven by the strength of the predation (our model \eqref{2patch}) and  $1-s$ indicates the portion of predator using the dispersal strategy driven by the density difference of  predator in two patches (the classical model \eqref{2patch2}). It will be even more interesting to develop a two patch model with adaptive dispersal strategies by letting $s$ change over time and depend on the fitness of predator. These are ongoing research projects by the authors. We would also like to point out that Bruder \emph{et al.} (Bruder, A., Thompson, H., Brown, D. \& Kummel, M: Pattern formation in a two-patch predator prey-model with diffusion and attraction to predation, working on progressing) is working on Model \eqref{2patch2d} focusing on the spatial patterns generated by these two strategies of dispersal in predator.

%%%%%%%%%%%%%%%%%%%%%%%%%%%%%%%%%%%%%%%%%%%%%%%%%%%%%%%%%%%%%%%
%%%%%%%%%%%%%%%%%%%%%%%%%%%%%%%%%%%%%%%%%%%%%%%%%%%%%%%%%%%%%%%

\section{Proofs}
\subsection*{Proof of Theorem \ref{th1:dp}}

\begin{proof}
Notice that both $\frac{dx_i}{dt} \big\vert_{x_i=0}=0$ and $\frac{dy_i}{dt} \big\vert_{y_i=0}=0$ for $i=1,2$, thus according to Theorem A.4 (p.423) in \cite{thieme2003mathematics}, we can conclude that the model \eqref{2patch} is positive invariant in $\mathbb R^4_+$. Now we can go ahead to show the boundedness of the system. First, we have the following inequalities due to the property of positive invariance:
$$\frac{dx_i}{dt}= r_i x_i\left(1 - \frac{x_i}{K_i}\right) - \frac{a_i x_i y_i}{1 + x_i}\leq r_i x_i\left(1 - \frac{x_i}{K_i}\right)$$ which implies that
$$\limsup_{t\rightarrow\infty} x_i(t)\leq K_i.$$
Define $V=\rho_2 (x_1+y_1)+\rho_1 (x_2+y_2)$, then we have
$$\begin{array}{lcl}
\frac{dV}{dt}&=&\rho_2\frac{d(x_1+y_1)}{dt}+\rho_1\frac{d(x_2+y_2)}{dt}\\
&=&\rho_2x_1\left(1 - \frac{x_1}{K_1}\right) +\rho_1 rx_2\left(1 - \frac{x_2}{K_2}\right)  -\rho_2 d_1 y_1 -\rho_1 d_2 y_2\\
&=&\rho_2x_1\left(1 - \frac{x_1}{K_1}\right) +\rho_1 rx_2\left(1 - \frac{x_2}{K_2}\right)+\rho_2 d_1 x_1 +\rho_1 d_2 x_2  -\rho_2 d_1 (x_1+y_1) -\rho_1 d_2 (x_2+y_2)\\
&\leq& M-d \left[\rho_2 (x_1+y_1)+\rho_1 (x_2+y_2)\right]=M-d V
\end{array}$$ where $d=\min\{d_1,d_2\}$ and
$$M=\max_{0\leq x_1\leq K_1}\Big\{\rho_2x_1\left(1 - \frac{x_1}{K_1}+ d_1  \right)\Big\} +\max_{0\leq x_2\leq K_2}\Big\{\rho_1x_2\left(1 - \frac{x_2}{K_2}+d_2  \right)\Big\}.$$
Therefore, we have
$$\limsup_{t\rightarrow\infty} V(t)=\limsup_{t\rightarrow\infty} \rho_2 (x_1(t)+y_1(t))+\rho_1 (x_2(t)+y_2(t))\leq \frac{M}{d}$$which implies that Model \eqref{Onepatch} is bounded in $\mathbb R^4_+$.

If there is no dispersal in predator, i.e., $\rho_i=0, i=1,2$, we can easily check that Model \eqref{2patch} is reduced to the two uncoupled Rosenzweig-MacArthur prey-predator single patch models \eqref{Onepatch} with $r_1=1$ and $r_2=r$. The global dynamics of the single patch model \eqref{Onepatch} can be summarized from the work of \cite{liu2003complex, Liuyuanyuan2010,hsu1977mathematical, hsu1978global}. Thus, we omit the detailed proof here.

Recall that both $\frac{dx_i}{dt} \big\vert_{x_i=0}=0$ and $\frac{dy_i}{dt} \big\vert_{y_i=0}=0$ for $i=1,2$, therefore, the sets $\{(x_1,y_1,x_2,y_2)\in\mathbb R^4_+:x_i=0 \}$ and $\{(x_1,y_1,x_2,y_2)\in\mathbb R^4_+: y_i=0 \}$ are invariant. This indicates that if $x_j(0)=0$, then $x_j(t)=0$ for all $t>0$. Therefore, the population of $y_j$ converges to 0 since
$$\frac{dy_j}{dt} =  -d_j y_j-\rho_j\frac{a_i x_i y_i}{1 + x_i} y_j \leq 0\Rightarrow \limsup_{t\rightarrow\infty} y_j(t)=0.$$
Applying the results in \cite{Markus1956}, we can conclude that Model \eqref{2patch} is reduced to the single patch model \eqref{Onepatch} when $x_j=0$.

%Markus, L. (1956): Asymptotically autonomous di®erential systems. Contributions to the Theory of Nonlinear Oscillations III (S. Lefschetz, ed.), 17-29. Annals of Mathematics Studies 36, Princeton Univ. Press
In the case that $y_j=0$,  Model \eqref{2patch} is reduced to Model \eqref{yj0} by replacing $y_j=0$ in Model \eqref{2patch}.
%we have the dynamics of $x_j$ can be described by  $$\frac{dx_j}{dt} =r_j x_j\left(1 - \frac{x_j}{K_j}\right) .$$And

Summarizing the discussions above, we can conclude that the statement of Theorem \ref{th1:dp} holds.
 \end{proof}
\subsection*{Proof of Theorem \ref{th2:be}}
\begin{proof}According Theorem \ref{th1:dp}, sufficient condition for the single patch model \eqref{Onepatch} having the unique interior equilibrium $(\mu_i,\nu_i),i=1,2$ is $\mu_i<K_i$. Therefore, sufficient condition for Model \eqref{2patch} having boundary equilibria $E_{\mu_1\nu_100}$ and$E_{\mu_1\nu_1K_20}$ is $\mu_1<K_1$. Similarly,  sufficient condition for Model \eqref{2patch} having boundary equilibria $E_{00\mu_2\nu_2}$ and $E_{K_10\mu_2\nu_2}$ is $\mu_2<K_2$.\\

The local stability of an equilibrium $(x_1^*,y_1^*,x_2^*,y_2^*)$ can be determined by the eigenvalues $\lambda_i$, $i = 1,~2,~3,~4$ of the Jacobian matrix $J_{(x_1^*,y_1^*,x_2^*,y_2^*)}$ \eqref{JE} of Model \eqref{2patch} evaluated at the equilibrium.

{\small\bae\label{JE}
\begin{array}{ll}
&J_{(x_1^*,y_1^*,x_2^*,y_2^*)}=\\\\
&\left[\begin{array}{cccc}
\lx1-\frac{2x_1^*}{K_1}\rx-\frac{a_1y_1^*}{\lx1+x_1^*\rx^2} & -\frac{a_1x_1^*}{1+x_1^*} & 0 & 0\\

\frac{a_1y_1^*\lx 1+\rho_1y_2^*\rx}{\lx1+x_1^*\rx^2} & \frac{a_1x_1^*\lx 1+\rho_1y_2^*\rx}{1+x_1^*}-\frac{\rho_1a_2x_2^*y_2^*}{1+x_2^*}-d_1 & -\frac{\rho_1a_2y_1^*y_2^*}{\lx1+x_2^*\rx^2} & \rho_1y_1^*\lx\frac{a_1x_1^*}{1+x_1^*} - \frac{a_2x_2^*}{1+x_2^*}\rx \\

0 & 0 & r\lx1-\frac{2x_2^*}{K_2}\rx-\frac{a_2y_2^*}{\lx1+x_2^*\rx^2} & -\frac{a_2x_2^*}{1+x_2^*} \\

-\frac{\rho_2a_1y_1^*y_2^*}{\lx1+x_1^*\rx^2} & \rho_2y_2^*\lx\frac{a_2x_2^*}{1+x_2^*} - \frac{a_1x_1^*}{1+x_1^*}\rx & \frac{a_2y_2^*\lx 1+\rho_2y_1^*\rx}{\lx1+x_2^*\rx^2} & \frac{a_2x_2^*\lx 1+\rho_2y_1^*\rx}{1+x_2^*}-\frac{\rho_2a_1x_1^*y_1^*}{1+x_1^*}-d_2 \\
\end{array}\right]
\end{array}
\eae}

After substituting the boundary equilibria $E_{0000}$, $E_{K_1000}$, $E_{00K_20}$, $E_{\mu_1\nu_100}$ and $E_{00\mu_2\nu_2}$ into the Jacobian Matrix \eqref{JE}, we can conclude that these equilibria are saddles since they have both positive and negative eigenvalues. \\

The eigenvalues of \eqref{JE} evaluated at $E_{K_10K_20}$ are as follows:
$$\lambda_1 = -1~(<0),~~\lambda_2 = \frac{a_1K_1}{1+K_1}-d_1<0\Leftrightarrow\mu_1>K_1,~~\lambda_3 = -r~(<0),~~\lambda_4 = \frac{a_2K_2}{1+K_2}-d_2<0\Leftrightarrow\mu_2>K_2.$$
Therefore, $E_{K_10K_20}$ is locally asymptotically stable if $\mu_i>K_i,i=1,2$ while it is a saddle if either $\left(\mu_1-K_1\right)\left(\mu_2-K_2\right)<0$ or $\mu_i<K_i, i=1,2$ holds. \\

Now we focus on the local stability of $E_{\mu_1\nu_1K_20}$ and $E_{K_10\mu_2\nu_2}$ when they exist.
After substituting the boundary equilibrium $E_{\mu_1\nu_1K_20}$ to \eqref{JE}, we can obtain the eigenvalues of the Jacobian matrix evaluated at this boundary equilibrium as follows:

 $$\lambda_1 = -r,\,\lambda_2=\frac{K_2(a_2-d_2)-d_2}{1+K_2}+\rho_2\frac{\nu_1\left[K_2(a_2-d_1)-d_1\right]}{(1+K_2)}$$ and
  %$$\lambda_1 = -r,\,\lambda_2=\frac{K_2(a_2-d_2)-d_2}{1+K_2}+\rho_2\frac{(K_2(a_2-d_1)-d_1)(K_1(a_1-d_1)-d_1)}{K_1(1+K_2)(a_1-d_1)^2}$$
 $$\lambda_3+\lambda_4=\frac{K_1(a_1-d_1)-(a_1+d_1)}{a_1K_1(a_1-d_1)},\,\lambda_3\lambda_4=\frac{d_1(K_1(a_1-d_1)-d_1)}{a_1K_1}.$$
 Notice that the eigenvalues of $\lambda_3$ and $\lambda_4$ being negative is equivalent to the case that the unique interior equilibrium $(\mu_1,\nu_1)$ being locally asymptotically stable for the single patch model \eqref{Onepatch} when $i=1$. Thus, we can conclude that $\frac{K_1-1}{2}<\mu_1<K_1$ are sufficient conditions for $\lambda_3$ and $\lambda_4$ being negative.  Now we explore sufficient conditions for $\lambda_2$ being negative. First, we have $\mu_1<K_1$ %$\Leftrightarrow K_1(a_1-d_1)-d_1>0$
 due to the existence of $E_{\mu_1\nu_1K_20}$. We have the following three cases:
 \begin{enumerate}
  \item If $\mu_2>K_2 \Leftrightarrow K_2(a_2-d_2)-d_2<0$, then the first term of $\lambda_2$ is negative. This also implies that  Model \eqref{2patch} has no boundary equilibria  of $E_{00\mu_2\nu_2}$ and $E_{K_10\mu_2\nu_2}$. Since $\mu_1<K_1$, %$K_1(a_1-d_1)-d_1>0$,
  therefore, we have $\lambda_2<0$ for all $\rho_2>0$ if $a_2\leq d_1$ or $K_2<\frac{d_1}{a_2-d_1}$ since
  $$K_2(a_2-d_1)-d_1<0 \Leftrightarrow \mbox{ either } a_2\leq d_1 \mbox{ or }K_2<\frac{d_1}{a_2-d_1}.$$Therefore, we can conclude that $\lambda_2$ is negative if \mbox{ either } $a_2\leq d_1, K_2<\mu_2$ \mbox{ or } $K_2<\min\Big\{\mu_2,\frac{d_1}{a_2-d_1}\Big\}.$
 Assume that $K_2(a_2-d_1)-d_1>0$ (i.e., $K_2>\frac{d_1}{a_2-d_1}>0$), then $\lambda_2<0$ if $\rho_2$ is small enough, i.e., satisfies the condition  of $\rho_2<\frac{d_2-K_2(a_2-d_2)}{\nu_1\left[K_2(a_2-d_1)-d_1\right]}.$ In this case, we can conclude that $\lambda_2$ is negative if
 $$0<\frac{d_1}{a_2-d_1}<K_2<\mu_2 \mbox{ and }\rho_2<\frac{d_2-K_2(a_2-d_2)}{\nu_1\left[K_2(a_2-d_1)-d_1\right]}.$$
 \item If $\mu_2<K_2 \Leftrightarrow K_2(a_2-d_2)-d_2>0$, then the first term of $\lambda_2$ is positive. This also implies that  Model \eqref{2patch} has two boundary equilibria of $E_{00\mu_2\nu_2}$ and $E_{K_10\mu_2\nu_2}$. In this case, sufficient conditions for $\lambda_2$ being negative are $K_2<\frac{d_1}{a_2-d_1}$ and $\rho_2$ large enough. More specifically, $\rho_2$ has to satisfy the following inequality:
 $$\rho_2>\frac{K_2(a_2-d_2)-d_2}{\nu_1\left[d_1-K_2(a_2-d_1)\right]}.$$
 Therefore, $\lambda_2$ is negative if $\mu_2<K_2<\frac{d_1}{a_2-d_1} \mbox{ and }\rho_2>\frac{K_2(a_2-d_2)-d_2}{\nu_1\left[d_1-K_2(a_2-d_1)\right]}.$
 \end{enumerate}
 Summarizing the discussions above, we can conclude that the boundary equilibrium $E_{\mu_1\nu_1K_20}$ is locally asymptotically stable if $\frac{K_1-1}{2}<\mu_1<K_1$ and one of the following conditions holds:
 \begin{enumerate}
 \item $a_2\leq d_1, K_2<\mu_2$.
 \item $K_2<\min\Big\{\mu_2,\frac{d_1}{a_2-d_1}\Big\}.$
 \item $0<\frac{d_1}{a_2-d_1}<K_2<\mu_2 \mbox{ and }\rho_2<\frac{d_2-K_2(a_2-d_2)}{\nu_1\left[K_2(a_2-d_1)-d_1\right]}.$
 \item $\mu_2<K_2<\frac{d_1}{a_2-d_1} \mbox{ and }\rho_2>\frac{K_2(a_2-d_2)-d_2}{\nu_1\left[d_1-K_2(a_2-d_1)\right]}.$
 \end{enumerate}
 And $E_{\mu_1\nu_1K_20}$ is a saddle if $\mu_1<\frac{K_1-1}{2}$ or one of the following conditions holds:
 \begin{enumerate}
 \item $K_2>\max\Big\{\mu_2,\frac{d_1}{a_2-d_1}\Big\}.$
  \item $0<\frac{d_1}{a_2-d_1}<K_2<\mu_2 \mbox{ and }\rho_2>\frac{d_2-K_2(a_2-d_2)}{\nu_1\left[K_2(a_2-d_1)-d_1\right]}.$
 \item $\mu_2<K_2<\frac{d_1}{a_2-d_1} \mbox{ and }\rho_2<\frac{K_2(a_2-d_2)-d_2}{\nu_1\left[d_1-K_2(a_2-d_1)\right]}.$
 \end{enumerate}
Similarly, we can obtain sufficient conditions for the local stability of the boundary equilibrium $E_{K_10\mu_2\nu_2}$ as the statement.

If $\mu_i<K_i$, then Model \eqref{2patch} has the boundary equilibria $E_{\mu_1\nu_1K_20}$ and $E_{K_10\mu_2\nu_2}$ according to Theorem \ref{th2:be} and the discussions above. If both $E_{\mu_1\nu_1K_20}$ and $E_{K_10\mu_2\nu_2}$ are locally stable, then the following inequalities are satisfied:
$$\mu_1<K_1<\frac{d_2}{a_1-d_2} \Rightarrow \frac{d_1}{a_1-d_1}<\frac{d_2}{a_1-d_2}\Rightarrow d_1<d_2$$
and
$$\mu_2<K_2<\frac{d_1}{a_2-d_1} \Rightarrow \frac{d_2}{a_2-d_2}<\frac{d_1}{a_2-d_1}\Rightarrow d_2<d_1$$ which are contradiction. Therefore, $E_{\mu_1\nu_1K_20}$ and $E_{K_10\mu_2\nu_2}$ can not be local stable at the same time.

Now if $r_i=1, a_i=a,d_i=d, K_i=d$ for both $i=1,2$, then $E_{i2}^b, i=1,2$ exist if $Ka-Kd-d>0$.
This implies that one of the eigenvalues of the Jacobian matrix of Model \eqref{2patch} evaluated at  $E_{i2}^b$ is positive, i.e.,
$$\frac{K(a-d)^(Ka-Kd-d)+\rho_j (Ka-Kd-d)^2}{K(K+1)(a-d)^2}>0$$ which indicates that $E_{i2}^b$ can not be stable for both $i=1$ and $i=2$.
\end{proof}
%%%%%%%%%%%%%%%%%%%%%%%%%%%%%%%%%%%%%%%%%%%%%%%%%%%%%%%%%%%%%%%

\subsection*{Proof of Theorem \ref{th3:gb}}
\begin{proof}
Let $p_i(x)=\frac{a_i x}{1+x}$ and $q_i(x)=\frac{r_i(K_i-x)(1+x)}{a_iK_i}$, then we have
$$r_ix_i\left(1-\frac{x_i}{K_i}\right)-\frac{a_i x_iy_i}{(1+x_i)}=\frac{a_i x_i}{1+x_i}\left[\frac{r_i(K_i-x_i)(1+x_i)}{a_iK_i}-y_i\right]=p_i(x_i)\left[q_i(x_i)-y_i\right].$$
We construct the following Lyapunov functions
\begin{equation}\label{Lyapunov1}
\begin{aligned}
V_1(x_1,y_1) =\rho_2 \int^{x_1}_{K_1} \frac{p_1(\xi)-p_1(K_1)}{p_1(\xi)}d\xi + \rho_2y_1
\end{aligned}
\end{equation}
and
\begin{equation}\label{Lyapunov2}
\begin{aligned}
V_2(x_2,y_2) =\rho_1 \int^{x_2}_{K_2} \frac{p_2(\xi)-p_2(K_2)}{p_2(\xi)}d\xi + \rho_1y_2
\end{aligned}
\end{equation}
Now taking derivatives of the functions \eqref{Lyapunov1} and \eqref{Lyapunov2} with respect to time $t$, we get
\begin{equation}\label{L1}
\begin{array}{lcl}
\frac{d}{dt}V_1(x_1(t),y_1(t)) &=& \rho_2\frac{p_1(x_1)-p_1(K_1)}{p_i(x_1)}\frac{dx_1}{dt} + \rho_2\frac{dy_1}{dt} \\
                               &=& \rho_2\lz p_1(x_1)-p_1(K_1)\rz\lz q_1(x_1)-y_1\rz +\rho_2 y_1\lz p_1(x_1)-d_1\rz + \rho_1\rho_2 y_1y_2\lz p_1(x_1)-p_2(x_2)\rz\\
                               &=&\rho_2\lz p_1(x_1)-p_1(K_1)\rz q_1(x_1)+ \rho_2 y_1\lz p_1(K_1)-d_1\rz+\rho_1\rho_2 y_1y_2\lz p_1(x_1)-p_2(x_2)\rz
\end{array}
\end{equation}
and
\begin{equation}\label{L2}
\begin{array}{lcl}
\frac{d}{dt}V_2(x_2(t),y_2(t)) &=&\rho_1 \frac{p_2(x_2)-p_2(K_2)}{p_2(x_2)}\frac{dx_2}{dt} +\rho_1 \frac{dy_2}{dt} \\
                               &=& \rho_1\lz p_2(x_2)-p_2(K_2)\rz\lz q_2(x_2)-y_2\rz +\rho_1 y_2\lz p_2(x_2)-d_2\rz + \rho_1\rho_2 y_1y_2\lz p_2(x_2)-p_1(x_1)\rz\\
                               &=&\rho_1\lz p_2(x_2)-p_2(K_2)\rz q_2(x_2)+ \rho_1 y_2\lz p_2(K_2)-d_2\rz+\rho_1\rho_2 y_1y_2\lz p_2(x_2)-p_1(x_1)\rz
\end{array}.
\end{equation}
Let $V=V_1+V_2$. Now adding \eqref{L1} and \eqref{L2}, we get
\begin{equation}\nonumber
\begin{aligned}
\frac{d}{dt}V =\frac{d}{dt}V_1(x_1(t),y_1(t)) + \frac{d}{dt}V_2(x_2(t),y_2(t)) &=&  \rho_2\lz p_1(x_1)-p_1(K_1)\rz q_1(x_1)+ \rho_2 y_1\lz p_1(K_1)-d_1\rz \\
  &+& \rho_1\lz p_2(x_2)-p_2(K_2)\rz q_2(x_2)+ \rho_1 y_2\lz p_2(K_2)-d_2\rz.
\end{aligned}
\end{equation}
Since $\mu_i> K_i\Leftrightarrow \frac{d_i}{a_i-d_i}>K_i\Leftrightarrow \frac{a_iK_i}{1+K_i}=p_i(K_i)<d_i$, thus, we have $p_i(K_i)-d_i<0$. Notice that $p_i(x_i)$ is an increasing function in $x_i$, therefore, $\lz p_i(x_i)-p_i(K_i)\rz$ is positive for $x_i>K_i$ and it is negative for $x_i<K_i$. At the mean time, we have $q_i(x_i)$ is positive for $x_i<K_i$ and it is negative for $x_i>K_i$. This implies that $\lz p_i(x_i)-p_i(K_i)\rz q_i(x_i) \leq 0$ for all $x_i\geq0$. Therefore, we have $\frac{d}{dt}V<0$ in $\mathbb R^4_+$. This implies that both $V_1$ and $V_2$ are Lyapunov functions, and the boundary equilibrium $E_{K_10K_20} = \lx K_1,0,K_2,0\rx$ is globally stable when $\mu_i> K_i$ according to Theorem $3.2$ in \cite{hsu1978global}.
\end{proof}
%%%%%%%%%%%%%%%%%%%%%%%%%%%%%%%%%%%%%%%%%%%%%%%%%%%%%%%%%%%%%%%
\subsection*{Proof of Theorem \ref{th4:p}}
\begin{proof}
According to Theorem \ref{th1:dp}, we know that Model \eqref{2patch} is attracted to a compact set $C$ in $\mathbb R^4_+$. Moreover, if $y_j=0$, Model \eqref{2patch} is reduced to the two uncoupled models \eqref{yj0} while if $x_j=0$ it is reduced to a single patch model \eqref{Onepatch}.

First we focus on the persistence conditions for prey $x_1$. Model \eqref{2patch} is reduced to a single patch model \eqref{Onepatch} when $x_1(0)=0$, i.e., we have $x_1=y_1=0$. Notice that
$$\frac{dx_1}{x_1dt}\Big\vert_{x1=0,y_1=0} = \left(1 - \frac{x_1}{K_1}\right) - \frac{a_1  y_1}{1 + x_1} \Big\vert_{x1=0,y_1=0} =1>0.$$
According to Theorem 2.5 of \cite{hutson1984theorem}, we can conclude that prey $x_1$ is persistent. Similarly, we can show that prey $x_1$ is persistent for all $r>0$.

Since both $x_1$ and $x_2$ are persistent, then we can conclude  that Model \eqref{2patch} is attracted to a subcompact set $C_s$ of $C$ that excludes $E_{0000}, \,E_{00K_20}\,\mbox{ and }\,E_{00K_20}$. Therefore, we can restrict the dynamics of Model \eqref{2patch} on the compact set $C_s$. Now we focus on the persistence conditions for predator $y_1$. According to Theorem \ref{th1:dp}, if $y_1=0$, Model \eqref{2patch} is reduced to the two uncoupled models \eqref{yj0}.
In this case, according to both Theorem \ref{th1:dp} and \ref{th2:be}, the omega limit sets of \eqref{2patch} on the compact set $C_s$ are $E_{K_1000},\,,\,E_{K_10K_20},\,E_{K_10\mu_2\nu_2}$ if $\frac{K_2-1}{2}<\mu_2<K_2$ while they are $E_{K_1000},\,E_{K_10K_20}$ if $\mu_2>K_2$. Now we consider the following two cases:
\begin{enumerate}
\item If $\mu_2>K_2$, according to Theorem 2.5 of  \cite{hutson1984theorem}, we can conclude that predator $y_1$ is persistent if all of the following equations are strictly positive:
$$\begin{array}{lcl}
%\frac{dy_1}{y_1dt} \Big\vert_{E_{0000}}& =& \left[ \frac{a_1 x_1 }{1 + x_1} -d_1+ \rho_1\left(   \frac{a_1 x_1 y_1}{1 + x_1}y_2- \frac{a_2 x_2 y_2}{1 + x_2}\right)\right] \Big\vert_{E_{0000}}=\\
\frac{dy_1}{y_1dt} \Big\vert_{E_{K_1000}}& =& \left[ \frac{a_1 x_1 }{1 + x_1} -d_1+ \rho_1\left(   \frac{a_1 x_1 y_2}{1 + x_1}- \frac{a_2 x_2 y_2}{1 + x_2}\right)\right] \Big\vert_{E_{K_1000}}=\frac{a_1 K_1 }{1 + K_1} -d_1\\
%\frac{dy_1}{y_1dt} \Big\vert_{E_{00K_20}}& =& \left[ \frac{a_1 x_1 }{1 + x_1} -d_1+ \rho_1\left(   \frac{a_1 x_1 y_1}{1 + x_1}y_2- \frac{a_2 x_2 y_2}{1 + x_2}\right)\right] \Big\vert_{E_{00K_20}}=\\
\frac{dy_1}{y_1dt} \Big\vert_{E_{K_10K_20}}& =& \left[ \frac{a_1 x_1 }{1 + x_1} -d_1+ \rho_1\left(   \frac{a_1 x_1 y_2}{1 + x_1}- \frac{a_2 x_2 y_2}{1 + x_2}\right)\right] \Big\vert_{E_{K_10K_20}}=\frac{a_1 K_1 }{1 + K_1} -d_1
\end{array}.$$
Since $\frac{a_1 K_1 }{1 + K_1} -d_1>0 \Leftrightarrow \mu_1<K_1$, therefore, we can conclude that predator $y_1$ is persistent if $\mu_1<K_1$ and $\mu_2>K_2$.
\item If $\frac{K_2-1}{2}<\mu_2<K_2$, according to Theorem 2.5 in \cite{hutson1984theorem} and discussions above, we can conclude that predator $y_1$ is persistent if $\mu_1<K_1$ and the following equation is strictly positive:
$$\begin{array}{lcl}
%\frac{dy_1}{y_1dt} \Big\vert_{E_{0000}}& =& \left[ \frac{a_1 x_1 }{1 + x_1} -d_1+ \rho_1\left(   \frac{a_1 x_1 y_1}{1 + x_1}y_2- \frac{a_2 x_2 y_2}{1 + x_2}\right)\right] \Big\vert_{E_{0000}}=\\
%\frac{dy_1}{y_1dt} \Big\vert_{E_{K_1000}}& =& \left[ \frac{a_1 x_1 }{1 + x_1} -d_1+ \rho_1\left(   \frac{a_1 x_1 y_1}{1 + x_1}y_2- \frac{a_2 x_2 y_2}{1 + x_2}\right)\right] \Big\vert_{E_{K_1000}}=\frac{a_1 K_1 }{1 + K_1} -d_1\\
%\frac{dy_1}{y_1dt} \Big\vert_{E_{00K_20}}& =& \left[ \frac{a_1 x_1 }{1 + x_1} -d_1+ \rho_1\left(   \frac{a_1 x_1 y_1}{1 + x_1}y_2- \frac{a_2 x_2 y_2}{1 + x_2}\right)\right] \Big\vert_{E_{00K_20}}=\\
\frac{dy_1}{y_1dt} \Big\vert_{E_{K_10\mu_2\nu_2}}& =& \left[ \frac{a_1 x_1 }{1 + x_1} -d_1+ \rho_1\left(   \frac{a_1 x_1 y_2}{1 + x_1}- \frac{a_2 x_2 y_2}{1 + x_2}\right)\right] \Big\vert_{E_{K_10\mu_2\nu_2}}=\frac{a_1 K_1 }{1 + K_1} -d_1+ \rho_1\left(   \frac{a_1 K_1 \nu_2}{1 + K_1}- \frac{a_2 \mu_2 \nu_2}{1 + \nu_2}\right)\\
&=&\frac{a_1 K_1 }{1 + K_1} -d_1+ \rho_1 \nu_2\left(   \frac{a_1 K_1}{1 + K_1}- \frac{a_2 \mu_2}{1 + \mu_2}\right)=\frac{a_1 K_1 }{1 + K_1} -d_1+ \rho_1 \nu_2\left(   \frac{a_1 K_1}{1 + K_1}-d_2\right)>0
%\frac{dy_1}{y_1dt} \Big\vert_{E_{00K_20}}& =& \left[ \frac{a_1 x_1 }{1 + x_1} -d_1+ \rho_1\left(   \frac{a_1 x_1 y_1}{1 + x_1}y_2- \frac{a_2 x_2 y_2}{1 + x_2}\right)\right] \Big\vert_{E_{00K_20}}=\\
\end{array}.$$
According to the proof of Theorem \ref{th2:be}, we can see that sufficient condition that $\frac{dy_1}{y_1dt} \Big\vert_{E_{K_10\mu_2\nu_2}}>0$ holds is the same as sufficient condition for the boundary equilibrium $E_{K_10\mu_2\nu_2}$ being unstable when $\mu_1<K_1$. Therefore, we can conclude that predator $y_1$ is persistent if one of the following inequalities hold
 \begin{enumerate}
\item $\mu_j<K_j, \mu_i>K_i$. Or
\item $\frac{K_i-1}{2}<\mu_i<K_i \mbox{ and } K_j>\max\Big\{\mu_j,\frac{d_i}{a_j-d_i}\Big\}$. Or
\item $\frac{K_i-1}{2}<\mu_i<K_i,\,\,\mu_j<K_j<\frac{d_i}{a_j-d_i} \mbox{ and }\rho_j<\frac{K_j(a_j-d_j)-d_j}{\nu_i\left[d_i-K_j(a_j-d_i)\right]}$. % where $i=1,j=2$ or $i=2,j=1$.
\end{enumerate} where $i=1,j=2$ or $i=2,j=1$.one of the following inequalities hold.
\end{enumerate}

Based on the discussion above, we can conclude that the statement of Theorem \ref{th4:p} holds.

\end{proof}

%%%%%%%%%%%%%%%%%%%%%%%%%%%%%%%%%%%%%%%%%%%%%%%%%%%%%%%%%%%%%%%
\subsection*{Proof of Theorem \ref{th5:perm}}
\begin{proof}If $\frac{K_j-1}{2}<\mu_j<K_j, \mu_i>K_i$, then according to Theorem \ref{th4:p}, we can conclude that prey $x_i$ for both $i=1,2$ and predator $y_j$ is persistent. This implies that Model \eqref{2patch} is permanent if predator $y_i$ is persistent.
Since  $\frac{K_j-1}{2}<\mu_j<K_j$, then Theorem \ref{th2:be} indicates that the omega limit set of Model \eqref{2patch} when $y_i=0$ is $E_{\mu_1\nu_1K_20}$ when $i=2,j=1$ while its omega limit set is $E_{K_10\mu_2\nu_2}$ when $i=2,j=1$. Now let $i=1,j=2$, then according to Theorem 2.5 of  \cite{hutson1984theorem}, we can conclude that predator $y_1$ is persistent if the following equation is strictly positive:
$$\begin{array}{lcl}
\frac{dy_1}{y_1dt} \Big\vert_{E_{K_10\mu_2\nu_2}}& =& \left[ \frac{a_1 x_1 }{1 + x_1} -d_1+ \rho_1\left(   \frac{a_1 x_1 y_2}{1 + x_1}- \frac{a_2 x_2 y_2}{1 + x_2}\right)\right] \Big\vert_{E_{K_10\mu_2\nu_2}}=\frac{a_1 K_1 }{1 + K_1} -d_1+ \rho_1\left(   \frac{a_1 K_1 \nu_2}{1 + K_1}- \frac{a_2 \mu_2 \nu_2}{1 + \nu_2}\right)\\
&=&\frac{a_1 K_1 }{1 + K_1} -d_1+ \rho_1 \nu_2\left(   \frac{a_1 K_1}{1 + K_1}- \frac{a_2 \mu_2}{1 + \mu_2}\right)=\frac{a_1 K_1 }{1 + K_1} -d_1+ \rho_1 \nu_2\left(   \frac{a_1 K_1}{1 + K_1}-d_2\right)>0
\end{array}.$$Since $\mu_i>K_i\Leftrightarrow \frac{a_1 K_1 }{1 + K_1} -d_1<0$, therefore, sufficient conditions for $\frac{dy_1}{y_1dt} \Big\vert_{E_{K_10\mu_2\nu_2}}>0$ is
$$\frac{a_1 K_1}{1 + K_1}-d_2>0 \Leftrightarrow K_1>\frac{d_2}{a_1-d_2} \mbox{ and } \rho_1>\frac{d_1-K_1(a_1-d_1)}{v_2\left[K_1(a_1-d_2)-d_2\right]}.$$
Similarly, we can show that predator $y_2$ is persistent when $i=2,j=1$. Therefore,  Model \eqref{2patch} is permanent if the following inequalities hold for either $i=2, j=1$ or $i=1,j=2$,
$$\frac{K_j-1}{2}<\mu_j<K_j, \mu_i>K_i,\,0<\frac{d_j}{a_i-d_j}<K_i<\mu_i \mbox{ and } \rho_i>\frac{d_i-K_i(a_i-d_i)}{\nu_j\left[K_i(a_i-d_j)-d_j\right]}.$$

According to Theorem \ref{th4:p}, we can conclude that prey $x_i$ for both $i=1,2$ and predator $y_i$ is persistent if the following inequalities hold
$$\frac{K_i-1}{2}<\mu_i<K_i,\, \mu_j>\frac{K_j-1}{2} \mbox{ and } K_j>\max\Big\{\mu_j,\frac{d_i}{a_j-d_i}\Big\}.$$
Therefore, Model \eqref{2patch} is permanent if the above inequalities hold for both $i=1,j=2$ and $i=2,j=1$.
On the other hand, predator $y_j$ is persistent if the following inequalities hold
$$\frac{K_i-1}{2}<\mu_i<K_i,\,\frac{K_j-1}{2}<\mu_j<K_j<\frac{d_i}{a_j-d_i} \mbox{ and }\rho_j<\frac{K_j(a_j-d_j)-d_j}{\nu_i\left[d_i-K_j(a_j-d_i)\right]}.$$
Therefore, both predator $y_i$ and $y_j$ are persistent if the following inequalities hold for either $i=1,j=2$ or $i=2,j=1$,
$$\frac{K_i-1}{2}<\mu_i<K_i,\,K_i>\max\Big\{\mu_i,\frac{d_j}{a_i-d_j}\Big\},\,\frac{K_j-1}{2}<\mu_j<K_j<\frac{d_i}{a_j-d_i} \mbox{ and }\rho_j<\frac{K_j(a_j-d_j)-d_j}{\nu_i\left[d_i-K_j(a_j-d_i)\right]}.$$

Based on the discussion above, we can conclude that the statement of Theorem \ref{th5:perm} holds.

\end{proof}
%%%%%%%%%%%%%%%%%%%%%%%%%%%%%%%%%%%%%%%%%%%%%%%%%%%%%%%%%%%%%%%
\subsection*{Proof of Theorem \ref{th6:interior}}
\begin{proof}If $\mu_i>K_i$ for both $i=1,2$, then Model \eqref{2patch} has global stability at $(K_1,0,K_2,0)$ according to Theorem \ref{th3:gb}. This implies that Model \eqref{2patch} has no interior equilibrium when $\mu_i>K_i$ for both $i=1,2$.

The interior equilibrium $(x^*_1,y^*_1,x^*_2,y_2^*)$ is determined by the positive intersections of the nullclines \eqref{interior-eq3}
$x_1=F(x_2)=\frac{f_t(x_2)}{f_b(x_2)}$ and $x_2=G(x_1)=\frac{g_t(x_1)}{g_b(x_1)}$ where
$$f_t(x_2)=a_2\left[r_2\rho_1x_2\left(K_2-x_2\right)+K_2d_1\right],\,f_b(x_2)=r_2\rho_1x_2\left(K_2a_1-K_2a_2-a_1\right)-r_2\rho_1x_2^2(a_1-a_2)+K_2(a_1r_2\rho_1+a_1a_2-a_2d_1)$$ and
$$g_t(x_1)=a_1\left[r_1\rho_2x_1\left(K_1-x_1\right)+K_1d_2\right],\, g_b(x_1)=r_1\rho_2x_1\left(K_1a_2-K_1a_1-a_2\right)-r_1\rho_2x_1^2(a_2-a_1)+K_1(a_2r_1\rho_2+a_1a_2-a_1d_2).$$
Notice that the nullclines  $x_1=F(x_2)=\frac{f_t(x_2)}{f_b(x_2)}$ and $x_2=G(x_1)=\frac{g_t(x_1)}{g_b(x_1)}$ has the following properties:
\begin{enumerate}
\item $F(0)=\frac{f_t(0)}{f_b(0)}=\frac{a_2K_2d_1}{K_2(a_1r_2\rho_1+a_1a_2-a_2d_1)}=\frac{a_2d_1}{a_1r_2\rho_1+a_1a_2-a_2d_1}$ and $F(K_2)=\frac{f_t(K_2)}{f_b(K_2)}=\frac{d_1}{a_1-d_1}=\mu_1$.
\item $f_t(x_2)=a_2\left[r_2\rho_1x_2\left(K_2-x_2\right)+K_2d_1\right]\geq a_2K_2 d_1>0 \mbox{ for } x_2\in [0, K_2]$ and
$$f_b(x_2)\big\vert_{a_1=a_2=a}=a\left[r_2\rho_1(K_2-x_2)+K_2(a-d_1)\right].$$
\item $G(0)=\frac{g_t(0)}{g_b(0)}=\frac{a_1K_1d_2}{K_1(a_2r_1\rho_2+a_1a_2-a_1d_2)}=\frac{a_1d_2}{a_2r_1\rho_2+a_1a_2-a_1d_2}$ and $G(K_1)=\frac{g_t(K_1)}{g_b(K_1)}=\frac{d_2}{a_2-d_2}=\mu_2$.
\item $g_t(x_1)=a_1\left[r_1\rho_2x_1\left(K_1-x_1\right)+K_1d_2\right]\geq a_1 K_1d_1>0\mbox{ for } x_1\in [0, K_1]$ and
$$g_b(x_1)\big\vert_{a_1=a_2=a}=a\left[r_1\rho_2(K_1-x_1)+K_1(a-d_2)\right].$$
\end{enumerate}
According to Theorem \ref{th1:dp}, we know that population of prey $x_i$ for $i=1,2$ has the following properties:
$$\limsup_{t\rightarrow\infty} x_i(t)\leq K_i.$$
Thus, we can restrict the function $F(x_2)$ on the domain of $[0, K_2]$ and $G(x_1)$ on the domain of $[0, K_1]$. Since $f_t(x_2)\geq a_2K_2 d_1>0 \mbox{ for } x_2\in [0, K_2]$ and $g_t(x_1)\geq a_1 K_1d_1>0\mbox{ for } x_1\in [0, K_1]$, thus, Model \eqref{2patch} has no interior equilibrium if $f_b(x_2)<0$ \mbox{ for } $x_2\in [0, K_2]$ or $g_b(x_1)<0$ \mbox{ for } $x_1\in [0, K_1]$. \\

Now we assume that $a_1>a_2$, then we have $f_b(x_2)<0$ \mbox{ for } $x_2\in [0, K_2]$ if
{\small$$r_2^2\rho_1^2\left(K_2a_1-K_2a_2-a_1\right)^2<4K_2r_2\rho_1(a_2-a_1)(a_1r_2\rho_1+a_1a_2-a_2d_1)\Leftrightarrow \rho_1<\frac{4K_2a_2(a_1-a_2)(d_1-a_1)}{r_2\left(K_2a_1-K_2a_2+a_1\right)^2}$$} while $f_b(x_2)>0$ \mbox{ for } $x_2\in [0, K_2]$ if $a_1>d_1$ since
$$f_b(0)=K_2(a_1r_2\rho_1+a_2(a_1-d_1))>0 \mbox{ and } f_b(K_2)=a_2K_2(a_1-d_1)>0.$$
And $g_b(x_1)>0$ \mbox{ for } $x_1\in [0, K_1]$ if
{\small$$r_1^2\rho_2^2\left(K_1a_2-K_1a_1-a_2\right)^2<4K_1r_1\rho_2(a_1-a_2)(a_2r_1\rho_2+a_1a_2-a_1d_2)\Leftrightarrow \rho_2<\frac{4K_1a_1(a_1-a_2)(a_2-d_2)}{r_1\left(K_1a_2-K_1a_1+a_2\right)^2}.$$}

Similar cases can be made for $a_1<a_2$, therefore we can conclude that Model \eqref{2patch} has no interior equilibrium if  either $$ a_1>a_2, \,\rho_1<\frac{4K_2a_2(a_1-a_2)(d_1-a_1)}{r_2\left(K_2a_1-K_2a_2+a_1\right)^2}$$%a_1>\max\{a_2,d_1\},
%, \,\rho_2>\frac{4K_1a_1(a_1-a_2)(a_2-d_2)}{r_1\left(K_1a_2-K_1a_1+a_2\right)^2}$$
or %\rho_1>\frac{4K_2a_2(a_2-a_1)(a_1-d_1)}{r_2\left(K_2a_1-K_2a_2+a_1\right)^2},
$$a_2>a_1,\,\rho_2<\frac{4K_1a_1(a_2-a_1)(d_2-a_2)}{r_1\left(K_1a_2-K_1a_1+a_2\right)^2}$$hold.

Now we focus on sufficient conditions lead to both $f_b(x_2)>0$ \mbox{ for } $x_2\in [0, K_2]$ and $g_b(x_1)>0$ \mbox{ for } $x_1\in [0, K_1]$. Assume that $a_1>a_2$ and $a_i>d_i$ for both $i=1,2$. Notice that $f_b(0)=K_2(a_1r_2\rho_1+a_2(a_1-d_1))> f_b(K_2)=a_2K_2(a_1-d_1)>0$ with
$$f_b(x_2)=r_2\rho_1x_2\left(K_2a_1-K_2a_2-a_1\right)-r_2\rho_1x_2^2(a_1-a_2)+K_2(a_1r_2\rho_1+a_1a_2-a_2d_1).$$
Therefore, we have $f_b(x_2)>f_b(K_2)$ \mbox{ for all } $x_2\in [0, K_2]$. Since $g_b(x_1)$ is a degree 2  polynomial with the positive coefficient in the degree 2 and $g_b(0)=K_1(a_2r_1\rho_2+a_1(a_2-d_2))> g_b(K_1)=a_1K_1(a_2-d_2)>0$. Therefore, we have $g_b(x_1)>0$ \mbox{ for all } $x_1\in [0, K_1]$ if the following conditions hold
$$\rho_2<\frac{4K_1a_1(a_2-a_1)(d_2-a_2)}{r_1\left(K_1a_2-K_1a_1+a_2\right)^2}.$$
The discussion so far also indicates that we have both $f_b(x_2)>0$ \mbox{ for } $x_2\in [0, K_2]$ and $g_b(x_1)>0$ \mbox{ for } $x_1\in [0, K_1]$ if the following inequalities hold for $i=1,j=2$ or $i=2,j=1$
$$a_i>\max\{a_j,d_i\},\, \rho_j<\frac{4K_ia_i(a_j-a_i)(d_j-a_j)}{r_i\left(K_ia_j-K_ia_i+a_j\right)^2} .$$
Now assume that these conditions hold, then we have $F(x_2)$ and $G(x_1)$ are positive on their restricted domain. By algebraic calculations, if $a_i>\max\{d_1,d_2\}$ for both $i=1,2$, then both $F(x_2)$ and $G(x_1)$ have its unique critical points $x^c_i,i=1,2$ in their restricted domain where
$$x^c_1=\frac{K_1\left(r_1\rho_2+a_1-d_2-\sqrt{(a_1-d_2)(r_1\rho_2+a_1-d_2)}\right)}{r_1\rho_2} \in (0,K_1)$$ \mbox{ and }$$x^c_2=\frac{K_2\left(r_2\rho_1+a_2-d_1-\sqrt{(a_2-d_1)(r_2\rho_1+a_2-d_1)}\right)}{r_2\rho_1} \in (0,K_2).$$
If $F(x_2^c)<K_1 \mbox{ and } G(x_1^c)<K_2$, then we can conclude that both maps $x_1=F(x_2)$ and $x_2=G(x_1)$ are unimode and the skew product of $F\times G$ maps $[0,K_2]\times[0,K_1]$ to its compact subset. Since both $F$ and $G$ are continuous and differentiable, therefore, $x_1=F(x_2)$ and $x_2=G(x_1)$ has at least one positive intersection for $x_2\in [0,K_2], \, x_1\in [0,K_1]$. Now we focus on sufficient condition that leads to $F(x_2^c)<K_1 \mbox{ and } G(x_1^c)<K_2$ when $a_1>a_2$.
Since $$\max_{0\leq x_2\leq K_2}\{F(x_2)\}=F(x_2^c)\leq\frac{\max_{0\leq x_2\leq K_2}\{f_t(x_2)\}}{\min_{0\leq x_2\leq K_2}\{f_b(x_2)\}}=\frac{f_b(K_2/2)}{f_b(K_2)}=\frac{K_2r_2\rho_1+4d_1}{4(a_1-d_1)}$$ and
{\small$$\max_{0\leq x_1\leq K_1}\{G(x_1)\}=G(x_1^c)\leq\frac{\max_{0\leq x_1\leq K_1}\{g_t(x_1)\}}{\min_{0\leq x_1\leq K_1}\{g_b(x_1)\}}=\frac{g_t(K_1/2)}{g_b\left(\frac{K_2a_1-K_2a_2-a_1}{2(a_1-a_2)}\right)}=\frac{a_1K_1(r_1K_1\rho_2/4+d_2)}{K_1a_1(a_2-d_2)-\frac{r_1\rho_2(K_1a_1-K_1a_2-a_2)^2}{4(a_1-a_2)}},$$}
therefore, we have $F(x_2^c)<K_1 \mbox{ and } G(x_1^c)<K_2$ when $a_1>a_2$ if the following inequalities hold
$$\frac{K_2r_2\rho_1+4d_1}{4(a_1-d_1)}\leq K_1 \Leftrightarrow \rho_1\leq \frac{4(K_1a_1-K_1d_1-d_1)}{K_2r_2}.$$
and $$\begin{array}{lcl}
K_2&\geq&\frac{a_1K_1(r_1K_1\rho_2/4+d_2)}{K_1a_1(a_2-d_2)-\frac{r_1\rho_2(K_1a_1-K_1a_2-a_2)^2}{4(a_1-a_2)}}\\
&\Leftrightarrow&\\
\rho_2&\leq&\frac{4K_1a_1(K_2a_2-K_2d_2-d_2)}{a_1r_1K_1^2+r_2K_2(K_1a_1-K_1a_2-a_2)^2}.
\end{array}$$Therefore, we can conclude that Model \eqref{2patch} has at least one interior equilibrium $(x_1^*,y_1^*,x_2^*,y_2^*)$ if the following inequalities hold
$$a_i>\max\{a_j,d_1,d_2\}, a_j>\max\{d_1,d_2\}, \,\rho_i\leq \frac{4(K_ia_i-K_id_i-d_i)}{K_jr_j}$$ and
$$\rho_j<\min\Big\{\frac{4K_ia_i(a_j-a_i)(d_j-a_j)}{r_i\left(K_ia_j-K_ia_i+a_j\right)^2},\frac{4K_ia_i(K_ja_j-K_jd_j-d_j)}{a_ir_iK_i^2+r_jK_j(K_ia_i-K_ia_j-a_j)^2}\Big\}.$$
In addition, since both $F(x_2)$ and $G(x_1)$ are unimode maps in their domain with unique local maximum, thus, we have
$$x_1=F(x_2)\geq F(0)=\frac{a_2d_1}{a_1r_2\rho_1+a_1a_2-a_2d_1}\mbox{ and }x_2=G(x_1)\geq G(0)=\frac{a_1d_2}{a_2r_1\rho_2+a_1a_2-a_1d_2}.$$Therefore, we have  $\frac{a_id_j}{a_jr_i\rho_j+a_ia_j-a_id_j}<x_j^*<K_j$ for both $i=1,j=2$ and $i=2,j=1$.\\

Now assume that $a_1=a_2=a$, then both $f_b$ and $g_b$ are reduced to linear decreasing functions, i.e.,
$$f_b(x_2)\big\vert_{a_1=a_2=a}=a\left[r_2\rho_1(K_2-x_2)+K_2(a-d_1)\right] \mbox{ and }g_b(x_1)\big\vert_{a_1=a_2=a}=a\left[r_1\rho_2(K_1-x_1)+K_1(a-d_2)\right]$$which indicates that  $f_b(x_2)<0$ \mbox{ for } $x_2\in [0, K_2]$ if $d_1>a+r_2\rho_1$ and $g_b(x_1)<0$ \mbox{ for } $x_1\in [0, K_1]$ if $d_2>a+r_1\rho_2$. Therefore, if $a_1=a_2=a$ and either $d_1>a+r_2\rho_1$ or $d_2>a+r_1\rho_2$ holds, then Model \eqref{2patch} has no interior equilibrium.
On the other hand, both $f_b(x_2)>0$ \mbox{ for } $x_2\in [0, K_2]$ and $g_b(x_1)>0$ \mbox{ for } $x_1\in [0, K_1]$ if $a>\max\{d_1,d_2\}$. Then apply the discussions for the case $a_1\neq a_2$, then we can conclude that Model \eqref{2patch} has at least one interior equilibrium if  $$a>\max\{d_1,d_2\},\,F(x^c_2)<K_1,\,\mbox{ and } G(x^c_1)<K_2.$$
Applying the similar arguments for the case $a_i>a_j$, we have
$$\max_{0\leq x_2\leq K_2}\{F(x_2)\}=F(x_2^c)\leq\frac{\max_{0\leq x_2\leq K_2}\{f_t(x_2)\}}{\min_{0\leq x_2\leq K_2}\{f_b(x_2)\}}=\frac{f_b(K_2/2)}{f_b(K_2)}=\frac{K_2r_2\rho_1+4d_1}{4(a_1-d_1)}$$ and
{\small$$\max_{0\leq x_1\leq K_1}\{G(x_1)\}=G(x_1^c)\leq\frac{\max_{0\leq x_1\leq K_1}\{g_t(x_1)\}}{\min_{0\leq x_1\leq K_1}\{g_b(x_1)\}}=\frac{g_t(K_1/2)}{g_b\left(K_1\right)}=\frac{K_1r_1\rho_2+4d_2}{4(a_2-d_2)}.$$}
Therefore, we can conclude that Model \eqref{2patch} has at least one interior equilibrium $(x_1^*,y_1^*,x_2^*,y_2^*)$ if the following inequalities hold
 $$a_1=a_2=a>\max\{d_1,d_2\}, \rho_i<\frac{4(K_ia-K_id_i-d_i)}{K_jr_j}$$ for both $i=1,j=2$ and $i=2,j=1$. In addition, $\frac{d_j}{r_i\rho_j+a-d_j}<x_j^*<K_j$ hold for both $i=1,j=2$ and $i=2,j=1$.

%\begin{enumerate}
%\item If $a_1>a_2$, then
%\end{enumerate}
\end{proof}
%%%%%%%%%%%%%%%%%%%%%%%%%%%%%%%%%%%%%%%%%%%%%%%%%%%%%%%%%%%%%%%
\subsection*{Proof of Theorem \ref{th7:interior}}
\begin{proof}
Suppose that $\min\{a_1,a_2\}>d$, then we have $\mu_i=\frac{d}{a_i-d}, \nu_i=\frac{r_i(K_i-\mu_i)(1+\mu_i)}{a_iK_i}$ where $(\mu_i,\nu_i)$ is the unique interior equilibrium of the single patch model \eqref{Onepatch} in the absence of the dispersal in predator for both $i=1,2$. Now recall from the null clines \eqref{interior-eq2}, we have
$$\begin{array}{lcl}
x_i(x_j)&=&\frac{\rho_i q_j(x_j)p_j(x_j)+d_i}{a_i(1+\rho_i q_j(x_j))-(\rho_i q_j(x_j)p_j(x_j)+d_i)}=\frac{\rho_i q_j(x_j)p_j(x_j)+d_i}{\rho_i q_j(x_j)\left[a_i-p_j(x_j)\right]+a_i-d_i}%\\\\%=\frac{f_t(x_2)}{f_b(x_2)}=F(x_2)\\\\
%x_2&=&G(x_1)=\frac{a_1\left[r_1\rho_2x_1\left(K_1-x_1\right)+K_1d_2\right]}{r_1\rho_2x_1\left(K_1a_2-K_1a_1-a_2\right)-r_1\rho_2x_1^2(a_2-a_1)+K_1(a_2r_1\rho_2+a_1a_2-a_1d_2)}%=\frac{g_t(x_1)}{g_b(x_1)}=G(x_1)
\end{array}$$which indicates that
$$x_i(\mu_j)=%\frac{\rho_i q_j(\mu_j)p_j(\mu_j)+d}{a_i(1+\rho_i q_j(\mu_j))-(\rho_i q_j(\mu_j)p_j(\mu_j)+d)}=
\frac{\rho_i q_j(\mu_j)p_j(\mu_j)+d}{\rho_i q_j(\mu_j)\left[a_i-p_j(\mu_j)\right]+a_i-d}=\frac{d}{a_i-d}=\mu_i.$$
This implies that $x_i=\mu_i$ for $i=1,2$ is the positive solutions of the null clines \eqref{interior-eq2}. Therefore, we can solve that $(\mu_1,\nu_1,\mu_2,\nu_2)$ is an interior solution of Model \eqref{2patch} if $\min\{a_1,a_2\}>d=d_1=d_2$.
By substituting the equilibrium $(\mu_1,\nu_1,\mu_2,\nu_2)$ into the Jacobian matrix \eqref{JE}, we obtain its characteristic equation as follows:
$$H(\lambda)=\lambda^4+(\alpha_1+\alpha_2)\lambda^3+[\alpha_1\alpha_2+d(\beta_1+\beta_2)]\lambda^2+d(\alpha_1\beta_2+\alpha_2\beta_1)\lambda+d^2(\beta_1\beta_2-\gamma_1\gamma_2=0$$where
$$\begin{array}{lcl}
\alpha_i&=&-\frac{r_i\mu_i(K_ia_i-K_id-a_i-d)}{K_ia_i}\Rightarrow [\alpha_i>0\Leftrightarrow  K_ia_1-K_id-a_i-d<0\Leftrightarrow \frac{K_i-1}{2}<\mu_i<K_i]\\
\beta_i&=&\frac{\nu_i(\nu_j\rho_i+1)(a_i-d)^2}{a_i}>0\\
\gamma_i&=&\frac{\rho_i\nu_i\nu_j(a_j-d)^2}{a_j}>0\\
\beta_1\beta_2-\gamma_1\gamma_2&=&\frac{\nu_1\nu_2(a_1-d)^(a_2-d)^2(\nu_1\rho_2+\nu_2\rho_1+1)}{a_1a_2}>0
\end{array}.$$
Then the real parts of the solutions of $H(\lambda)$ are all negative if $\alpha_1+\alpha_2>0$ while the solutions of $H(\lambda)$ has positive if $\alpha_1+\alpha_2<0$. Notice that the single patch $i$ \eqref{Onepatch} has global stability at $(\mu_i,\nu_i)$ if $\alpha_i>0\Leftrightarrow \frac{K_i-1}{2}<\mu_i<K_i.$ Therefore, the interior equilibrium is locally asymptotically stable of $\frac{K_i-1}{2}<\mu_i<K_i$ for both $i=1,2$ while it is unstable if
$$\frac{\mu_1(K_1a_1-K_1d-a_1d)}{K_1a_1}+\frac{r\mu_2(K_2a_2-K_2d-a_2-d)}{K_2a_2}>0.$$
Assume that $\alpha_1\alpha_2<0$\mbox{ and } $\alpha_1+\alpha_2>0$. Then the real parts of all solutions of $H(\lambda)$ can be still negative if $\alpha_1\alpha_2+d(\beta_1+\beta_2)>0$ and $d(\alpha_1\beta_2+\alpha_2\beta_1)$. By algebraic calculations, we can conclude that if $\alpha_i<0$ and the dispersal of predator $y_i$ is large enough, then the interior equilibrium $(\mu_1,\nu_1,\mu_2,\nu_2)$ can still be locally asymptotically stable, where $\rho_i$ should satisfy the following condition:
{\small$$\rho_i>\max\Big\{\frac{-\nu_j-r_j\mu_i\mu_j(K_ia_i-K_id-a_i-d)(K_ja_j-K_jd-a_j-d)}{(K_iK_ja_j\nu_jd\nu_i(a_i-d)^2)},\frac{-\frac{\mu_i\nu_jK_j(\nu_i\rho_j+1)(a_j-d)^2(K_ia_i-K_id-a_i-d)}{r_j\mu_j\nu_iK_i(a_i-d)^2(K_ja_j-K_jd-a_j-d)}-1}{\nu_j}\Big\} .$$}
Now if $a_1=a_2=a, r_1=r_2=1, K_1=K_2=K, d_1=d_2=d$. The discussions above implies that Model \eqref{2patch} has the same stability at $E^i=(\mu,\nu,\mu,\nu)$ as the stability of the single patch model \eqref{Onepatch} at $(\mu, \nu)$ where $\mu=\frac{d}{a-d}$ and $\nu=\frac{(K-\mu)(1+\mu)}{aK}$. Therefore, $E^i$ is locally asymptotically stable if $\frac{K-1}{2}<\mu<K$ while it is unstable if $\mu>\frac{K-1}{2}$.

Now we should show that Model \eqref{2patch} has the unique $E^i=(\mu,\nu,\mu,\nu)$ whenever $a>d$. Notice that $F(x_2)$ and $G(x_1)$ have the following properties in the symmetric case (i.e., $a_1=a_2=a, r_1=r_2=1, K_1=K_2=K, d_1=d_2=d$):
\begin{enumerate}
\item $F(x_2)=\frac{\rho_1x_2\left(K-x_2\right)+Kd}{\rho_1(K-x_2)+aK(1-d)}$ with $F(0)=\frac{d}{\rho_1+a -d}$ and $F(\mu)=F(K)=\mu$. In addition, $F(x_2)>0$ for all $x_2\in [0,K]$ and $F(x_2)$ has a unique critical point $x^c_2$ for $x_2\in [0,K]$ where
$$\mu<x^c_2=\frac{K\left(\rho_1+a-d-\sqrt{(a-d)(\rho_1+a-d)}\right)}{\rho_1} \in (0,K).$$
\item $G(x_1)=\frac{\rho_2x_1\left(K-x_1\right)+Kd}{\rho_2(K-x_1)+aK(1-d)}$ with $G(0)=\frac{d}{\rho_2+a -d}$ and $G(\mu)=G(K)=\mu$. In addition, $G(x_1)>0$ for all $x_1\in [0,K]$ and $G(x_1)$ has a unique critical point $x^c_1$ for $x_1\in [0,K]$ where
$$\mu<x^c_1=\frac{K\left(\rho_2+a-d-\sqrt{(a-d)(\rho_2+a-d)}\right)}{\rho_2} \in (0,K).$$
\end{enumerate}
The discussions above indicate that both $F(x_2)$ and $G(x_1)$ are unimode maps with a unique interception at $x_1=x_2=\mu$.

%we discuss the following three cases:
\end{proof}
%%%%%%%%%%%%%%%%%%%%%%%%%%%%%%%%%%%%%%%%%%%%%%%%%%%%%%%%%%%%%%%
\subsection*{Proof of Theorem \ref{th8:ds}}
\begin{proof}Proof of Item 1 can be obtained by adopting the proof provided in Theorem \ref{th1:dp}. We omit details.\\

The stability of $E_{0000},\,E_{K_1000},\,E_{00K_20},\,E_{K_10K_20}$ can be obtained from eigenvalues of the Jacobian matrix of Model \eqref{2patch2} evaluated at these equilibria through simple algebraic calculations. We omit details. But we will return to the local stability of $E^b_i$ when we prove Item 4.\\

Item 3(a):  If $x_i=0$, then we have $\lim_{t\rightarrow\infty} y_1(t)=\lim_{t\rightarrow\infty} y_2(t)=0$. This implies that the omega limit set of Model \eqref{2patch2s} is $y_1=y_2=0$. Therefore, prey $x_i$ persists by applying Theorem 2.5 of \cite{hutson1984theorem} since
$$\frac{dx_i}{x_idt}\Big\vert_{x_i=0}=r_i>0.$$
Item 3(b): Recall $p_i(x)=\frac{a_i x}{1+x}$ and $q_i(x)=\frac{r_i(K_i-x)(1+x)}{a_iK_i}$ for $i=1,2$. Then we construct the following Lyapunov function
$$\begin{aligned}
V_{ij}(x_i,y_i,y_j) =(\rho_j +d_j)\int^{x_i}_{K_i} \frac{p_i(\xi)-p_i(K_i)}{p_i(\xi)}d\xi + (\rho_j +d_j)y_i+\rho_iy_j.
\end{aligned}$$If $ \hat{\mu_i}>K_i$, then we have
$$\frac{dV_{ij}(x_i,y_i,y_j)}{dt} =(\rho_j+d_j)\left[(p_i(x_i)-p_i(K_i)\right]q_i(x_i)+y_i(p_i(K_i)-\hat{d}_i)<0$$ since
$p_i(K_i)-\hat{d}_i<0\Leftrightarrow \hat{\mu_i}>K_i.$ Therefore, Model \eqref{2patch2s} has global stability at $(K_i,0,0)$ if $\hat{\mu_i}>K_i$.\\\\
Item 3(c): We construct the following Lyapunov function
$$\begin{aligned}
V_{ij}(x_i,y_i,y_j) =(\rho_j +d_j)\int^{x_i}_{\hat{\mu}_i} \frac{p_i(\xi)-p_i(\hat{\mu}_i)}{p_i(\xi)}d\xi + (\rho_j +d_j)\int^{y_i}_{\hat{\nu}_i} \frac{\eta_i-\hat{\nu}_i}{\eta_i}d\eta_i+\rho_i\int^{y_j}_{\hat{\nu}_j^i} \frac{\eta_j-\hat{\nu}_j^i}{\eta_j}d\eta_j.
\end{aligned}$$If $\frac{K_i-1}{2}< \hat{\mu_i}<K_i$, then we have
$$\frac{dV_{ij}(x_i,y_i,y_j)}{dt} =(\rho_j+d_j)\left[(p_i(x_i)-p_i(\hat{\mu}_i)\right]\left[q_i(x_i)-\hat{\nu}_i\right]-\frac{\rho_i\hat{\nu}_i((\rho_j+d_j)y_j-\rho_jy_i)^2}{(\rho_j+d_j)y_iy_j}<0.$$ Therefore, Model \eqref{2patch2s} has global stability at $(\hat{\mu_i},\hat{\nu_i},\hat{\nu}_j^i)$ if $\frac{K_i-1}{2}<\hat{\mu_i}<K_i.$\\\\

Item 4: If $x_j=0$, then Model \eqref{2patch2} reduces to Model \eqref{2patch2s} who has global stability at $(K_i,0,0)$ if $\hat{\mu_i}>K_i$ while has $(\hat{\mu_i},\hat{\nu_i},\hat{\nu}_j^i)$ if
$\frac{K_i-1}{2}<\hat{\mu_i}<K_i.$
Therefore, by applying Theorem 2.5 of \cite{hutson1984theorem}, prey $x_j$ persists if
$$\frac{dx_j}{x_jdt}\Big\vert_{x_j=0,y_j=0}=r_j>0 \mbox{ when } \hat{\mu_i}>K_i$$
and $$\frac{dx_j}{x_jdt}\Big\vert_{x_j=0,y_j=\hat{\nu}_j^i}=r_j-a_j\hat{\nu}_j^i>0 \mbox{ when } \frac{K_i-1}{2}<\hat{\mu_i}<K_i.$$
The persistence of both prey can be easily obtained from the persistence of one prey.\\
 If $\frac{K_i-1}{2}<\hat{\mu_i}<K_i$ and $r_j-a_j\hat{\nu}_j^i<0$, then $\frac{dx_j}{x_jdt}\Big\vert_{x_j=0,y_j=\hat{\nu}_j^i}=r_j-a_j\hat{\nu}_j^i<0$. According to Theorem 2.18 by \cite{hutson1992permanence}, we can conclude that the boundary equilibrium $(\hat{\mu_i},\hat{\nu_i},0,\hat{\nu}_j^i)$ is locally asymptotically stable. This proves the stability condition of $E^b_i$ of Item 2.
\\\\

Item 5: We construct the following Lyapunov function
$$\begin{aligned}
V(x_i,y_i,x_j,y_j) =(\rho_i +d_i)\int^{x_j}_{\hat{\mu}_j} \frac{p_j(\xi_j)-p_j(\hat{\mu}_j)}{p_j(\xi_j)}d\xi_j + (\rho_i +d_i)\int^{y_j}_{\hat{\nu}_j} \frac{\eta_j-\hat{\nu}_j}{\eta_j}d\eta_j+\rho_i x_i+\int^{y_i}_{\hat{\nu}_i^j} \frac{\eta_i-\hat{\nu}_i^j}{\eta_i}d\eta_i.
\end{aligned}$$
Then we have
{\small$$\begin{array}{lcl}
\frac{dV(x_i,y_i,x_j,y_j)}{dt} &=&(\rho_j+d_j)\left[(p_i(x_i)-p_i(\hat{\mu}_i)\right]\left[q_i(x_i)-\hat{\nu}_i\right]-\frac{\rho_i\hat{\nu}_i((\rho_j+d_j)y_j-\rho_jy_i)^2}{(\rho_j+d_j)y_iy_j}+\rho_jp_i(x_i)\left[q_i(x_i)-\hat{\nu}_i^j\right]\\
&<&(\rho_j+d_j)\left[(p_i(x_i)-p_i(\hat{\mu}_i)\right]\left[q_i(x_i)-\hat{\nu}_i\right]-\frac{\rho_i\hat{\nu}_i((\rho_j+d_j)y_j-\rho_jy_i)^2}{(\rho_j+d_j)y_iy_j}+\rho_jp_i(x_i)\left[q_i(\frac{K_i-1}{2})-\hat{\nu}_i^j\right]\\
&=&(\rho_j+d_j)\left[(p_i(x_i)-p_i(\hat{\mu}_i)\right]\left[q_i(x_i)-\hat{\nu}_i\right]-\frac{\rho_i\hat{\nu}_i((\rho_j+d_j)y_j-\rho_jy_i)^2}{(\rho_j+d_j)y_iy_j}+\rho_jp_i(x_i)\left[\frac{r_i (K_i+1)^2}{4a_iK_i}-\hat{\nu}_i^j\right]\end{array}
.$$}
Therefore, if $\frac{K_j-1}{2}<\hat{\mu_j}<K_j \mbox{ and } \frac{r_i (K_i+1)^2}{4a_iK_i}<\hat{\nu}_j^i$ hold, then we have
$\frac{dV(x_i,y_i,x_j,y_j)}{dt}<0$. This implies that Model \eqref{2patch2s} has global stability at
$(0,\hat{\nu}_i^j,\hat{\mu_j},\hat{\nu_j})$.\\

Item 6: Define $V(y_1,y_2)=\rho_2y_1+\rho_1y_2$, then we have
$$\frac{dV}{dt}=\rho_2(p_1(x_1)-d_1)y_1+\rho_1(p_2(x_2)-d_2)y_2.$$
Notice that $\limsup_{t\rightarrow\infty} x_i(t)\leq K_i$ for both $i=1,2$. Then if $\mu_i>K_i$ for both $i=1,2$, then we have
$\max\{p_1(K_1)-d_1,p_2(K_2)-d_2\}<-\delta<0$. This implies that
$$\frac{dV}{dt}=\rho_2(p_1(x_1)-d_1)y_1+\rho_1(p_2(x_2)-d_2)y_2<-\delta(\rho_2y_1+\rho_1y_2).$$Therefore,
both predators go extinct if  $\mu_i>K_i$ for both $i=1$ and $i=2$. Since both $\limsup_{t\rightarrow\infty} y_i(t)=0$ for both $i=1,2$. Then we have Model \eqref{2patch2} reduced to the following uncoupled prey model
$$x_i'=r_i x_i\left(1-\frac{x_i}{K_i}\right)$$which converges to $x_i=K_i$. Thus, Model \eqref{2patch2s} has global stability at $(K_1,0,K_2,0)$ when $\mu_i>K_i$ for both $i=1,2$.\\
On the other hand, if $0<\mu_i<K_i$ for both $i=1,2$, then we have
$\max\{p_1(K_1)-d_1,p_2(K_2)-d_2\}>\delta>0$. This implies that
$$\frac{dV}{dt}=\rho_2(p_1(x_1)-d_1)y_1+\rho_1(p_2(x_2)-d_2)y_2>\delta(\rho_2y_1+\rho_1y_2).$$
Therefore, both predators persist if  $0<\mu_i<K_i$ for both $i=1$ and $i=2$.\\\\

Item 7 can be obtained from Item 4 and Item 6.\\

Item 8 can be obtained from eigenvalues of the Jacobian matrix of Model \eqref{2patch2} evaluated at the symmetric interior equilibrium $(\mu,\nu,\mu,\nu)$ through simple algebraic calculations. We omit details. The global stability of $(\mu,\nu,\mu,\nu)$ when $\frac{K-1}{2}<\mu<K$ can be obtained by constructing the following Lyapunov function
{\small$$\begin{aligned}
V(x_1,y_1,x_2,y_2) =\rho_2\int^{x_1}_{\mu} \frac{p_1(\xi_1)-p_1(\mu)}{p_1(\xi_1)}d\xi_1 + \rho_2\int^{y_1}_{\nu} \frac{\eta_1-\nu}{\eta_1}d\eta_1+\rho_1\int^{x_2}_{\mu} \frac{p_2(\xi_2)-p_2(\mu)}{p_2(\xi_2)}d\xi_2+\rho_1\int^{y_2}_{\nu} \frac{\eta_2-\nu}{\eta_2}d\eta_2\end{aligned}$$}which gives
$$\frac{dV(x_1,y_1,x_2,y_2)}{dt} =\rho_2(p_1(\xi_1)-p_1(\mu))(q_1(x_1)-K)+\rho_1(p_2(\xi_2)-p_2(\mu))(q_2(x_2)-K)+\frac{\rho_1\rho_2\nu(y_1-y_2)\left(\frac{1}{y_1}-\frac{1}{y_2}\right)}{y_1y_2}.$$
\end{proof}

%%%%%%%%%%%%%%%%%%%%%%%%%%%%%%%%%%%%%%%%%%%%%%%%%%%%%%%%%%%%%%%
%%%%%%%%%%%%%%%%%%%%%%%%%%%%%%%%%%%%%%%%%%%%%%%%%%%%%%%%%%%%%%%
\section*{Acknowledgement}
%\begin{acknowledgements}
The research of Y.K. and S.K.S. are partially supported by Simons Collaboration Grants for Mathematicians (208902), NSF-DMS(1313312) and Research Fellowship from School of Letters and Sciences. S.K.S. is also partially supported by senior research fellowship from the Council for Scientific and Industrial Research, Government of India. All authors would like to thank Dr. Andrea Bruder for the discussions on the modeling dispersal strategies in the early stage of this manuscript.
%\end{acknowledgements}
%%%%%%%%%%%%%%%%%%%%%%%%%%%%%%%%%%%%%%%%%%%%%%%%%%%%%%%%%%%%%%%
%%%%%%%%%%%%%%%%%%%%%%%%%%%%%%%%%%%%%%%%%%%%%%%%%%%%%%%%%%%%%%%

\section*{References}
\bibliographystyle{abbrvnat}
\bibliographystyle{apalike}
%\bibliography{references}
%\bibliographystyle{acm}	
\bibliography{YK-work27Aug2014}

\begin{thebibliography}{74}
\providecommand{\natexlab}[1]{#1}
\providecommand{\url}[1]{\texttt{#1}}
\expandafter\ifx\csname urlstyle\endcsname\relax
  \providecommand{\doi}[1]{doi: #1}\else
  \providecommand{\doi}{doi: \begingroup \urlstyle{rm}\Url}\fi

\bibitem[Aarssen and Turkington(1985)]{aarssen1985biotic}
L.~Aarssen and R.~Turkington.
\newblock Biotic specialization between neighbouring genotypes in lolium
  perenne and trifolium repens from a permanent pasture.
\newblock \emph{The Journal of Ecology}, pages 605--614, 1985.

\bibitem[Alder(1993)]{Adler1993}
R.~F. Alder.
\newblock Migration alone can produce persistence of host-parasitoid models.
\newblock \emph{The American Naturalist}, 141:\penalty0 642--650, 1993.

\bibitem[Bascompte and Sol{\'e}(1994)]{bascompte1994spatially}
J.~Bascompte and R.~V. Sol{\'e}.
\newblock Spatially induced bifurcations in single-species population dynamics.
\newblock \emph{Journal of Animal Ecology}, pages 256--264, 1994.

\bibitem[Bolker and Pacala(1999)]{bolker1999spatial}
B.~M. Bolker and S.~W. Pacala.
\newblock Spatial moment equations for plant competition: understanding spatial
  strategies and the advantages of short dispersal.
\newblock \emph{The American Naturalist}, 153\penalty0 (6):\penalty0 575--602,
  1999.

\bibitem[Carroll and Janzen(1973)]{carroll1973ecology}
C.~Carroll and D.~H. Janzen.
\newblock Ecology of foraging by ants.
\newblock \emph{Annual Review of Ecology and Systematics}, pages 231--257,
  1973.

\bibitem[Casal et~al.(1994)Casal, Eilbeck, L{\'o}pez-G{\'o}mez,
  et~al.]{casal1994existence}
A.~Casal, J.~Eilbeck, J.~L{\'o}pez-G{\'o}mez, et~al.
\newblock Existence and uniqueness of coexistence states for a predator-prey
  model with diffusion.
\newblock \emph{Differential and Integral Equations}, 7\penalty0 (2):\penalty0
  411--439, 1994.

\bibitem[Chesson and Murdoch(1986)]{chesson1986aggregation}
P.~L. Chesson and W.~W. Murdoch.
\newblock Aggregation of risk: relationships among host-parasitoid models.
\newblock \emph{American Naturalist}, 127\penalty0 (5):\penalty0 696--715,
  1986.

\bibitem[Chewning(1975)]{chewning1975migratory}
W.~C. Chewning.
\newblock Migratory effects in predator-prey models.
\newblock \emph{Mathematical Biosciences}, 23\penalty0 (3):\penalty0 253--262,
  1975.

\bibitem[Cressman and K{\v{r}}ivan(2013)]{cressman2013two}
R.~Cressman and V.~K{\v{r}}ivan.
\newblock Two-patch population models with adaptive dispersal: the effects of
  varying dispersal speeds.
\newblock \emph{Journal of mathematical biology}, 67\penalty0 (2):\penalty0
  329--358, 2013.

\bibitem[Curio(1976)]{Curio1976}
E.~Curio.
\newblock \emph{The ethology of predation}, volume~7.
\newblock Springer-Verlag Berlin Heidelberg, 1976.

\bibitem[Doebli(1995)]{Doebli1995}
M.~Doebli.
\newblock Dispersal and dynamics.
\newblock \emph{Theoretical Population Biology}, 47:\penalty0 82--106, 1995.

\bibitem[Feng et~al.(2011)Feng, Rock, and Hinson]{feng2011new}
W.~Feng, B.~Rock, and J.~Hinson.
\newblock On a new model of two-patch predator-prey system with migration of
  both species.
\newblock \emph{Journal of Applied Analysis and Computation}, 1\penalty0
  (2):\penalty0 193--203, 2011.

\bibitem[Ford(1971)]{Ford1971}
J.~Ford.
\newblock \emph{The Role of the Trypanosomiasis in African Ecology}.
\newblock Oxford University Press. Oxford, 1971.

\bibitem[Gatehouse(1972)]{Gatehouse1972}
A.~G. Gatehouse.
\newblock Permanence and the dynamics of biological systems.
\newblock \emph{Host finding behaviour of tsetse flies. See Ref. 60}, pages
  83--95, 1972.

\bibitem[Ghosh and Bhattacharyya(2011)]{ghosh2011two}
S.~Ghosh and S.~Bhattacharyya.
\newblock A two-patch prey-predator model with food-gathering activity.
\newblock \emph{Journal of Applied Mathematics and Computing}, 37\penalty0
  (1-2):\penalty0 497--521, 2011.

\bibitem[Gillies and Wilkes(1970)]{gillies1970range}
M.~Gillies and T.~Wilkes.
\newblock The range of attraction of single baits for some west african
  mosquitoes.
\newblock \emph{Bulletin of Entomological Research}, 60\penalty0 (02):\penalty0
  225--235, 1970.

\bibitem[Gillies and Wilkes(1972)]{gillies1972range}
M.~Gillies and T.~Wilkes.
\newblock The range of attraction of animal baits and carbon dioxide for
  mosquitoes. studies in a freshwater area of west africa.
\newblock \emph{Bulletin of Entomological Research}, 61\penalty0 (03):\penalty0
  389--404, 1972.

\bibitem[Gillies and Wilkes(1974)]{gillies1974range}
M.~Gillies and T.~Wilkes.
\newblock The range of attraction of birds as baits for some west african
  mosquitoes (diptera, culicidae).
\newblock \emph{Bulletin of Entomological Research}, 63\penalty0 (04):\penalty0
  573--582, 1974.

\bibitem[Hanski and Hanski(1999)]{hanski1999metapopulation}
I.~Hanski and I.~A. Hanski.
\newblock \emph{Metapopulation ecology}, volume 232.
\newblock Oxford University Press Oxford, 1999.

\bibitem[Hanski and (eds)(1997)]{hanski1997metapopulationbiology}
I.~A. Hanski and G.~M.~E. (eds).
\newblock \emph{Metapopulation biology: ecology, genetics, and evolution.}
\newblock Academic Press, San Diego, 1997.

\bibitem[Hassell and May(1974)]{hassell1974aggregation}
M.~Hassell and R.~May.
\newblock Aggregation of predators and insect parasites and its effect on
  stability.
\newblock \emph{The Journal of Animal Ecology}, 43\penalty0 (2):\penalty0
  567--594, 1974.

\bibitem[Hassell and Southwood(1978)]{hassell1978foraging}
M.~Hassell and T.~Southwood.
\newblock Foraging strategies of insects.
\newblock \emph{Annual Review of Ecology and Systematics}, pages 75--98, 1978.

\bibitem[Hassell et~al.(1995)Hassell, Miramontes, Rohani, and
  May]{hassell1995appropriate}
M.~Hassell, O.~Miramontes, P.~Rohani, and R.~May.
\newblock Appropriate formulations for dispersal in spatially structured
  models: comments on bascompte \& sol{\'e}.
\newblock \emph{Journal of Animal Ecology}, pages 662--664, 1995.

\bibitem[Hassell et~al.(1991)Hassell, Comins, and May]{hassell1991spatial}
M.~P. Hassell, H.~N. Comins, and R.~M. May.
\newblock Spatial structure and chaos in insect population dynamics.
\newblock \emph{Nature}, 353\penalty0 (6341):\penalty0 255--258, 1991.

\bibitem[Hastings(1983)]{Hastings1983}
A.~Hastings.
\newblock Can spatial variation along lead to selection for dispersal?
\newblock \emph{Theoretical Population Biology}, 24:\penalty0 244--251, 1983.

\bibitem[Hastings(1995)]{Hastings1993}
A.~Hastings.
\newblock Complex interactions between dispersal and dynamics: \textrm{L}essons
  from coupled logistic equations.
\newblock \emph{Ecology}, 74:\penalty0 1362--1372, 1995.

\bibitem[Hauzy et~al.(2010)Hauzy, Gauduchon, Hulot, and
  Loreau]{hauzy2010density}
C.~Hauzy, M.~Gauduchon, F.~D. Hulot, and M.~Loreau.
\newblock Density-dependent dispersal and relative dispersal affect the
  stability of predator--prey metacommunities.
\newblock \emph{Journal of theoretical biology}, 266\penalty0 (3):\penalty0
  458--469, 2010.

\bibitem[Holt(1985)]{holt1985population}
R.~D. Holt.
\newblock Population dynamics in two-patch environments: some anomalous
  consequences of an optimal habitat distribution.
\newblock \emph{Theoretical population biology}, 28\penalty0 (2):\penalty0
  181--208, 1985.

\bibitem[Hsu et~al.(1977)Hsu, Hubbell, and Waltman]{hsu1977mathematical}
S.~Hsu, S.~Hubbell, and P.~Waltman.
\newblock A mathematical theory for single-nutrient competition in continuous
  cultures of micro-organisms.
\newblock \emph{SIAM Journal on Applied Mathematics}, 32\penalty0 (2):\penalty0
  366--383, 1977.

\bibitem[Hsu(1978)]{hsu1978global}
S.-B. Hsu.
\newblock On global stability of a predator-prey system.
\newblock \emph{Mathematical Biosciences}, 39\penalty0 (1):\penalty0 1--10,
  1978.

\bibitem[Huang and Diekmann(2001)]{huang2001predator}
Y.~Huang and O.~Diekmann.
\newblock Predator migration in response to prey density: what are the
  consequences?
\newblock \emph{Journal of mathematical biology}, 43\penalty0 (6):\penalty0
  561--581, 2001.

\bibitem[Hutson(1984)]{hutson1984theorem}
V.~Hutson.
\newblock A theorem on average liapunov functions.
\newblock \emph{Monatshefte f{\"u}r Mathematik}, 98\penalty0 (4):\penalty0
  267--275, 1984.

\bibitem[Hutson and Schmitt(1992)]{hutson1992permanence}
V.~Hutson and K.~Schmitt.
\newblock Permanence and the dynamics of biological systems.
\newblock \emph{Mathematical Biosciences}, 111\penalty0 (1):\penalty0 1--71,
  1992.

\bibitem[Jansen(1995)]{jansen1995regulation}
V.~A. Jansen.
\newblock Regulation of predator-prey systems through spatial interactions: a
  possible solution to the paradox of enrichment.
\newblock \emph{Oikos}, 74\penalty0 (345):\penalty0 384--390, 1995.

\bibitem[Jansen(2001)]{jansen2001dynamics}
V.~A. Jansen.
\newblock The dynamics of two diffusively coupled predator--prey populations.
\newblock \emph{Theoretical Population Biology}, 59\penalty0 (2):\penalty0
  119--131, 2001.

\bibitem[Jansen(1994)]{jansen1994theoretical}
V.~A.~A. Jansen.
\newblock \emph{Theoretical aspects of metapopulation dynamics}.
\newblock PhD thesis, Ph. D. thesis, Leiden University, The Netherlands, 1994.

\bibitem[Kang and Armbruster(2011)]{Kang2011a}
Y.~Kang and D.~Armbruster.
\newblock Dispersal effects on a discrete two-patch model for plant-insect
  interactions.
\newblock \emph{J Theor Biol}, 268\penalty0 (1):\penalty0 84--97, 1 2011.
\newblock ISSN 1095-8541.
\newblock \doi{10.1016/j.jtbi.2010.09.033}.

\bibitem[Kang and Castillo-Chavez(2012)]{Kang_C2012a}
Y.~Kang and C.~Castillo-Chavez.
\newblock Multiscale analysis of compartment models with dispersal.
\newblock \emph{Journal of Biological Dynamics}, 2012.
\newblock Epub ahead of print.

\bibitem[Kareiva and Odell(1987)]{kareiva1987swarms}
P.~Kareiva and G.~Odell.
\newblock Swarms of predators exhibit" preytaxis" if individual predators use
  area-restricted search.
\newblock \emph{American Naturalist}, 130\penalty0 (2):\penalty0 233--270,
  1987.

\bibitem[Kareiva et~al.(1990)Kareiva, Mullen, and
  Southwood]{kareiva1990population}
P.~Kareiva, A.~Mullen, and R.~Southwood.
\newblock Population dynamics in spatially complex environments: theory and
  data [and discussion].
\newblock \emph{Philosophical Transactions of the Royal Society of London.
  Series B: Biological Sciences}, 330\penalty0 (1257):\penalty0 175--190, 1990.

\bibitem[K{\'e}fi et~al.(2007)K{\'e}fi, Rietkerk, van Baalen, and
  Loreau]{kefi2007local}
S.~K{\'e}fi, M.~Rietkerk, M.~van Baalen, and M.~Loreau.
\newblock Local facilitation, bistability and transitions in arid ecosystems.
\newblock \emph{Theoretical Population Biology}, 71\penalty0 (3):\penalty0
  367--379, 2007.

\bibitem[Klepac et~al.(2007)Klepac, Neubert, and van~den
  Driessche]{klepac2007dispersal}
P.~Klepac, M.~G. Neubert, and P.~van~den Driessche.
\newblock Dispersal delays, predator--prey stability, and the paradox of
  enrichment.
\newblock \emph{Theoretical population biology}, 71\penalty0 (4):\penalty0
  436--444, 2007.

\bibitem[Kummel et~al.(2013)Kummel, Brown, and Bruder]{kummel2013aphids}
M.~Kummel, D.~Brown, and A.~Bruder.
\newblock How the aphids got their spots: predation drives self-organization of
  aphid colonies in a patchy habitat.
\newblock \emph{Oikos}, 122\penalty0 (6):\penalty0 896--906, 2013.

\bibitem[Kuto and Yamada(2004)]{kuto2004multiple}
K.~Kuto and Y.~Yamada.
\newblock Multiple coexistence states for a prey--predator system with
  cross-diffusion.
\newblock \emph{Journal of Differential Equations}, 197\penalty0 (2):\penalty0
  315--348, 2004.

\bibitem[Lengyel and Epstein(1991)]{lengyel1991diffusion}
I.~Lengyel and I.~R. Epstein.
\newblock Diffusion-induced instability in chemically reacting systems:
  Steady-state multiplicity, oscillation, and chaos.
\newblock \emph{Chaos: An Interdisciplinary Journal of Nonlinear Science},
  1\penalty0 (1):\penalty0 69--76, 1991.

\bibitem[Levin(1974)]{levin1974dispersion}
S.~A. Levin.
\newblock Dispersion and population interactions.
\newblock \emph{American Naturalist}, pages 207--228, 1974.

\bibitem[Levins(1969)]{levins1969some}
R.~Levins.
\newblock Some demographic and genetic consequences of environmental
  heterogeneity for biological control.
\newblock \emph{Bulletin of the ESA}, 15\penalty0 (3):\penalty0 237--240, 1969.

\bibitem[Li et~al.(2005)Li, Gao, Hui, Han, and Shi]{li2005impact}
Z.-z. Li, M.~Gao, C.~Hui, X.-z. Han, and H.~Shi.
\newblock Impact of predator pursuit and prey evasion on synchrony and spatial
  patterns in metapopulation.
\newblock \emph{Ecological Modelling}, 185\penalty0 (2):\penalty0 245--254,
  2005.

\bibitem[Liu and Chen(2003)]{liu2003complex}
X.~Liu and L.~Chen.
\newblock Complex dynamics of holling type ii lotka--volterra predator--prey
  system with impulsive perturbations on the predator.
\newblock \emph{Chaos, Solitons \& Fractals}, 16\penalty0 (2):\penalty0
  311--320, 2003.

\bibitem[Liu(2010)]{Liuyuanyuan2010}
Y.~Liu.
\newblock The dynamical behavior of a two patch predator-prey model, 2010.

\bibitem[Loughrin et~al.(1996)Loughrin, Potter, Hamilton-Kemp, and
  Byers]{loughrin1996role}
J.~H. Loughrin, D.~A. Potter, T.~R. Hamilton-Kemp, and M.~E. Byers.
\newblock Role of feeding--induced plant volatiles in aggregative behavior of
  the japanese beetle (coleoptera: Scarabaeidae).
\newblock \emph{Environmental Entomology}, 25\penalty0 (5):\penalty0
  1188--1191, 1996.

\bibitem[Madden(1977)]{madden1977physiological}
J.~Madden.
\newblock Physiological reactions of pinus radiata to attack by woodwasp, sirex
  noctilio f.(hymenoptera: Siricidae).
\newblock \emph{Bulletin of Entomological Research}, 67\penalty0 (03):\penalty0
  405--426, 1977.

\bibitem[Markus(1956)]{Markus1956}
L.~Markus.
\newblock \emph{Asymptotically autonomous differential systems. Contributions
  to the Theory of Nonlinear Oscillations III (S. Lefschetz, ed.), 17-29.
  Annals of Mathematics Studies 36}.
\newblock Princeton University Press, 1956.

\bibitem[May(1978)]{May1978}
R.~M. May.
\newblock Host-parasitoid systems in patchy environments: a phenomenological
  model.
\newblock \emph{Journal of Animal Ecology}, 47:\penalty0 833--843, 1978.

\bibitem[McMurtrie(1978)]{mcmurtrie1978persistence}
R.~McMurtrie.
\newblock Persistence and stability of single-species and prey-predator systems
  in spatially heterogeneous environments.
\newblock \emph{Mathematical Biosciences}, 39\penalty0 (1):\penalty0 11--51,
  1978.

\bibitem[Miller et~al.(2002)Miller, Mladenoff, and Clayton]{miller2002old}
T.~F. Miller, D.~J. Mladenoff, and M.~K. Clayton.
\newblock Old-growth northern hardwood forests: spatial autocorrelation and
  patterns of understory vegetation.
\newblock \emph{Ecological Monographs}, 72\penalty0 (4):\penalty0 487--503,
  2002.

\bibitem[Murdoch et~al.(1992)Murdoch, Briggs, Nisbet, Gurney, and
  Stewart-Oaten]{murdoch1992aggregation}
W.~W. Murdoch, C.~J. Briggs, R.~M. Nisbet, W.~S. Gurney, and A.~Stewart-Oaten.
\newblock Aggregation and stability in metapopulation models.
\newblock \emph{American Naturalist}, pages 41--58, 1992.

\bibitem[Pascual(1993)]{pascual1993diffusion}
M.~Pascual.
\newblock Diffusion-induced chaos in a spatial predator--prey system.
\newblock \emph{Proceedings of the Royal Society of London. Series B:
  Biological Sciences}, 251\penalty0 (1330):\penalty0 1--7, 1993.

\bibitem[Rees et~al.(1996)Rees, Grubb, and Kelly]{rees1996quantifying}
M.~Rees, P.~J. Grubb, and D.~Kelly.
\newblock Quantifying the impact of competition and spatial heterogeneity on
  the structure and dynamics of a four-species guild of winter annuals.
\newblock \emph{American Naturalist}, pages 1--32, 1996.

\bibitem[Rietkerk and Van~de Koppel(2008)]{rietkerk2008regular}
M.~Rietkerk and J.~Van~de Koppel.
\newblock Regular pattern formation in real ecosystems.
\newblock \emph{Trends in Ecology \& Evolution}, 23\penalty0 (3):\penalty0
  169--175, 2008.

\bibitem[Rohani and Ruxton(1999)]{Rohani1999dispersal}
P.~Rohani and G.~D. Ruxton.
\newblock Dispersal and stability in metapopulations.
\newblock \emph{IMA journal of mathematics applied in medicine and biology},
  16:\penalty0 297--306, 1999.

\bibitem[Rosenzweig and MacArthur(1963)]{rosenzweig1963graphical}
M.~L. Rosenzweig and R.~H. MacArthur.
\newblock Graphical representation and stability conditions of predator-prey
  interactions.
\newblock \emph{American Naturalist}, 97\penalty0 (895):\penalty0 209--223,
  1963.

\bibitem[Ruxton(1996)]{ruxton1996density}
G.~D. Ruxton.
\newblock Density-dependent migration and stability in a system of linked
  populations.
\newblock \emph{Bulletin of Mathematical Biology}, 58\penalty0 (4):\penalty0
  643--660, 1996.

\bibitem[Schoonhoven(1974)]{Schoonhoven1974}
L.~M. Schoonhoven.
\newblock Plant recognition by lepidopterous larvae.
\newblock \emph{See ref. 99}, pages 87--99, 1974.

\bibitem[Schoonhoven(1976)]{Schoonhoven1976}
L.~M. Schoonhoven.
\newblock On the variability of chemosensory information.
\newblock \emph{Symp. Biol. Hung.}, 16:\penalty0 261--266, 1976.

\bibitem[Schoonhoven(1977)]{Schoonhoven1977}
L.~M. Schoonhoven.
\newblock Chemosensory systems and feeding behavior in phytophagous insects.
\newblock \emph{See Ref. 19}, pages 391--398, 1977.

\bibitem[Seabloom et~al.(2005)Seabloom, BJ{\o}RNSTAD, Bolker, and
  Reichman]{seabloom2005spatial}
E.~W. Seabloom, O.~N. BJ{\o}RNSTAD, B.~M. Bolker, and O.~Reichman.
\newblock Spatial signature of environmental heterogeneity, dispersal, and
  competition in successional grasslands.
\newblock \emph{Ecological Monographs}, 75\penalty0 (2):\penalty0 199--214,
  2005.

\bibitem[Shahak et~al.(2008)Shahak, Gal, Offir, and
  Ben-Yakir]{shahak2008photoselective}
Y.~Shahak, E.~Gal, Y.~Offir, and D.~Ben-Yakir.
\newblock Photoselective shade netting integrated with greenhouse technologies
  for improved performance of vegetable and ornamental crops.
\newblock In \emph{International Workshop on Greenhouse Environmental Control
  and Crop Production in Semi-Arid Regions 797}, pages 75--80, 2008.

\bibitem[Sol{\'e} and Bascompte(2006)]{sole2006self}
R.~V. Sol{\'e} and J.~Bascompte.
\newblock \emph{Self-Organization in Complex Ecosystems.(MPB-42)}, volume~42.
\newblock Princeton University Press, 2006.

\bibitem[Soro et~al.(1999)Soro, Sundberg, and Rydin]{soro1999species}
A.~Soro, S.~Sundberg, and H.~Rydin.
\newblock Species diversity, niche metrics and species associations in
  harvested and undisturbed bogs.
\newblock \emph{Journal of Vegetation Science}, 10\penalty0 (4):\penalty0
  549--560, 1999.

\bibitem[Thieme(2003)]{thieme2003mathematics}
H.~R. Thieme.
\newblock \emph{Mathematics in population biology}.
\newblock Princeton University Press, 2003.

\bibitem[Tilman and Kareiva(1997)]{tilman1997spatial}
D.~Tilman and P.~M. Kareiva.
\newblock \emph{Spatial ecology: the role of space in population dynamics and
  interspecific interactions}, volume~30.
\newblock Princeton University Press, 1997.

\bibitem[van~de Koppel et~al.(2008)van~de Koppel, Gascoigne, Theraulaz,
  Rietkerk, Mooij, and Herman]{van2008experimental}
J.~van~de Koppel, J.~C. Gascoigne, G.~Theraulaz, M.~Rietkerk, W.~M. Mooij, and
  P.~M. Herman.
\newblock Experimental evidence for spatial self-organization and its emergent
  effects in mussel bed ecosystems.
\newblock \emph{Science}, 322\penalty0 (5902):\penalty0 739--742, 2008.

\bibitem[Waage(1977)]{Waage1977}
J.~K. Waage.
\newblock \emph{Behavioral aspects of foraging in the parasitoid, Nemeritis
  canescens (Grav.)}.
\newblock PhD thesis, PhD thesis. University London, 1977.

\end{thebibliography}
%\bibliography{TwoStrategyTwoPatchBib}

\end{document}